\documentclass[12pt]{article}

\usepackage{amsmath}
\usepackage{amssymb}
\usepackage{booktabs}
\usepackage[dvipdfmx]{color}
\usepackage[dvipdfmx]{graphicx}

\usepackage{bm}
\usepackage{url}
\usepackage{subfigure}
\usepackage{algorithm}
\usepackage{algorithmic}
\usepackage{tabulary}
\usepackage[dvipdfmx]{graphicx}

\newtheorem{theorem}{Theorem}[section]
\newtheorem{prop}{Proposition}[section]
\newtheorem{lemma}{Lemma}[section]
\newtheorem{corollary}{Corollary}[section]

\newtheorem{example}{Example}[section]
\newtheorem{problem}{Problem}[section]

\newtheorem{assum}{Assumption}[section]

\renewcommand{\theenumi}{\roman{enumi}}

\newcommand{\argmin}{\operatornamewithlimits{argmin}}

\numberwithin{equation}{section}

\begin{document}
\makeatletter

\begin{center}
%%%%Performance Optimization of a Dynamic Channel Bonding Strategy in Cognitive Radio Networks
\large{\bf Riemannian Stochastic Fixed Point Optimization Algorithm}\\
\small{This work was supported by JSPS KAKENHI Grant Number JP18K11184.}
\end{center}\vspace{3mm}

\begin{center}
\textsc{Hideaki Iiduka and Hiroyuki Sakai}\\
Department of Computer Science, 
Meiji University,
1-1-1 Higashimita, Tama-ku, Kawasaki-shi, Kanagawa 214-8571 Japan. 
(iiduka@cs.meiji.ac.jp)
\end{center}

\vspace{2mm}

\footnotesize{
\noindent\begin{minipage}{14cm}
{\bf Abstract:}
This paper considers a stochastic optimization problem over the fixed point sets of quasinonexpansive mappings on Riemannian manifolds.
The problem enables us to consider Riemannian hierarchical optimization problems 
over complicated sets, such as the intersection of many closed convex sets,
the set of all minimizers of a nonsmooth convex function,
and the intersection of sublevel sets of nonsmooth convex functions. 
We focus on adaptive learning rate optimization algorithms, which adapt step-sizes (referred to as learning rates in the machine learning field) to find optimal solutions quickly. 
We then propose a Riemannian stochastic fixed point optimization algorithm, which combines fixed point approximation methods on Riemannian manifolds
with the adaptive learning rate optimization algorithms.
We also give convergence analyses of the proposed algorithm
for nonsmooth convex and smooth nonconvex optimization.
The analysis results indicate that, with small constant step-sizes,
the proposed algorithm approximates a solution to the problem.
Consideration of the case in which step-size sequences are diminishing demonstrates 
that the proposed algorithm solves the problem with a guaranteed convergence rate.
This paper also provides numerical comparisons that
demonstrate the effectiveness of the proposed algorithms with formulas 
based on the adaptive learning rate optimization algorithms, such as Adam and AMSGrad. 

\end{minipage}
 \\[5mm]

\noindent{\bf Keywords:} {adaptive learning rate optimization algorithm, fixed point, hierarchical optimization, quasinonexpansive mapping, Riemannian stochastic fixed point optimization algorithm, Riemannian stochastic optimization}\\
\noindent{\bf Mathematics Subject Classification:} {65K05, 90C15, 90C25, 90C26}

\hbox to14cm{\hrulefill}\par

%%%%%%%%%%%%%%%%%%%%%%%%%%%%%%%%%%%%%%%%%%%%%%%%%%%%%%%%%%%%%%%%%%%%%%%
%%%%%%%%%%%%%

\section{Introduction}\label{sec:1}
In light of developments in machine learning and image/signal processing
(see, e.g.,
%\cite{gec2019,hawe2013,liu2017,nickel2017,sakai2020,sel2012,zhang2016}
\cite{gec2019,hawe2013,nickel2017,sakai2020}
and references therein),
Riemannian optimization has attracted a great deal of attention.
Useful iterative algorithms thus have been presented for Riemannian optimization.
For example, 
nonlinear Riemannian conjugate gradient methods 
have been widely studied in 
%\cite{hawe2013,ring2012,sakai_coap,sato2016,sato2015,sel2012} 
\cite{hawe2013,ring2012,sato2016,sato2015}
for unconstrained optimization.
First-order methods 
%\cite{bento2012,liu2017,wang2015_1,zhang2016} 
\cite{bento2012,wang2015_1}
and proximal point algorithms \cite{ferr2002,li2011} have been reported for 
unconstrained/constrained Riemannian optimization.
Riemannian stochastic gradient methods were proposed in \cite{bon2013,kasai2019,sato2019} for Riemannian stochastic optimization.

For training deep neural networks, 
{\em adaptive learning rate optimization algorithms} 
based on using stochastic subgradients and exponential moving averages 
have a strong presence
since they have fast convergence for stochastic optimization 
in Euclidean space.
For example, AdaGrad \cite{adagrad} and RMSProp \cite{tie2012} take advantage of efficient learning rates 
(referred to as step-sizes in the field of optimization)
derived from element-wise squared stochastic gradients. 
Adam \cite{adam} and AMSGrad \cite{reddi2018} are also useful algorithms
using exponential moving averages of stochastic gradients and of 
element-wise squared stochastic gradients.

Recently, Riemannian AMSGrad (RAMSGrad) was studied in \cite{gec2019},
which is a modification of AMSGrad for Euclidean space 
to be applicable to a product of Riemannian manifolds. 
RAMSGrad uses the metric projection onto a constraint convex set
so as to satisfy that the sequence generated by RAMSGrad belongs to the 
constraint set.
Accordingly, RAMSGrad can be applied to only Riemannian convex optimization
with {\em simple} constraints in the sense that the metric projection 
can be easily computed. 

In contrast to \cite{gec2019}, this paper tries to consider 
a Riemannian optimization problem with {\em complicated} constraints, 
such as the intersection of many convex sets \cite{bau,bento2012,wang2015_1}, 
the set of minimizers of a convex function \cite{ferr2002,li2011}, 
and the intersection of sublevel sets of convex functions \cite{wang2015_1}. 
The problem is a {\em hierarchical constrained optimization problem} with three stages, as follows.
The first stage is to find points in Riemannian manifolds (e.g., to find points in nonconvex constraints in Euclidean space). 
The second stage is to find fixed points of quasinonexpansive mappings on Riemannian manifolds.  
Complicated convex sets, such as those mentioned above, can be expressed as the fixed point set of a quasinonexpansive
mapping on a Riemannian manifold (Proposition \ref{mappings}). 
The third stage is to optimize an objective function over the second stage. 
For example, the third stage includes the case of trying to find a stationary point
of a smooth nonconvex function over the set of minimizers of a convex
function over a Riemannian manifold.

The reason why the above problem should be considered is to enable us to resolve
unsolved optimization problems on Riemannian manifolds.  
For example, in the natural language processing for hierarchical representations of symbolic data, 
embeddings into a Poincar\'e ball perform better than embeddings into
a Euclidean space \cite{nickel2017}.
This implies that a Riemannian optimization problem should be considered 
for natural language processing.
As seen above, we expect to gain new insights from
re-considering 
several problems with complicated constraints in the Hilbert/Euclidean space setting
as Riemannian optimization. 
In addition, a classifier ensemble problem with sparsity leaning can be expressed as a Euclidean convex optimization problem over the sublevel set of a convex function \cite{iiduka_cyb}.
Since the results in this paper enable us to consider
a Riemannian optimization problem over the sublevel set of a convex function, there is a possibility that the results will lead to new Riemannian learning methods which can outperform the existing methods in \cite{iiduka_cyb} by a wide margin. 

%In the Hilbert/Euclidean space setting, 
%there have been useful methods, such as incremental methods \cite{neto2009,nedic2001}, (parallel) proximal methods \cite{comb2007,com2008,pes2012}, and 
%fixed point methods \cite{com1999,yamada},
%for convex optimization with complicated constraints.
%This paper studies a Riemannian optimization problem with complicated constraints
%by using the ideas of fixed point methods, which have been applied to significant 
%applications in the Hilbert/Euclidean space setting,
%such as signal processing \cite{com1999}, network resource allocation \cite{iiduka_siopt2013,iiduka_TCNS}, and ensemble learning \cite{iiduka_cyb}.

We first define {\em quasinonexpansive mappings} of which 
{\em fixed point sets} are equal to complicated constraint sets.
Thanks to the previously reported results in \cite{ferr2002,li2010,li2011},
we can define quasinonexpansive mappings for the cases where 
the constraint sets are those mentioned above: 
the intersection of many closed convex sets, 
the set of all minimizers of a nonsmooth convex function,
and the intersection of sublevel sets of convex functions  
(Proposition \ref{examples}). 
Accordingly, the Riemannian optimization problem with such complicated constraints can be expressed as a Riemannian optimization problem 
over the fixed point sets of quasinonexpansive mappings (Problem \ref{problem:1}).
Next, we combine the ideas of adaptive learning rate optimization algorithms (see the second and third paragraphs of this section)
with the fixed point methods \cite{li2010}. 
We then propose a {\em Riemannian stochastic fixed point optimization algorithm} (Algorithm \ref{algo:1}) for solving the problem. 

The intellectual contribution of this paper is that the proposed methodology enables one to deal with
Riemannian optimization over the fixed point sets of quasinonexpansive mappings, especially in contrast to recent papers 
\cite{gec2019,sakai2020} that discussed Riemannian convex optimization over simple constraints.
To clarify this contribution, let us consider the case where 
a constraint set is the intersection of many closed convex sets 
on a Riemannian manifold (Proposition \ref{examples}(ii), Example \ref{exp:1}).
Even if the metric projection onto each closed convex set can be easily
computed within a finite number of arithmetic operations, 
the metric projection onto the intersection of many closed convex sets
would not be implemented in practice.  
This is because, for each iteration, we must solve a subproblem
of minimizing the distance function over the intersection to find the nearest point to the intersection. 
Meanwhile, we can use a {\em computable} mapping which consists of 
the product of the metric projections (see Example \ref{exp:1} for the details), since the metric projection onto each closed convex set can be easily implemented.
This computable mapping satisfies the quasinonexpansivity condition 
and that the fixed point set of this mapping coincides with 
the intersection of many closed convex sets. 
Therefore, the proposed algorithm (Algorithm \ref{algo:1})
using this computable quasinonexpansive mapping can be applied to
Riemannian optimization over the intersection of many closed convex sets, in contrast to \cite{gec2019,sakai2020},
which discussed Riemannian convex optimization over simple constraints.
This paper also gives other examples of Problem \ref{problem:1},
namely, Riemannian optimization over the set of minimizers of 
a convex function (Example \ref{exp:2}) and 
Riemannian optimization over the intersection of sublevel sets of 
convex functions (Example \ref{exp:3}).

% nonconvex, smooth case mo
The theoretical contribution of this paper is its analysis of the proposed algorithm (Algorithm \ref{algo:1}) for solving the Riemannian optimization problem
over the fixed point sets of quasinonexpansive mappings (Problem \ref{problem:1}).
The analysis indicates that the proposed algorithm with small constant step-sizes can approximate a solution to the main problem (Theorems \ref{cor:1} and \ref{cor:2_1}). 
The analysis also shows that the proposed algorithm with diminishing step-sizes
can solve the main problem with a guaranteed convergence rate (Theorems \ref{cor:2} and \ref{cor:2_2}, and Corollary \ref{COR:1_1}). 
%Optimization algorithms with diminishing step-sizes would not be 
%implementable in practice because 
%the step-sizes are approximately zero for a number of iterations.  
%Therefore, 
%the advantage of the proposed algorithm is that it uses constant step-sizes rather than diminishing step-sizes, in contrast to the convergence analysis of RAMSGrad \cite[Theorem 1]{gec2019} with only diminishing step-sizes.
%We also discuss comparisons of the proposed algorithm with RAMSGrad (Remark \ref{rem:1}).

The practical contribution of this paper is a presentation of numerical results demonstrating that the proposed algorithm can be applied to Riemannian optimization over fixed point constraints. 
In this paper, we consider two cases for the constraint conditions. 
The first case is a consistent case such that the intersection of finite closed balls on the Poincar\'e disk is nonempty (Subsection \ref{consistent}). 
The second case is an inconsistent case such that the intersection is empty
(Subsection \ref{inconsistent}).
For the second case, we define a {\em generalized convex feasible set}
as a subset of the absolute constrained set with the elements closest to 
the subsidiary constraint set.
Numerical results show that the proposed algorithms with formulas 
based on Adam and AMSGrad perform well. 

%This paper is organized as follows. 
%Section \ref{sec:2} gives the mathematical preliminaries.
%Section \ref{sec:3} defines the Riemannian convex stochastic optimization
%problem over the fixed point sets of quasinonexpansive mappings (Problem \ref{problem:1}) and its examples (Examples \ref{exp:1}, \ref{exp:2}, and \ref{exp:3}).
%Section \ref{sec:4} presents the Riemannian stochastic fixed point optimization algorithm for solving the main problem.
%Section \ref{sec:5} analyzes the proposed algorithm's convergence. 
%Section \ref{sec:6} numerically compares the behaviors
%of the proposed algorithm with .... 
%Section \ref{sec:7} concludes the paper with a brief summary. 

\section{Mathematical Preliminaries}\label{sec:2}
Let $\mathbb{N}$ be the set of all positive integers including zero,
$\mathbb{R}^I$ be an $I$-dimensional Euclidean space, 
$\mathbb{R}_+^I := \{ (x_i)_{i=1}^I \in \mathbb{R}^I \colon x_i \geq 0 \text{ } (i=1,2,\ldots,I)  \}$,
and $\mathbb{R}_{++}^I := \{ (x_i)_{i=1}^I \in \mathbb{R}^I \colon x_i > 0 \text{ } (i=1,2,\ldots,I)  \}$.
%Let $\mathbb{S}^I$ be the set of $I \times I$ symmetric matrices, i.e., 
%$\mathbb{S}^I = \{ X \in \mathbb{R}^{I \times I} \colon X=X^\top   \}$.
%Let $\mathbb{S}^I_{++}$ denote the set of symmetric positive-definite matrices, i.e., 
%$\mathbb{S}^I_{++} = \{ X \in \mathbb{S}^{I} \colon X \succ O   \}$.
%Let $\mathsf{diag}(x_i)$ be an $I \times I$ diagonal matrix with diagonal components  
%$x_i \in \mathbb{R}$ $(i=1,2,\ldots,I)$ and 
%let $\mathbb{D}^I$ be the set of $I \times I$ diagonal matrices, i.e., 
%$\mathbb{D}^I = \{ X \in \mathbb{R}^{I \times I} \colon X = \mathsf{diag}(x_i), \text{ } x_i \in \mathbb{R} \text{ } (i=1,2,\ldots,I)    \}$.
Let $\mathbb{E}[X]$ denote the expectation of random variable $X$.
%The history of the process $\xi_0,\xi_1,\ldots$ up to time $n$ is denoted by 
%$\xi_{[n]} = (\xi_0,\xi_1,\ldots,\xi_n)$.
%Let $\mathbb{E}[X|\xi_{[n]}]$ denote the conditional expectation of $X$ given 
%$\xi_{[n]} = (\xi_0,\xi_1,\ldots,\xi_n)$.
Unless stated otherwise, all relations between random variables are supported to 
hold almost surely.

\subsection{Riemannian manifold and Hadamard manifold}\label{subsec:2.1}
Let $M$ be a connected $m$-dimensional smooth manifold. 
Let $T_x M$ be the tangent space of $M$ at $x \in M$ and $TM = \bigcup_{x\in M} T_x M$ be the tangent bundle 
of $M$. 
A Riemannian metric at $x \in M$ is denoted by 
$\langle \cdot, \cdot \rangle_x \colon T_x M \times T_x M \to \mathbb{R}$
and its induced norm is defined for all $u \in T_x M$ by $\| u \|_x := \sqrt{\langle u, u \rangle_x}$.
Manifold $M$ endowed with Riemannian metric $\langle \cdot,\cdot \rangle := (\langle \cdot, \cdot \rangle_x)_{x\in M}$ is called a Riemannian manifold.

Given a piecewise smooth curve $\gamma \colon [a,b] \to M$
joining $p$ to $q$ (i.e., $\gamma(a)=p$ and $\gamma(b)=q$),
the length $L(\gamma)$ of $\gamma$ is defined by 
$L(\gamma) := \int_a^b \| \dot{\gamma}(t)  \|_{\gamma(t)} \mathrm{d}t$,
where $\dot{\gamma}$ denotes the derivative of $\gamma$.
The distance function $\mathrm{d} \colon M \times M \to \mathbb{R}_+$
is defined for all $p,q \in M$ by the minimal length over the set of 
all such curves joining $p$ to $q$.

A complete, simply connected Riemannian manifold of nonpositive sectional curvature is called an Hadamard manifold. 
An $m$-dimensional Hadamard manifold $M$ is diffeomorphic to the Euclidean space $\mathbb{R}^m$ \cite[Chapter V, Corollary 3.5]{sakai1996}.
An exponential mapping at a point $x$ in an Hadamard manifold $M$ is denoted by $\exp_x \colon T_x M \to M$. 
The mapping $\exp_x$ is well-defined on $T_x M$, which is guaranteed by the Hopf-Rinow theorem \cite[Chapter III, Theorem 1.1]{sakai1996}.
The mapping $\exp_x$ maps $u \in T_x M$ to $y := \exp_x (u) \in M$ 
such that there exists a geodesic $\gamma \colon [a,b] \to M$
satisfying $\gamma(a) =x$, $\gamma(b)= y$, and $\dot{\gamma}(a) = u$.
The Hadamard-Cartan theorem \cite[Chapter V, Theorem 4.1]{sakai1996}
guarantees that $\exp_x$ is diffeomorphic, 
that is, there exists an inverse mapping $\exp_x^{-1} \colon M \to T_x M$.
For all $x,y\in M$, $\varphi_{x \to y}$ denotes an isometry from 
$T_x M$ to $T_y M$.

Let $M^i$ be an $m^i$-dimensional Hadamard space and 
$M$ be the Cartesian product of the $M^i$s, i.e., $M := M^1 \times M^2 \times \cdots \times M^I$.
The tangent space of $M$ at $x = (x^1, x^2, \ldots,x^I) \in M$ is defined by $T_x M := T_{x^1} M^1 \oplus T_{x^2} M^2 \oplus \cdots \oplus T_{x^I} M^I$, where $\oplus$ stands for the direct sum of vector spaces.
For all $x = (x_i)_{i=1}^I \in M$, 
we define $\psi \in T_x M$ by $\psi = (\psi^i)_{i=1}^I = (\psi^1, \psi^2, \ldots, \psi^I)$,
where $\psi^i \in T_{x^i} M^i$.
An exponential mapping at a point $x^i \in M^i$ is denoted by 
$\exp_{x^i}^i$, and an isometry from $T_{x^i} M^i$ to $T_{y^i} M^i$ 
is denoted by $\varphi^i_{x^i \to y^i}$.

\subsection{Convexity, monotonicity, and related mappings}\label{subsec:2.2}
Let $M$ be an Hadamard manifold.
A set $C \subset M$ is referred to as a convex set  
(see, e.g., \cite[Subsection 3.1]{ferr2002} and references therein) if, for any pair of points in $C$, the geodesic joining those two points is contained in $C$.
Suppose that $C \subset M$ is nonempty, closed, and convex,
and $x\in M$.
Then there exists a unique point 
\cite[Corollary 3.1]{ferr2002}, denoted by $P_C(x)$, 
such that 
\begin{align*}
P_C(x) \in C \text{ and }
\mathrm{d}(x,P_C(x))= \inf_{y\in C} \mathrm{d}(x,y) =: \mathrm{d}(x,C).
\end{align*}
We call $P_C$ the {\em metric projection} onto $C$.

A function $f \colon M \to \mathbb{R}$ is said to be convex (see, e.g., \cite[Subsection 3.2]{ferr2002} and references therein)
if, for any geodesic $\gamma$ of $M$, 
$f \circ \gamma \colon \mathbb{R} \to \mathbb{R}$ is convex.
Accordingly, any convex function on $M$ is continuous.
Suppose that $f \colon M \to \mathbb{R}$ is convex.
Theorem 3.3 in \cite{ferr2002} guarantees that, for all $x\in M$, there exists $u_x \in T_x M$ such that,
for all $y\in M$,
\begin{align*}
f(y) \geq f(x) + \langle u_x, \exp_x^{-1} (y) \rangle_x.
\end{align*} 
The tangent vector $u_x$ is called a {\em subgradient} of $f$ at $x$.
When $f$ is smooth, the vector $u_x$ is called the Riemannian gradient of $f$ at $x$
and is denoted by $\mathrm{grad}f(x)$.
The {\em subdifferential vector field} $\partial f \colon M \rightrightarrows TM$ of a convex function 
$f \colon M \to \mathbb{R}$ is defined by
the set of all subgradients of $f$, i.e., for all $x\in M$,
\begin{align*}
\partial f(x) := 
\left\{ u \in T_x M \colon f(y) \geq f(x) + \langle u, \exp_x^{-1} (y) \rangle_x \text{ } (y\in M)
\right\} \neq \emptyset.
\end{align*}
The {\em subgradient projection} $P_{f, \lambda}$ relative to a convex function $f \colon M \to \mathbb{R}$ and $\lambda > 0$ 
is defined for all $x\in M$ by
\begin{align*}
P_{f,\lambda}(x) :=
\begin{cases}
x \qquad &\left(x \in \mathrm{lev}_{\leq 0} (f) := 
\left\{x \in M \colon f(x) \leq 0 \right\} \right),\\
\displaystyle{\exp_x \left( - \lambda \frac{f(x)}{\|u_x\|_x} u_x \right)}
\qquad &(x \notin \mathrm{lev}_{\leq 0} (f)),
\end{cases}
\end{align*}
where $u_x$ is any tangent vector in $\partial f (x)$.
The results in \cite[Lemma 3.1]{b-c2014}, \cite[Proposition 2.3]{b-c2001}, 
and \cite[Subchapter 4.3]{vasin} provide the definition and properties of the subgradient projection under the Hilbert space setting.
 
Let $A \colon M \rightrightarrows TM$ be a set-valued vector field 
such that, for all $x\in D(A) := \{x \in M \colon A(x) \neq \emptyset \}$,
$A(x) \subset T_x M$. 
$A$ is said to be {\em monotone} (see, e.g., \cite[Definition 2]{li2011}
and references therein) if, for all $x,y \in D(A)$,
all $u \in A(x)$, and all $v\in A(y)$,
%\begin{align*}
$\langle u, \exp_x^{-1} (y) \rangle_x \leq \langle v, - \exp_y^{-1} (x) \rangle_y$.
%\end{align*}
$A$ is said to be {\em maximal} (see, e.g., \cite[Definition 2]{li2011}
and references therein) if $A$ is monotone and the following holds: 
for all $x\in M$ and all $u\in T_x M$,
$\langle u, \exp_x^{-1} (y) \rangle_x \leq \langle v, - \exp_y^{-1} (x) \rangle_y$ ($y\in D(A), v\in A(y)$) 
implies that $u \in A(x)$.
The subdifferential vector field $\partial f$ of a convex function $f \colon M \to \mathbb{R}$ with 
$D(f) := \{ x \in M \colon f(x) < + \infty  \}= M$
is maximal monotone \cite[Theorem 5.1]{li2009}.
We call the set of zeros of a set-valued vector field $A\colon M \rightrightarrows TM$ 
the zero point set, which is defined by  
$\mathrm{zer}(A) := \{ x\in D(A) \colon 0 \in A(x) \}$.

Let $\lambda > 0$.
The {\em resolvent} $J_\lambda \colon M \rightrightarrows M$ 
\cite[Definition 6]{li2011} of a set-valued vector field $A \colon M \rightrightarrows TM$ 
is defined for all $x \in M$ by 
\begin{align*} 
J_\lambda (x) := \left\{ z\in M \colon x \in \exp_z \left(\lambda A(z) \right)\right\}.
\end{align*} 
$J_\lambda$ is single-valued when $A$ is monotone \cite[Theorem 4]{li2011}.
The {\em Moreau-Yosida regularization} 
$R_\lambda^f \colon M \rightrightarrows M$ \cite[(20)]{ferr2002}, \cite[(60)]{li2011}
of a convex function $f \colon M \to \mathbb{R}$ is defined 
for all $x\in M$ by 
\begin{align*}
R_\lambda^f (x) := \argmin_{y\in M} 
\left\{ f(y) + \frac{1}{2 \lambda} \mathrm{d}(x,y)^2 \right\}.
\end{align*} 
$R_\lambda^f$ is single-valued and $D(R_\lambda^f) = M$ \cite[Lemma 4.2]{ferr2002}.

\subsection{Nonexpansivity and fixed point set}\label{subsec:2.3}
Let $C$ be a nonempty subset of a Riemannian manifold $M$ with the distance function $\mathrm{d}$
and let $T \colon C \to M$ be a mapping.
The {\em fixed point set} of $T$ is defined by 
\begin{align*}
\mathrm{Fix}(T) := \left\{ x \in C \colon T(x) = x \right\}.
\end{align*}
$T$ is said to be {\em firmly nonexpansive} \cite[Definition 1]{li2011}, \cite[Subchapter 1.11]{goebel2} if, for all $x,y\in C$,
the function $\Phi \colon [0,1] \to \mathbb{R}_+$ defined by 
\begin{align}\label{firm_non}
\Phi(t) := \mathrm{d}\left( \exp_x [t \exp_x^{-1} (T(x))], 
\exp_y [t \exp_y^{-1} (T(y))] \right)
\text{ is decreasing.}
\end{align} 
$T$ is said to be {\em nonexpansive} if
\begin{align}\label{non}
\mathrm{d}(T(x),T(y)) \leq \mathrm{d}(x,y)
\quad (x,y \in C).
\end{align}
$T$ is said to be {\em quasinonexpansive} if
\begin{align}\label{quasi}
\mathrm{d}(T(x),y) \leq \mathrm{d}(x,y)
\quad (x \in C, y\in \mathrm{Fix}(T)).
\end{align}
$T$ is said to be {\em strictly quasinonexpansive} if
\begin{align}\label{squasi}
\mathrm{d}(T(x),y) <  \mathrm{d}(x,y)
\quad (x \in C \backslash \mathrm{Fix}(T), y\in \mathrm{Fix}(T)).
\end{align}
Finally, $T$ is said to be {\em firmly quasinonexpansive} if
\begin{align}\label{fquasi}
\mathrm{d}(T(x),y)^2 + \mathrm{d}(T(x),x)^2 \leq \mathrm{d}(x,y)^2
\quad (x \in C, y\in \mathrm{Fix}(T)).
\end{align}

The following proposition is true. 

\begin{prop}\label{fixed}
Suppose that $C$ is a nonempty, closed convex set of an Hadamard manifold $M$ 
and $T \colon C \to C$ is quasinonexpansive.
\begin{enumerate}
\item[{\em (i)}]{\em \cite[Theorem 1.3]{chao2006}}
If $\mathrm{Fix}(T)$ is nonempty, then $\mathrm{Fix}(T)$ is closed and convex;
\item[{\em (ii)}]{\em \cite[Theorem 13]{kirk2003_1}}
If $C$ is bounded and $T$ is nonexpansive, then $\mathrm{Fix}(T)$ is nonempty.
\end{enumerate} 
\end{prop}

The relationships between the above mappings are given in the following 
proposition (the proof is given in Supplementary Material).

\begin{prop}\label{mappings}
Let $C$ be a nonempty subset of an $m$-dimensional Hadamard manifold $M$
with the distance function $\mathrm{d}$
and let $T \colon C \to M$ be a mapping.
Then, 
\begin{enumerate}
\item[{\em (i)}]
\eqref{firm_non} implies \eqref{non}, and \eqref{non} implies 
\eqref{quasi};
\item[{\em (ii)}]
\eqref{firm_non} implies \eqref{fquasi},
\eqref{fquasi} implies \eqref{squasi}, and 
\eqref{squasi} implies \eqref{quasi}.
\end{enumerate}
\end{prop}

Let $T \colon C \to M$ be quasinonexpansive and $\alpha \in (0,1)$.
Here, we define $S_\alpha \colon C \to M$ as follows: 
for all $x\in M$,
\begin{align*}
S_\alpha (x) := \exp_x [(1 - \alpha) \exp_x^{-1} (T(x))].
\end{align*}
Then, the condition $\mathrm{Fix}(T) = \mathrm{Fix}(S_\alpha)$
holds from the facts that $\exp_x$ is bijective and $\exp_x (0_x) = x$,
where $0_x$ denotes the zero element of $T_x M$.
Moreover, the discussion in \cite[p.553]{li2010} guarantees that,
for all $x\in C \backslash \mathrm{Fix}(T)$ and all $y\in \mathrm{Fix}(T)$, 
\begin{align}\label{keyineq}
\mathrm{d}(S_\alpha (x), y)^2
\leq 
\mathrm{d}(x, y)^2 - \alpha (1-\alpha) \mathrm{d}(T(x), x)^2
< 
\mathrm{d}(x, y)^2,
\end{align} 
that is, $S_\alpha$ is strictly quasinonexpansive.

The following proposition suggests some examples of quasinonexpansive mappings (the proof is given in Supplementary Material).

\begin{prop}\label{examples}
Let $M$ be an $m$-dimensional Hadamard manifold, $C_j$ ($j=1,\ldots,J$) be a nonempty, closed convex subset of $M$,
and $\lambda > 0$.
Suppose that $P_j := P_{C_j}$ is the metric projection onto $C_j$,
$A \colon M \rightrightarrows TM$ is monotone, 
$g \colon M \to \mathbb{R}$ is convex, and 
$P_{g,\lambda}$ is the subgradient projection.
Then, the following hold:
\begin{enumerate}
\item[{\em (i)}]
The metric projection $P_j$ is firmly nonexpansive with $\mathrm{Fix}(P_j) = C_j$;
\item[{\em (ii)}]
Under $\bigcap_{j=1}^J C_j \neq \emptyset$, 
the mapping $T := P_1 P_2 \cdots P_J$ is nonexpansive
with $\mathrm{Fix}(T) = \bigcap_{i=1}^J C_i$;
\item[{\em (iii)}]
The resolvent $J_\lambda$ of $A$ is firmly nonexpansive with $\mathrm{Fix}(J_\lambda) = \mathrm{zer}(A)$;
\item[{\em (iv)}]
The Moreau-Yosida regularization $R_\lambda^g$ of $g$ is firmly nonexpansive with $\mathrm{Fix}(R_\lambda^g) = \argmin_{x\in M} g(x)$;
\item[{\em (v)}]
The subgradient projection $P_{g,\lambda}$ satisfies that
$\mathrm{Fix}(P_{g,\lambda}) = \mathrm{lev}_{\leq 0} (g)$.
\end{enumerate}
Moreover, suppose that $M$ has its sectional curvature lower-bounded by $\kappa \leq 0$,
that $C$ has a diameter bounded by $D$, and that $h_j \colon C \to \mathbb{R}$
($j=1,2,\ldots,J$) is convex with $D(h_j) = C$. 
Then, the following also hold:
\begin{enumerate}
\item[{\em (vi)}]
The subgradient projection $P_{h_j,\lambda}$ with $\lambda \in (0,2/\zeta)$
is strictly quasinonexpansive with
$\mathrm{Fix}(P_{h_j,\lambda}) = \mathrm{lev}_{\leq 0} (h_j)$,
where $\zeta$ is a positive number depending on $\kappa$ and $D$;
\item[{\em (vii)}]
Under $\bigcap_{j=1}^J \mathrm{lev}_{\leq 0} (h_j) \neq \emptyset$, 
the mapping $T := P_{h_1,\lambda}P_{h_2,\lambda} \cdots P_{h_J,\lambda}$ with 
$\lambda \in (0,2/\zeta)$
is strictly quasinonexpansive with 
$\mathrm{Fix}(T) = \bigcap_{j=1}^J \mathrm{lev}_{\leq 0} (h_j)$.
\end{enumerate}
\end{prop}

\section{Stochastic Optimization over Fixed Point Set on 
Riemannian Manifold}
\label{sec:3}
This paper considers the following problem.

\begin{problem}\label{problem:1}
Let $M^i$ ($i\in \mathcal{I} := \{1,2,\ldots,I\}$) be an $m^i$-dimensional Hadamard manifold with sectional curvature lower-bounded by $\kappa^i \leq 0$
and distance function $\mathrm{d}^i$ 
and $M$ be the Cartesian product of the $M^i$s, i.e., $M := M^1 \times M^2 \times \cdots \times M^I$.
Assume that
\begin{enumerate}
\item[(A1)]
$T^i \colon M^i \to M^i$ ($i\in\mathcal{I}$) is quasinonexpansive
with  
$\mathrm{Fix}(T^i) \neq \emptyset$ ($i\in\mathcal{I}$), and 
$X := \mathrm{Fix}(T^1) \times \mathrm{Fix}(T^2) \times \cdots \times \mathrm{Fix}(T^I)$;
\item[(A2)]
A function $f \colon M \to \mathbb{R}$ is defined for all $x \in M$ by 
$f(x) := \mathbb{E}[F({x},\xi)]$, where $F(\cdot,\xi) \colon M \to \mathbb{R}$
and
$\xi$ is a random vector whose probability distribution $P$ is supported on a set $\Xi \subset \mathbb{R}^M$. 
\end{enumerate}
Then, we would like to find a point $x_\star$ in $X_\star$ 
defined by 
\begin{align*}
X_\star := 
\left\{ x_\star \in X \colon 
\left\langle \exp_{x_\star}^{-1} (x), \mathsf{g} (x_\star) \right\rangle_{x_\star} \geq 0 \text{ } (x\in X) \right\},
\end{align*}
where $\mathsf{g}(x) = (\mathsf{g}^i(x))_{i\in\mathcal{I}}$ denotes the (sub)gradient of $f$.
\end{problem}

The relationship between Problem \ref{problem:1} and the problem of minimizing $f$ 
over $X$ is expressed by the following proposition.

\begin{prop}\label{vi}
Suppose that Assumptions (A1) and (A2) hold. 
\begin{enumerate}
\item[{\em (i)}]
If $f \colon M \to \mathbb{R}$ is smooth, 
%with the Riemannian gradient $\mathrm{grad}f$,  
then
\begin{align*}
X_\star \supset \argmin_{x\in X} f(x) := \left\{ x_\star \in X \colon f(x_\star)
= f_\star := \inf_{x\in X} f(x)  \right\};
\end{align*}
\item[{\em (ii)}]
If $f \colon M \to \mathbb{R}$ is convex, then
\begin{align*}
X_\star = \argmin_{x\in X} f(x).
\end{align*}
\end{enumerate}  
\end{prop}

Proposition \ref{fixed}(i) %(the closedness and convexity of $\mathrm{Fix}(T^i)$) 
and Proposition 3.1 in \cite{li2009_1} imply Proposition \ref{vi}(i), which in turn implies that Problem \ref{problem:1} when $f$ is smooth and nonconvex is a stationary point problem 
associated with the nonconvex optimization problem 
to minimize $f$ over $X$. 
Meanwhile, when $f$ is nonsmooth and convex, from the definition of the subdifferential vector field $\partial f$, 
we can prove Proposition \ref{vi}(ii), i.e., that Problem \ref{problem:1} coincides with the nonconvex optimization problem to minimize $f$ over $X$.
%Moreover, $X_\star$ in Problem \ref{problem:1} can be expressed as follows.
%
%\begin{prop}{\em \cite[Theorem 3.3]{li2009_1}}\label{vi1}
%Suppose that (A1) and (A2) hold and $f$ is smooth. 
%Then, there exists $\bar{r} > 0$ such that, for 
%all $r \in (0,\bar{r})$, 
%\begin{align*}
%X_\star = \mathrm{Fix}\left(P_X (\exp (- r \mathrm{grad} f)) \right),
%\end{align*}
%where $\exp (- r \mathrm{grad} f)(x) := \exp_x (- r \mathrm{grad} f(x))$ ($x\in M$).
%\end{prop}

Proposition \ref{examples}(i) and (ii) suggest the following example of Problem \ref{problem:1}.

\begin{example}[Optimization over the intersection of convex sets]
\label{exp:1}
Let $C_j^i$ ($i\in \mathcal{I}, j\in \mathcal{J}^i := \{1,2,\ldots,J^i\}$) be a nonempty, closed convex subset of 
$M^i$ with $\bigcap_{j\in \mathcal{J}^i} C_j^i \neq \emptyset$
and  
$P_j^i$ ($j\in \mathcal{J}^i$) be the metric projection onto $C_j^i$.
Then, Problem \ref{problem:1} with a mapping $T^i:= P_1^i P_2^i \cdots P_{J^i}^i$ ($i\in \mathcal{I}$) is to find a point $x_\star$ in $X_\star$ with 
\begin{align*}
X = \bigcap_{j\in \mathcal{J}^1} C_j^1 \times \bigcap_{j\in \mathcal{J}^2} C_j^2 \times \cdots \times \bigcap_{j\in \mathcal{J}^I} C_j^I.
\end{align*}
\end{example}

Let us compare the convex optimization problem considered in \cite[Section 4]{gec2019}
with Example \ref{exp:1}.
Example \ref{exp:1} when $J^i = 1$ ($i\in \mathcal{I}$) and $f$ is a convex function coincides with 
the problem in \cite[Section 4]{gec2019} that is to 
\begin{align}\label{Rie}
\text{minimize } f(x) \text{ subject to } 
x \in C^1 \times C^2 \times \cdots \times C^I,
\end{align}
where $C^i := C_1^i$ ($i\in \mathcal{I}$) is simple in the sense that 
$P^i := P_1^i$ can be easily computed.
Meanwhile, Example \ref{exp:1} has three stages as follows:
The first stage is to find points of $M$.
The second stage is to find points of complicated sets
$\bigcap_{j\in \mathcal{J}^i} C_j^i$ ($i\in\mathcal{I}$),
which are each the intersection of many convex sets.
The problem in the second stage is called a {\em convex feasibility problem} \cite{bau}, \cite[p.99]{b-c}, \cite{bento2012,wang2015_1}.
The third stage is to minimize a function over the second stage.
Hence, Problem \ref{problem:1} includes optimization problems with complicated constraints,
as seen in Example \ref{exp:1}. 

From Proposition \ref{examples}(iii) and (iv),
we also have the following.

\begin{example}[Optimization over the zero point sets]\label{exp:2}
Let $A^i \colon M^i \rightrightarrows TM^i$ ($i\in \mathcal{I}$) be a monotone set-valued vector field with 
$\mathrm{zer}(A^i) \neq \emptyset$
and $J^i_{\lambda^i} \colon M^i \to M^i$ ($i\in \mathcal{I}$) be the resolvent of $A^i$ with $\lambda^i > 0$.
Then, Problem \ref{problem:1} with a mapping $T^i:= J^i_{\lambda^i}$ ($i\in \mathcal{I}$) is to find a point $x_\star$ in $X_\star$ with
\begin{align*}
X = \mathrm{zer} \left(A^1 \right) \times \mathrm{zer} \left(A^2 \right)\times \cdots \times \mathrm{zer} \left(A^I \right).
\end{align*}
In the case where $A^i := \partial g^i$ ($i\in \mathcal{I}$), where $g^i \colon M \to \mathbb{R}$ is convex, 
the problem is to find a point in $X_\star$ with 
\begin{align*}
X = \argmin_{x^1\in M^1} g^1(x^1) \times \argmin_{x^2\in M^2} g^2(x^2) \times \cdots \times \argmin_{x^I\in M^I} g^I(x^I).
\end{align*}
\end{example}

References \cite{ferr2002} and \cite{li2011} presented proximal point algorithms
which use the resolvents of a monotone vector field $A$
for finding a zero of $A$, 
\begin{align*}
x^* \in \mathrm{zer} (A).
\end{align*}
Thanks to the results in \cite{ferr2002} and \cite{li2011} for the resolvents and Moreau-Yosida regularizations,
Problem \ref{problem:1} includes Example \ref{exp:2} that is to minimize not only convex functions 
$g^i$ (using the resolvents of $\partial g^i$) but also a function $f$ over 
the sets of minimizers of the $g^i$s.
%Meanwhile, references \cite{liu2017} and \cite{zhang2016} presented first-order methods for constrained optimization on Riemannian manifolds.

Proposition \ref{examples} (v)--(vii) suggest the following example:

\begin{example}[Optimization over the sublevel sets of convex functions]\label{exp:3}
Let $C^i$ be a nonempty, closed convex subset of $M^i$
which has a diameter bounded by $D^i$ and $g_j^i \colon C^i \to \mathbb{R}$ ($i\in\mathcal{I},j\in \mathcal{J}^i$) be a convex function with 
$D(g_j^i) = C^i$ ($j\in \mathcal{J}^i$)
and $\bigcap_{j\in \mathcal{J}^i} \mathrm{lev}_{\leq 0} (g_j^i) \neq \emptyset$.
Then, Problem \ref{problem:1} with a mapping 
$T^i := P_{g_1^i,\lambda^i} P_{g_2^i,\lambda^i} \cdots P_{g_{J^i}^i,\lambda^i}$
($i\in\mathcal{I}$) is to find a point $x_\star$ in $X_\star$ with
\begin{align*}
X = \bigcap_{j\in \mathcal{J}^1} \mathrm{lev}_{\leq 0} \left(g_j^1\right) \times 
\bigcap_{j\in \mathcal{J}^2} \mathrm{lev}_{\leq 0} \left(g_j^2 \right) 
\times \cdots \times \bigcap_{j\in \mathcal{J}^I} \mathrm{lev}_{\leq 0} \left(g_j^I \right) ,
\end{align*}
where $\lambda^i \in (0,2/\zeta^i)$ and $\zeta^i := \sqrt{|\kappa^i|}D^i/\tanh (\sqrt{|\kappa^i|}D^i)$ ($i\in\mathcal{I}$). 
\end{example}
 
References \cite{bento2012} and \cite{wang2015_1} proposed Riemannian subgradient algorithms for finding a point $x^*$ in the intersection of sublevel sets of convex functions
$g_j$ ($j=1,2,\ldots,J$) defined on a Riemannian manifold,
i.e., 
\begin{align*}
x^* \in \bigcap_{j=1}^J \mathrm{lev}_{\leq 0} \left(g_j\right).
\end{align*}
Algorithm 3.1 in \cite{wang2015_1} converges linearly to $x^*$ without assuming 
that the domain of $g_j$ has a bounded diameter.
Meanwhile, under the assumption that the domain of $g_j$ has a bounded diameter,
Example \ref{exp:3} enables us to consider the three-stage Riemannian optimization problem
such that the first stage is to find points in $M$, 
the second stage is to find points in sublevel sets of convex functions, 
and the third stage is to minimize a function over the second stage.

This section ends with a statement of the conditions for being able to solve Problem \ref{problem:1}
(see, e.g., \cite[(A1), (A2), (2.5)]{nem2009}). 

\begin{enumerate}
\item[(C1)]
There is an independent and identically distributed sample 
$\xi_0, \xi_1, \ldots$ of realizations of the random vector $\xi$;
\item[(C2)]
There is an oracle which, for a given input point $(x, \xi) \in M \times \Xi$,
returns a stochastic (sub)gradient $\mathsf{G}(x,\xi) = (\mathsf{G}^i(x,\xi))_{i\in \mathcal{I}}$ 
%\in \partial_x F(x,\xi)$ 
such that 
\begin{align*}
\mathsf{g}(x) = (\mathsf{g}^i(x))_{i\in \mathcal{I}} := \mathbb{E}[\mathsf{G}(x,\xi)]
\begin{cases}
\in \partial f (x) \text{ (} f \text{ is nonsmooth and convex),}\\
= \mathrm{grad} f(x) \text{ (} f \text{ is smooth and nonconvex);}
\end{cases}
\end{align*}
\item[(C3)]
For all $i\in\mathcal{I}$, there exists a positive number $B^i$ such that,
for all $x\in M$, $\mathbb{E}[\|\mathsf{G}^i(x,\xi)\|_{x^i}^2] \leq {B^i}^2$.
\end{enumerate}

%Suppose that $F(\cdot, \xi)$ $(\xi \in \Xi)$ is convex
%and consider the oracle which returns a stochastic subgradient 
%$\mathsf{G}(x,\xi) \in \partial_x F (x,\xi)$ for a given $(x,\xi) \in \mathbb{R}^N \times \Xi$.
%Then, $f (\cdot) = \mathbb{E}[F(\cdot,\xi)]$ is well defined and convex, and $\partial f (x) = \mathbb{E}[\partial_x F(x,\xi)]$ \cite[Theorem 7.51]{shap2014}, \cite[p.1575]{nem2009}. 

\section{Riemannian Stochastic Fixed Point Optimization Algorithm}\label{sec:4}
Let $i\in\mathcal{I}$. 
Given a quasinonexpansive mapping $T^i \colon M^i \to M^i$ in Problem \ref{problem:1} and $\alpha^i \in (0,1)$,
we define $S_{\alpha^i}^i \colon M^i \to M^i$ for all $x^i \in M^i$
by 
\begin{align}\label{1}
S_{\alpha^i}^i (x^i) := \exp_{x^i}^i \left[\left(1 - \alpha^i \right) \left(\exp_{x^i}^{i} \right)^{-1} \left(T^i(x^i) \right) \right].
\end{align}
The discussion in Subsection \ref{subsec:2.3} (see \eqref{keyineq}) ensures that 
$S_{\alpha^i}^i$ is strictly quasinonexpansive with 
$\mathrm{Fix}(T^i) = \mathrm{Fix}(S_{\alpha^i}^i)$.
Moreover, we define 
\begin{align}\label{2}
Q_{\alpha^i}^i := P^i S_{\alpha^i}^i,
\end{align} 
where $P^i$ is the metric projection onto a nonempty, closed convex
set $C^i$ satisfying 
\begin{align}\label{3}
C^i\supset \mathrm{Fix} \left(S_{\alpha^i}^i \right) 
= \mathrm{Fix} \left(T^i \right).
%= \mathrm{Fix} \left(Q_{\alpha^i}^i \right).
\end{align}
%The second equation in \eqref{3} comes from 
%$\mathrm{Fix}(P^i) \cap \mathrm{Fix}(T^i) \neq \emptyset$.

Algorithm \ref{algo:1} is the proposed algorithm for solving Problem \ref{problem:1}. 
The tangent vectors $m_n$ and $\hat{m}_n$ generated by steps 3 and 4 in Algorithm \ref{algo:1} are based on so-called momentum terms \cite[Subchapter 8.3.2]{deep}.
Step 8 in Algorithm \ref{algo:1} is expressed as 
\begin{align*}
x_{n+1}^i := Q_{\alpha^i}^i \left[ \exp_{x_n^i}^i 
\left(- \frac{\alpha_n}{(1 - \hat{\beta}^{n+1}) \mathsf{h}_n^i} 
\hat{m}_n^i  \right) \right],
\end{align*} 
which implies that Algorithm \ref{algo:1} adapts the step-size
$\alpha_n/((1 - \hat{\beta}^{n+1})\mathsf{h}_n^i)$
for each $i\in\mathcal{I}$ and each $n\in\mathbb{N}$.
Hence, we can see that Algorithm \ref{algo:1} is based on so-called {\em adaptive learning rate optimization algorithms}, such as AdaGrad \cite{adagrad}, Adam \cite{adam}, and AMSGrad \cite{reddi2018} defined on Euclidean space
and RAMSGrad \cite{gec2019} defined on a Riemannian manifold. 
Examples of $\mathsf{h}_n^i$ are included in Examples \ref{exp:adam} and \ref{exp:amsgrad}.

\begin{algorithm}
\caption{Riemannian stochastic fixed point optimization algorithm}
\label{algo:1}
\begin{algorithmic}[1]
\REQUIRE 
$(\alpha_n)_{n\in\mathbb{N}} \subset (0,1)$, $(\alpha^i)_{i\in\mathcal{I}} \subset (0,1)$, 
$(\beta_n)_{n\in\mathbb{N}} \subset [0,1)$, $\hat{\beta} \in [0,1)$
\STATE{$n \gets 0$, ${x}_0 \in M$, $\tau_{-1} = m_{-1}  \in T_{x_0} M$, 
$(\mathsf{h}_0^i)_{i\in \mathcal{I}} \subset \mathbb{R}_{++}^I$}
\LOOP
\STATE{$m_{n} := \beta_n \tau_{n-1} + (1 - \beta_n) \mathsf{G}(x_n,\xi_n)$}
\STATE{$\hat{m}_n := \left(1 - \hat{\beta}^{n+1} \right)^{-1} m_n$}
\STATE{$(\mathsf{h}_n^i)_{i\in\mathcal{I}} \subset \mathbb{R}^I_{++}$}
\FOR{$i=1,2,\ldots,I$}
\STATE{$\displaystyle{\mathsf{d}_n^i = - \frac{\hat{m}_{n}^i}{\mathsf{h}_n^i}}$}
\STATE{$\displaystyle{x_{n+1}^i := Q_{\alpha^i}^i \left[ \exp_{x_n^i}^i \left(\alpha_n \mathsf{d}_n^i \right) \right]}$}
\STATE{$\tau_n^i := \varphi_{x_n^i \to x_{n+1}^i}^i (m_n^i)$}
\STATE{$n \gets n+1$}
\ENDFOR
\ENDLOOP
%\RETURN $T$
\end{algorithmic}
\end{algorithm}
%\end{tcbverbatimwrite}

The following conditions are assumed to analyze Algorithm \ref{algo:1}.

\begin{assum}\label{assum:1}
The sequence $(\mathsf{H}_n)_{n\in\mathbb{N}} := ((\mathsf{h}_n^i)_{i\in\mathcal{I}})_{n\in\mathbb{N}}$
and a nonempty, closed convex set $C^i \supset \mathrm{Fix}(T^i)$
in Algorithm \ref{algo:1} satisfy the following conditions:
\begin{enumerate}
\item[{\em (A3)}]
For all $i\in \mathcal{I}$, $C^i$ has a diameter bounded by $D^i$;
\item[{\em (A4)}]
For all $n\in\mathbb{N}$ and all $i\in\mathcal{I}$, almost surely
$\mathsf{h}_{n+1}^i \geq \mathsf{h}_n^i$;
\item[{\em (A5)}]
For all $i\in\mathcal{I}$, there exists a positive number $\hat{B}^i$ such that, for all $n\in\mathbb{N}$, 
$\mathbb{E}[\mathsf{h}_n^i] \leq \hat{B}^i$.
\end{enumerate}
\end{assum}

Assumption (A3) will be needed to analyze Algorithm \ref{algo:1} 
since the previously reported results
were analyzed under Assumption (A3) (see, e.g., 
\cite[p.1574]{nem2009} and \cite[p.2]{reddi2018} for convex stochastic optimization on Euclidean space, and see, e.g., 
\cite[Section 4]{gec2019}, \cite[Subsection 3.2]{liu2017}, and \cite[Subsection 3.2]{zhang2016} for convex stochastic optimization on a Riemannian manifold). 
%If $C^i$ ($i\in\mathcal{I}$) is bounded, then Assumption (A3) holds.
For example, let us consider the Poincar\'e model of a hyperbolic space defined by a manifold $\mathcal{D}^{m^i} := \{ x^i \in \mathbb{R}^{m^i} \colon \|x^i \| < 1 \}$ 
equipped with the Riemannian metric $\rho_{x^i} := (1/(1-\|x^i\|^2)^2) \rho_{x^i}^{\mathrm{E}}$,
where $\|\cdot\|$ is the Euclidean norm, $x^i\in \mathcal{D}^{m^i}$, and $\rho_{x^i}^{\mathrm{E}}$ is the Euclidean metric tensor
(the Poincar\'e embedding has been used for natural language processing 
\cite[Section 5]{gec2019}, \cite[Section 4]{sakai2020}).
In \cite[Section 4]{sakai2020}, 
$C^i := \{ x^i \in  \mathcal{D}^{m^i} \colon \|x^i\| \leq 1 - 10^{-5} \}$,
which has a bounded diameter,
was used to evaluate the performance of RAMSGrad for natural language processing.
In the case of Example \ref{exp:1}, 
Assumption (A3) is satisfied 
when at least one of $C_j^i$ ($j\in \mathcal{J}^i$) 
has a bounded diameter.
In Example \ref{exp:3}, 
since $\mathrm{Fix}(T^i) = \bigcap_{j\in \mathcal{J}^i} \mathrm{lev}_{\leq 0}
(g_j^i) \subset C^i$ and $C^i$ has a bounded diameter, 
Assumption (A3) is satisfied.
Assumption (A3) implies that, for all $i\in \mathcal{I}$, 
\begin{align}\label{Di}
D^i := \sup \left\{ \mathrm{d}^i \left( x^i, y^i \right) \colon (x^i)_{i\in\mathcal{I}}, (y^i)_{i\in\mathcal{I}} \in X \right\}
< + \infty, 
\end{align}
where $\mathrm{d}^i \colon M^i \times M^i \to \mathbb{R}_+$ is the distance function of $M^i$.

Under Assumption (A3), we provide some examples of $(\mathsf{H}_n)_{n\in\mathbb{N}}$ 
satisfying Assumptions (A4) and (A5).
The following examples are based on adaptive learning rate 
optimization algorithms, such as Adam \cite{adam} and AMSGrad \cite{reddi2018}, defined on Euclidean space.

\begin{example}[$\mathsf{H}_n$ based on Adam \cite{adam}]\label{exp:adam} 
Let us define $\mathsf{h}_n^i$ and $v_n^i$ for all $i\in \mathcal{I}$ 
and all $n\in \mathbb{N}$ by 
\begin{align}\label{vn_adam}
\begin{split}
&v_n^i := \bar{\beta} v_{n-1}^i + (1 - \bar{\beta}) 
\left\|\mathsf{G}^i(x_n,\xi_n) \right\|_{x_n^i}^2,\\
&\bar{v}_n^i := \frac{v_n^i}{1 - \bar{\beta}^{n+1}}, \text{ }
\hat{v}_{n}^i := \max \left\{   
\hat{v}_{n-1}^i, \bar{v}_n^i \right\},\\
&\mathsf{h}_n^i := \sqrt{\hat{v}_{n}^i},
\end{split}
\end{align}
where $v_{-1}^i, \hat{v}_{-1}^i \in \mathbb{R}_+$ and $\bar{\beta} \in [0,1)$. 
From \eqref{vn_adam}, $\mathsf{H}_n$ satisfies Assumption (A4).
Moreover, \eqref{1}, \eqref{2}, and \eqref{3} mean that $(x_n^i)_{n\in\mathbb{N}} \subset C^i$,
which, together with Assumption (A3), implies that $(\|\mathsf{G}^i (x_n, \xi_n)\|_{x_n^i})_{n\in\mathbb{N}}$ is almost surely 
bounded \cite[Lemma 3.3]{gro2016}.
For all $i\in\mathcal{I}$, we define 
\begin{align*}
U^i := \max \left\{  v_{-1}^i , 
 \sup \left\{ \left\|\mathsf{G}^i (x_n, \xi_n) \right\|_{x_n^i}^2 \colon n\in\mathbb{N}  \right\} \right\} < + \infty.
\end{align*}
Induction, together with the definitions of $v_n^i$, $\bar{v}_n^i$, and $\bar{\beta} \in [0,1)$, implies that, for all $n\in\mathbb{N}$,
%\begin{align*}
$v_n^i \leq U^i$ and  
$\bar{v}_n^i \leq U^i/(1 - \bar{\beta})$.
%\end{align*} 
Accordingly, induction ensures that 
\begin{align*}
%\hat{v}_n^i \leq \max \left\{ \hat{v}_{-1}^i, \frac{U^i}{1 - \bar{\beta}}   \right\}, 
%\text{ i.e., } 
\mathbb{E}\left[\mathsf{h}_n^i \right] \leq \sqrt{\max \left\{ \hat{v}_{-1}^i, \frac{U^i}{1 - \bar{\beta}}   \right\}},
\end{align*}
which implies that Assumption (A5) holds. 
\end{example}

\begin{example}[$\mathsf{H}_n$ based on AMSGrad \cite{gec2019,reddi2018}]
\label{exp:amsgrad}
Let us define $\mathsf{h}_n^i$ and $v_n^i$ for all $i\in \mathcal{I}$ 
and all $n\in \mathbb{N}$ by 
\begin{align}\label{vn_amsgrad}
\begin{split}
&v_n^i := \bar{\beta} v_{n-1}^i + (1 - \bar{\beta}) 
\left\|\mathsf{G}^i(x_n,\xi_n) \right\|_{x_n^i}^2,\\
&\hat{v}_{n}^i := \max \left\{   
\hat{v}_{n-1}^i, {v}_n^i \right\},\\
&\mathsf{h}_n^i := \sqrt{\hat{v}_{n}^i},
\end{split}
\end{align}
where $v_{-1}^i, \hat{v}_{-1}^i \in \mathbb{R}_+$ and $\bar{\beta} \in [0,1)$. 
The same discussion as in Example \ref{exp:adam} ensures that 
$\mathsf{h}_n^i$ defined by \eqref{vn_amsgrad} satisfies Assumptions (A4) and (A5),
i.e.,
\begin{align*}
\mathbb{E}\left[\mathsf{h}_n^i \right] \leq \sqrt{\max \left\{ \hat{v}_{-1}^i, U^i   \right\}}.
\end{align*}
\end{example}

\section{Convergence analyses of Algorithm \ref{algo:1}}\label{sec:5}
\subsection{Nonsmooth convex optimization}\label{subsec:5.1}
This subsection considers Problem \ref{problem:1} when $f$ is nonsmooth and convex.
The following is a convergence analysis of Algorithm \ref{algo:1} with constant 
step-sizes (the proof of the following theorem is given in Supplementary Material).

\begin{theorem}\label{cor:1}
Suppose that Assumptions (A1)--(A5) and Conditions (C1)--(C3) hold.
Then, Algorithm \ref{algo:1} with $\alpha_n := \alpha$ and $\beta_n := \beta$ satisfies that, for all $i\in \mathcal{I}$,
\begin{align}
&\limsup_{n\to + \infty}
\mathbb{E}\left[ \mathrm{d}^i \left(y_n^i,x_n^i \right)^2 \right]
\leq
\frac{\tilde{{B}^i}^2}{(1-\hat{\beta})^2(\mathsf{h}_0^i)^2} \alpha^2,\label{limsup}\\
&\liminf_{n \to + \infty}
\mathbb{E} \left[
\mathrm{d}^i \left( T^i (y_n^i), y_n^i \right)^2 \right]
\leq 
\frac{1}{\hat{\alpha}^i}
\left\{
\frac{2 \tilde{B}^i D^i}{(1-\hat{\beta})\mathsf{h}_0^i} \alpha
+
\frac{\zeta^i  \tilde{{B}^i}^2}{(1-\hat{\beta})^2(\mathsf{h}_0^i)^2} \alpha^2
\right\},\label{liminf}
\end{align}
and 
\begin{align}\label{liminf_2}
\liminf_{n\to + \infty} \mathbb{E}\left[ f(x_n) - f^\star \right]
\leq \frac{\sum_{i\in\mathcal{I}} \zeta^i \tilde{{B}^i}^2 (\mathsf{h}_0^i)^{-1}}{2(1-\beta)(1-\hat{\beta})} \alpha
+
\frac{\sum_{i\in\mathcal{I}} \tilde{B}^i D^i}{(1-\beta)(1-\hat{\beta})} \beta,
\end{align}
where $\hat{\alpha}_i := \alpha^i(1-\alpha^i)$ and 
$\tilde{{B}^i}^2 := \max \{ \|\tau_{-1}\|_{x_0^i}^2, {{B}^i}^2 \}$.
Moreover, for all $i\in\mathcal{I}$ and all $n \geq 1$, 
\begin{align}\label{ave_1}
\mathbb{E} \left[ \frac{1}{n} \sum_{k=1}^n
\mathrm{d}^i \left(T^i (y_k^i), y_k^i \right)^2 \right]
\leq
\frac{D^i}{\hat{\alpha}_i}
\frac{1}{n}
+
\frac{2 \tilde{B}^i D^i}{\hat{\alpha}_i \hat{\mathsf{h}}_0^i} \alpha
+ 
\frac{\zeta^i \tilde{{B}^i}^2}{\hat{\alpha}_i(\hat{\mathsf{h}}_0^i)^2}
\alpha^2,
\end{align}
where $\hat{\mathsf{h}}_0^i := (1 - \hat{\beta}) \mathsf{h}_0^i$.
Let us define $\bar{x}_n$ for all $n\geq 1$ by  
\begin{align}\label{ave}
\bar{x}_n := \exp_{\bar{x}_{n-1}}\left( 
\frac{1}{n} \exp_{\bar{x}_{n-1}}^{-1}(x_n)    
\right),
\end{align}
where $\bar{x}_0 := x_0$. 
Then, for all $n\geq 1$,
\begin{align}\label{ave_2}
\mathbb{E} \left[ f (\bar{x}_n) - f_\star \right]
&\leq 
\frac{\sum_{i\in \mathcal{I}} \hat{B}^i {D^i}^2}{2 (1-\beta_1)} \frac{1}{\alpha n}
+
\frac{\sum_{i\in\mathcal{I}} \zeta^i \tilde{{B}^i}^2 (\mathsf{h}_0^i)^{-1}}{2(1-\hat{\beta})(1-\beta_1)} 
\alpha
+
\frac{\sum_{i\in\mathcal{I}} \tilde{B}^i D^i}{1 - \beta_1} 
\beta.
\end{align}
If
{\em (A1)'}
$T^i \colon M^i \to M^i$ ($i\in \mathcal{I}$) is nonexpansive with $\mathrm{Fix}(T^i) \neq
\emptyset$, then
\begin{align}\label{av_0}
\mathbb{E} \left[ \frac{1}{n} \sum_{k=1}^n
\mathrm{d}^i \left(T^i (x_k^i), x_k^i \right)^2 \right]
 \leq 
\frac{2 D^i}{\hat{\alpha}^i}
\frac{1}{n}
+
\frac{4 \tilde{B}^i D^i}{\hat{\alpha}^i \hat{\mathsf{h}}_0^i} \alpha
+ 
\frac{2 \tilde{{B}^i}^2}{(\hat{\mathsf{h}}_0^i)^2}
\left\{
\frac{\zeta^i}{\hat{\alpha}^i}
+ \frac{4}{(1-\hat{\beta})^2}
\right\}
\alpha^2.
\end{align}
\end{theorem}

%Theorem \ref{cor:1} ensures that Algorithm \ref{algo:1} satisfies that there exist positive real numbers $\mathcal{C}_j$
%($j=1,2,\ldots,9$) such that
%\begin{align}\label{constant}
%\begin{split}
%&\limsup_{n \to + \infty}
%\mathbb{E} \left[
%\mathrm{d}^i \left( y_n^i, x_n^i \right)^2 \right]
%\leq 
%\mathcal{C}_1 \alpha^2,\text{ }
%\liminf_{n \to + \infty}
%\mathbb{E} \left[
%\mathrm{d}^i \left( T^i (y_n^i), y_n^i \right)^2 \right]
%\leq 
%\mathcal{C}_2 \alpha
%+
%\mathcal{C}_3 \alpha^2,\\
%&\liminf_{n\to + \infty} \mathbb{E}\left[ f(x_n) - f^\star \right]
%\leq 
%\mathcal{C}_4 \alpha
%+
%\mathcal{C}_5 \beta,
%\end{split}
%\end{align}
%and
%\begin{align}\label{const_1}
%\begin{split}
%&\mathbb{E} \left[ \frac{1}{n} \sum_{k=1}^n \sum_{i\in\mathcal{I}}
%\mathrm{d}^i \left(T^i (y_k^i), y_k^i \right)^2 \right]
%\leq 
%\mathcal{O} \left(\frac{1}{n} \right)
%+
%\mathcal{C}_6 \alpha + \mathcal{C}_7 \alpha^2,\\
%&
%\mathbb{E} \left[ f (\bar{x}_n) - f_\star \right]
%\leq 
%\mathcal{O} \left( \frac{1}{n} \right)
%+
%\mathcal{C}_8 \alpha
%+
%\mathcal{C}_9 \beta.
%\end{split}
%\end{align}
%Therefore, \eqref{constant} and \eqref{const_1} indicate that 
%Algorithm \ref{algo:1} with sufficiently small step-sizes  
%$\alpha$ and $\beta$ approximates a solution
%to Problem \ref{problem:1}.

The following is a convergence analysis of Algorithm \ref{algo:1}
with diminishing step-sizes (the proof of the theorem is given in
Supplementary Material).

\begin{theorem}\label{cor:2}
Suppose that Assumptions (A1)--(A5) and Conditions (C1)--(C3) hold
and assume that $(\alpha_n)_{n\in\mathbb{N}}$ is monotone decreasing
and $(\alpha_n)_{n\in\mathbb{N}}$ and $(\beta_n)_{n\in\mathbb{N}}$ satisfy that
\begin{align}\label{step:0}
\sum_{n=0}^{+\infty} \alpha_n = + \infty,
\sum_{n=0}^{+\infty} \alpha_n^2 < + \infty,
\text{ and }
\sum_{n=0}^{+\infty} \alpha_n \beta_n <  + \infty.
\end{align}
Then, Algorithm \ref{algo:1} satisfies that, for all $i\in\mathcal{I}$,
\begin{align}\label{Ty}
&\lim_{n \to +\infty} \mathbb{E}\left[ \mathrm{d}^i \left(y_n^i, x_n^i \right)^2 \right] = 0,\text{ } \liminf_{n \to +\infty} \mathbb{E}\left[ \mathrm{d}^i \left(T^i(y_n^i), y_n^i \right)^2 \right] = 0,
\end{align}
and 
\begin{align*}
&\liminf_{n \to +\infty} \mathbb{E}\left[ f(x_n) - f_\star \right] \leq 0.
\end{align*}
Suppose that Assumptions (A1)--(A5) and Conditions (C1)--(C3) hold
and assume that $(\alpha_n(1-\beta_n))_{n\in\mathbb{N}}$ and $(\beta_n)_{n\in\mathbb{N}}$ are monotone decreasing
and satisfy the following:
\begin{align}\label{step}
\lim_{n\to + \infty} \frac{1}{n\alpha_n} = 0,
\text{ }
\lim_{n \to + \infty} \frac{1}{n} \sum_{k=1}^n \alpha_k = 0,
\text{ and }
\lim_{n \to + \infty} \frac{1}{n} \sum_{k=1}^n \beta_k = 0.
\end{align}
Then, Algorithm \ref{algo:1} satisfies that
\begin{align}\label{main:4}
\lim_{n\to + \infty} \mathbb{E} \left[ \frac{1}{n} \sum_{k=1}^n \sum_{i\in\mathcal{I}}
\mathrm{d}^i \left(T^i (y_k^i), y_k^i \right)^2 \right] = 0 
\end{align}
and
\begin{align}\label{main:4_1}
\limsup_{n\to +\infty} \mathbb{E} \left[  f (\bar{x}_n) - f_\star \right] \leq 0
\end{align} 
with the rate of convergence expressed as follows: 
\begin{align*}
&\mathbb{E} \left[ \frac{1}{n} \sum_{k=1}^n
\mathrm{d}^i \left(T^i (y_k^i), y_k^i \right)^2 \right]
\leq 
\frac{1}{\hat{\alpha}^i}
\left\{
\frac{D^i}{n}
+
\frac{2 \tilde{B}^i D^i}{\hat{\mathsf{h}}_0^i} \frac{1}{n} \sum_{k=1}^n \alpha_k
+ 
\frac{\zeta^i \tilde{{B}^i}^2}{(\hat{\mathsf{h}}_0^i)^2}
\frac{1}{n} \sum_{k=1}^n \alpha_k^2
\right\},\\
&\mathbb{E} \left[ \frac{1}{n} \sum_{k=1}^n \mathrm{d}^i \left(y_k^i, x_k^i \right)^2 \right] 
\leq 
\frac{\tilde{{B}^i}^2}{(1-\hat{\beta})^2(\hat{\mathsf{h}}_0^i)^2}
\frac{1}{n} \sum_{k=1}^n \alpha_k^2,
\end{align*}
and
\begin{align*}
\mathbb{E} \left[ f(\bar{x}_n) - f_\star  \right]
&\leq 
\frac{\sum_{i\in \mathcal{I}} \hat{B}^i {D^i}^2}{2 (1-\beta_1)} \frac{1}{n \alpha_n}
+
\frac{\sum_{i\in\mathcal{I}} \zeta^i \tilde{{B}^i}^2 (\mathsf{h}_0^i)^{-1}}{2(1-\hat{\beta})(1-\beta_1)
} 
\frac{1}{n} \sum_{k=1}^n \alpha_k\\
&\quad +
\frac{\sum_{i\in\mathcal{I}} \tilde{B}^i D^i}{1 - \beta_1} 
\frac{1}{n} \sum_{k=1}^n \beta_k,
\end{align*}
where $\bar{x}_n$ is defined by \eqref{ave}.
Under Assumption (A1)', we have 
\begin{align*}
&\mathbb{E} \left[ \frac{1}{n} \sum_{k=1}^n
\mathrm{d}^i \left(T^i (x_k^i), x_k^i \right)^2 \right]\\
&\quad \leq 
\frac{2}{\hat{\alpha}^i}
\frac{D^i}{n}
+
\frac{4 \tilde{B}^i D^i}{\hat{\alpha}^i \hat{\mathsf{h}}_0^i} \frac{1}{n} \sum_{k=1}^n \alpha_k
+ 
\frac{2 \tilde{{B}^i}^2}{(\hat{\mathsf{h}}_0^i)^2}
\left\{
\frac{\zeta^i}{\hat{\alpha}^i}
+ \frac{4}{(1-\hat{\beta})^2}
\right\}
\frac{1}{n} \sum_{k=1}^n \alpha_k^2.
\end{align*}
\end{theorem}

Theorem \ref{cor:2} yields the following corollary.

\begin{corollary}\label{COR:1_1}
Suppose that Assumptions (A1)--(A5) and Conditions (C1)--(C3) hold.
Then, Algorithm \ref{algo:1} with $\alpha_n := 1/n^\eta$ 
($\eta \in (1/2,1], n \geq 1$) and  
$(\beta_n)_{n\in\mathbb{N}}$ such that $\sum_{n=1}^{+\infty} \alpha_n \beta_n < + \infty$\footnote{The step-sizes 
$\beta_n := \lambda^{n}$ and $\alpha_n := 1/n^\eta$ ($n \geq 1, \lambda \in (0,1), \eta \in (1/2,1]$) satisfy $\sum_{n=1}^{+\infty} \alpha_n = + \infty$,
$\sum_{n=1}^{+\infty} \alpha_n^2 < + \infty$, and  
$\sum_{n=1}^{+\infty} \alpha_n \beta_n < + \infty$.}
satisfies that, for all $i\in\mathcal{I}$,  
$\mathbb{E} [\mathrm{d}^i (y_n^i, x_n^i)^2 ]
= \mathcal{O}( n^{- 2 \eta})$,
\begin{align}\label{main:5}
\liminf_{n \to + \infty} \mathbb{E} \left[ 
\mathrm{d}^i \left(T^i(y_n^i), y_n^i \right)^2  
\right]
= 0, 
\text{ and } 
\liminf_{n \to + \infty} \mathbb{E} \left[ 
f(x_n) - f_\star
\right]
\leq 0. 
\end{align}
Moreover, Algorithm \ref{algo:1} with $\alpha_n := 1/n^\eta$ ($\eta \in [1/2,1)$) and  
$(\beta_n)_{n\in\mathbb{N}}$ such that $\sum_{n=1}^{+\infty} \beta_n < + \infty$\footnote{The step-sizes 
$\beta_n := 1/2^{n}$ and $\alpha_n := 1/n^\eta$ ($n \geq 1, \eta \in [1/2,1)$) are used to implement adaptive learning rate optimization algorithms, such as Adam \cite{adam}, AMSGrad \cite{reddi2018},
and RAMSGrad \cite{gec2019}.
These step-sizes satisfy $\sum_{n=1}^{+\infty} \beta_n = 1$ and 
$(\alpha_n(1-\beta_n))_{n\in\mathbb{N}}$ is monotone decreasing.
}
satisfies that, for all $n \geq 1$,
\begin{align*}
%&\mathbb{E} \left[ \sum_{i\in\mathcal{I}}
%\mathrm{d}^i \left(y_k^i, x_k^i \right)^2 \right]
%= \mathcal{O}\left( \frac{1}{n^{2\eta}} \right),\\
\mathbb{E} \left[ \frac{1}{n} \sum_{k=1}^n \sum_{i\in\mathcal{I}}
\mathrm{d}^i \left(T^i (y_k^i), y_k^i \right)^2 \right]
=
\mathcal{O}\left( \frac{1}{n^\eta}\right)
\end{align*}
and  
\begin{align}\label{rate}
\mathbb{E} \left[ f (\bar{x}_n) - f_\star \right]
\leq 
\mathcal{O}\left( \frac{1}{n^{1-\eta}}\right),
\end{align}
where $\bar{x}_n$ is defined by \eqref{ave}.
Under Assumption (A1)', we have 
\begin{align*}
&\mathbb{E} \left[ \frac{1}{n} \sum_{k=1}^n \sum_{i\in\mathcal{I}}
\mathrm{d}^i \left(T^i (x_k^i), x_k^i \right)^2 \right]
=
\mathcal{O}\left( \frac{1}{n^\eta}\right).
\end{align*}
\end{corollary}

\subsection{Smooth nonconvex optimization}
This subsection considers Problem \ref{problem:1} when $f$ is smooth and nonconvex.
The following is a convergence analysis of Algorithm \ref{algo:1} with constant step-sizes (the proof of the theorem is given in Supplementary Material).

\begin{theorem}\label{cor:2_1}
Suppose that Assumptions (A1)--(A5) and Conditions (C1)--(C3) hold.
Then, Algorithm \ref{algo:1} with $\alpha_n := \alpha$ and $\beta_n := \beta$ satisfies that, for all $i\in \mathcal{I}$,
\eqref{limsup} and \eqref{liminf} hold,
and 
\begin{align}\label{liminf_3}
\begin{split}
\limsup_{n\to + \infty} 
\mathbb{E}\left[ \left\langle \exp_{x_n}^{-1}(x), \mathrm{grad}f(x_n) \right\rangle_{x_n} \right]
\geq - \frac{\sum_{i\in\mathcal{I}} \zeta^i \tilde{{B}^i}^2}{2\gamma \mathsf{h}_0^i} \alpha -
\frac{\sum_{i\in\mathcal{I}} \tilde{B}^i D^i}{\gamma} \beta,
\end{split}
\end{align}
where $\gamma := (1-\beta)(1-\hat{\beta})$.
Moreover, for all $i\in\mathcal{I}$ and  all $n \geq 1$, 
\eqref{ave_1} holds and
\begin{align}\label{ave_3}
\begin{split}
&\mathbb{E} \left[ \frac{1}{n} \sum_{k=1}^n 
\left\langle \exp_{x_k}^{-1}(x), \mathrm{grad}f(x_k) \right\rangle_{x_k}
\right]\\
&\geq 
- \frac{\sum_{i\in \mathcal{I}} \hat{B}^i {D^i}^2}{2 (1-\beta_1)} \frac{1}{\alpha n}
-
\frac{\sum_{i\in\mathcal{I}} \zeta^i \tilde{{B}^i}^2 (\mathsf{h}_0^i)^{-1}}{2(1-\hat{\beta})(1-\beta_1)} 
\alpha
-
\frac{\sum_{i\in\mathcal{I}} \tilde{B}^i D^i}{1 - \beta_1} 
\beta.
\end{split}
\end{align}
Under
Assumption (A1)', \eqref{av_0} holds.
\end{theorem}

The following is a convergence analysis of Algorithm \ref{algo:1} with diminishing step-sizes (the proof of the theorem is given in Supplementary Material).

\begin{theorem}\label{cor:2_2}
Suppose that Assumptions (A1)--(A5) and Conditions (C1)--(C3) hold
and assume that $(\alpha_n)_{n\in\mathbb{N}}$ is monotone decreasing
and $(\alpha_n)_{n\in\mathbb{N}}$ and $(\beta_n)_{n\in\mathbb{N}}$ satisfy
\eqref{step:0}.
Then, Algorithm \ref{algo:1} satisfies that, for all $i\in\mathcal{I}$,
\eqref{Ty} holds
and 
\begin{align*}
&\limsup_{n \to +\infty} \mathbb{E}\left[\left\langle \exp_{x_n}^{-1}(x), \mathrm{grad}f(x_n) \right\rangle_{x_n}  \right] \geq 0.
\end{align*}
Suppose that Assumptions (A1)--(A5) and Conditions (C1)--(C3) hold
and assume that $(\alpha_n(1-\beta_n))_{n\in\mathbb{N}}$ and $(\beta_n)_{n\in\mathbb{N}}$ are monotone decreasing
and satisfy \eqref{step}.
Then, Algorithm \ref{algo:1} satisfies that \eqref{main:4} holds and 
\begin{align*}
%\lim_{n\to + \infty} \mathbb{E} \left[ \frac{1}{n} \sum_{k=1}^n \sum_{i\in\mathcal{I}}
%\mathrm{d}^i \left(T^i (y_k^i), y_k^i \right)^2 \right] = 0
%\text{ and }
\liminf_{n\to +\infty} \mathbb{E} \left[ \frac{1}{n} \sum_{k=1}^n  
\left\langle \exp_{x_k}^{-1}(x), \mathrm{grad}f(x_k) \right\rangle_{x_k}  \right] 
\geq 0
\end{align*} 
with 
\begin{align*}
&\mathbb{E} \left[ \frac{1}{n} \sum_{k=1}^n  
\left\langle \exp_{x_k}^{-1}(x), \mathrm{grad}f(x_k) \right\rangle_{x_k}  \right]\\
&\geq 
-\frac{\sum_{i\in \mathcal{I}} \hat{B}^i {D^i}^2}{2 (1-\beta_1)} \frac{1}{n \alpha_n}
-
\frac{\sum_{i\in\mathcal{I}} \zeta^i \tilde{{B}^i}^2 (\mathsf{h}_0^i)^{-1}}{2(1-\hat{\beta})(1-\beta_1)
} 
\frac{1}{n} \sum_{k=1}^n \alpha_k
-
\frac{\sum_{i\in\mathcal{I}} \tilde{B}^i D^i}{1 - \beta_1} 
\frac{1}{n} \sum_{k=1}^n \beta_k
\end{align*}
and the same convergence rate of $\mathrm{d}^i(T^i(y_k^i),y_k^i)$ and $\mathrm{d}^i(T^i(x_k^i),x_k^i)$ (under Assumption (A1)') as in Theorem \ref{cor:2}.
\end{theorem}

A discussion similar to the one for obtaining Corollary \ref{COR:1_1} implies that Algorithm \ref{algo:1} with $\alpha_n := 1/n^\eta$
($\eta \in [1/2,1)$) and $(\beta_n)_{n\in \mathbb{N}}$ such that $\sum_{n=1}^{+\infty}\beta_n < + \infty$ satisfies that, under Assumption (A1)', 
\begin{align*}
\mathbb{E} \left[ \frac{1}{n} \sum_{k=1}^n \sum_{i\in\mathcal{I}}
\mathrm{d}^i \left(T^i (x_k^i), x_k^i \right)^2 \right]
=
\mathcal{O}\left( \frac{1}{n^\eta}\right)
\end{align*}
and
\begin{align*}
\mathbb{E} \left[ \frac{1}{n} \sum_{k=1}^n  
\left\langle \exp_{x_k}^{-1}(x), \mathrm{grad}f(x_k) \right\rangle_{x_k} \right]
\geq 
-\mathcal{O}\left( \frac{1}{n^{1-\eta}}\right).
\end{align*}

\section{Numerical Comparisons}\label{sec:6}
%This section considers two concrete examples of Problem \ref{problem:1} 
%and investigates 
%how Algorithm \ref{algo:1} works depending on the choices of 
%the parameters $\mathsf{h}_n^i$, $\alpha_n$, and $\beta_n$. 

\subsection{Preliminaries}
The $m$-dimensional Poincar\'e disk model of hyperbolic space 
is defined by 
\begin{align*}
\mathcal{D}^{m} := \left\{ x \in \mathbb{R}^{m} \colon \|x \| < 1 \right\},
\end{align*} 
where $\|\cdot\|$ denotes the Euclidean norm of $\mathbb{R}^{m}$. 
Let us also define $M := \underbrace{\mathcal{D}^{m} \times \mathcal{D}^{m} \times \cdots \times \mathcal{D}^{m}}_{I}$.
Let $j\in\mathcal{J}^i := \{1,2,\ldots, J^i\}$ ($i\in \mathcal{I} := \{1,2,\ldots,I\}$). 
We define a closed ball with center $c_j^i \in \mathcal{D}^{m}$ and radius $r_j^i > 0$ in $\mathcal{D}^{m}$ by 
\begin{align}\label{ball}
\mathrm{B}_j^i 
:= \left\{ x \in \mathcal{D}^{m} \colon
\mathrm{d} \left(c_j^i, x \right) \leq r_j^i   \right\},
\end{align}
where $\mathrm{d} \colon \mathcal{D}^{m} \times \mathcal{D}^{m}
\to \mathbb{R}$ denotes the distance function of $\mathcal{D}^{m}$.
Then, the metric projection onto the closed convex set
$\mathrm{B}_j^i$ can be expressed as follows:
\begin{align*}%\label{proj}
P_j^i (x) :=
\begin{cases} 
\exp_{c_j^i}^i \left( \frac{r_j^i \left(\exp^i_{c_j^i}\right)^{-1} (x)}
{\left\| \left(\exp_{c_j^i}^i \right)^{-1} (x) \right\|_{c_j^i}}  \right)
&\text{ if } x \notin \mathrm{B}_j^i,\\
x 
&\text{ if } x \in \mathrm{B}_j^i.
\end{cases}
\end{align*}
We used the nonexpansive mapping $T^i \colon D^{m} \to D^{m}$ ($i\in \mathcal{I}$) defined by
\begin{align}\label{nonexp}
T^i := P_1^i P_2^i \cdots P_{J^i}^i
\end{align}
and the smooth, nonconvex function $f \colon M \to \mathbb{R}$ defined for all $x\in M$ by
\begin{align*}
f(x) = \frac{1}{I} \sum_{i=1}^I \underbrace{\left\{ e^{(x^i)^\top x^j} + (x^i)^\top x^j  \right\}}_{= F(x, i)}, 
\text{ where } j := (i\bmod I) + 1.
\end{align*}
We implemented the following algorithms, all with $\alpha^i := 0.5$:\begin{itemize}
\item
Algorithm \ref{algo:1} with constant step-sizes
\begin{description}
\item[CSD: ]
Algorithm \ref{algo:1} with $\mathsf{h}_n^i$ defined by Stochastic Gradient Descent \cite{bon2013} (i.e., $\mathsf{h}_n^i := 1$), $\alpha_n := 10^{-2}$, and $\beta_n = \hat{\beta} := 0$
\item[CAG: ]
Algorithm \ref{algo:1} with $\mathsf{h}_n^i$ defined by AdaGrad \cite{adagrad}, 
$\alpha_n := 10^{-2}$, and $\beta_n = \hat{\beta} := 0$
\item[CAM1: ]
Algorithm \ref{algo:1} with $\mathsf{h}_n^i$ defined by AMSGrad \eqref{vn_amsgrad}, 
$\alpha_n := 10^{-2}$, $\beta_n = 0.9$, $\hat{\beta} := 0$, and $\bar{\beta} := 0.999$
\item[CAM2: ]
Algorithm \ref{algo:1} with $\mathsf{h}_n^i$ defined by AMSGrad \eqref{vn_amsgrad}, 
$\alpha_n := 10^{-2}$, $\beta_n = 10^{-3}$, $\hat{\beta} := 0$, and $\bar{\beta} := 0.999$
\item[CAD1: ]
Algorithm \ref{algo:1} with $\mathsf{h}_n^i$ defined by Adam \eqref{vn_adam}, 
$\alpha_n := 10^{-2}$, $\beta_n = 0.9$, $\hat{\beta} := 0.9$, and $\bar{\beta} := 0.999$ ($\hat{\beta} := 0.9$ and $\bar{\beta} := 0.999$ were used in \cite{adam,reddi2018})
\item[CAD2: ]
Algorithm \ref{algo:1} with $\mathsf{h}_n^i$ defined by Adam \eqref{vn_adam}, 
$\alpha_n := 10^{-2}$, $\beta_n = 10^{-3}$, $\hat{\beta} := 0.9$, and $\bar{\beta} := 0.999$
\end{description}
\item
Algorithm \ref{algo:1} with diminishing step-sizes
\begin{description}
\item[DSD: ]
Algorithm \ref{algo:1} with $\mathsf{h}_n^i$ defined by Stochastic Gradient Descent \cite{bon2013} (i.e., $\mathsf{h}_n^i := 1$), $\alpha_n := 10^{-1}/\sqrt{n}$, 
and $\beta_n = \hat{\beta} := 0$
\item[DAG: ]
Algorithm \ref{algo:1} with $\mathsf{h}_n^i$ defined by AdaGrad \cite{adagrad}, 
$\alpha_n := 10^{-1}/\sqrt{n}$, and $\beta_n = \hat{\beta} := 0$
\item[DAM1: ]
Algorithm \ref{algo:1} with $\mathsf{h}_n^i$ defined by AMSGrad \eqref{vn_amsgrad}, 
$\alpha_n := 10^{-1}/\sqrt{n}$, $\beta_n = 0.5^n$, $\hat{\beta} := 0$, and $\bar{\beta} := 0.999$
\item[DAM2: ]
Algorithm \ref{algo:1} with $\mathsf{h}_n^i$ defined by AMSGrad \eqref{vn_amsgrad}, 
$\alpha_n := 10^{-1}/\sqrt{n}$, $\beta_n = 0.9^n$, $\hat{\beta} := 0$, and $\bar{\beta} := 0.999$
\item[DAD1: ]
Algorithm \ref{algo:1} with $\mathsf{h}_n^i$ defined by Adam \eqref{vn_adam}, 
$\alpha_n := 10^{-1}/\sqrt{n}$, $\beta_n = 0.5^n$, $\hat{\beta} := 0.9$, and $\bar{\beta} := 0.999$
\item[DAD2: ]
Algorithm \ref{algo:1} with $\mathsf{h}_n^i$ defined by Adam \eqref{vn_adam}, 
$\alpha_n := 10^{-1}/\sqrt{n}$, $\beta_n = 0.9^n$, $\hat{\beta} := 0.9$, and $\bar{\beta} := 0.999$
\end{description}
\end{itemize}
The difference between CAM1 (resp. CAD1) and CAM2 (resp. CAD2) is the setting 
of $\beta_n$. 
The step-size $\beta_n = 0.9$ in CAM1 (resp. CAD1) is based on previously reported 
results (see, e.g., \cite[Section 5]{gec2019}),
while the step-size $\beta_n = 10^{-3}$ is based on Theorem \ref{cor:2_1}
indicating that a small step-size approximates a solution to Problem \ref{problem:1}.
The algorithms with diminishing step-sizes all used $\alpha_n := 10^{-1}/\sqrt{n}$, which is based on previously reported 
results (see, e.g., \cite[Theorems 1 and 2]{gec2019}). 

Ten samplings, each starting from a different randomly chosen
initial point $x_0(s) \in M$ ($s=1,2,\ldots,10$), were performed, and the results were averaged. 
The following two performance measures were used: for each $n\in\mathbb{N}$,
\begin{align*}
D_n := \frac{1}{10} \sum_{s=1}^{10} 
\sqrt{\sum_{i\in\mathcal{I}} \mathrm{d} \left(x_n^i(s), T^i(x_n^i(s)) \right)^2} \text{ and }
F_n := \frac{I}{10} \sum_{s=1}^{10} f(x_n(s)),
\end{align*}
where $(x_n^i(s))_{n\in\mathbb{N}}$ denotes the sequence generated by Algorithm \ref{algo:1} with 
an initial point $x_0(s)$.
If $(D_n)_{n\in\mathbb{N}}$ converges to $0$, then Algorithm \ref{algo:1} converges 
to a fixed point of $T^i$.

The experiments were conducted on a MacBook Air (2017) with a 1.8 GHz Intel Core i5 CPU, 8 GB
1600 MHz DDR3 memory, and the macOS Mojave version 10.14.5 operating system. The algorithms were written in Python 3.7.6 with the NumPy 1.19.2 package and the Matplotlib 3.1.2 package.
%Python implementations of the algorithms used in the numerical experiments are available at \url{https://github.com/iiduka-researches/202011-fixed-ropt}.

\subsection{Consistent case}\label{consistent}
We first consider the consistent case such that $\bigcap_{j \in \mathcal{J}^i} \mathrm{B}_j^i \neq \emptyset$ ($m=2,10,100; I = 5; J_i = 5$), where 
$c_j^i \in M$ and $r_j^i > 0$ in $\mathrm{B}_j^i$ defined by \eqref{ball} were randomly chosen.
%so that $\bigcap_{j \in \mathcal{J}^i} \mathrm{B}_j^i \neq \emptyset$ was satisfied. 
A nonexpansive mapping $T^i \colon D^{m} \to D^{m}$ ($i\in \mathcal{I}$) defined by \eqref{nonexp} satisfies 
$\mathrm{Fix}(T^i) = \bigcap_{j \in \mathcal{J}^i} \mathrm{B}_j^i$
(see also Proposition \ref{examples}(ii) and Example \ref{exp:1}).

Tables \ref{table:1} and \ref{table:2} show the average elapsed time (s) 
for the algorithms used in the experiment for $n = 500$ when $m=2$, $n = 1000$ when $m=10$, and $n=1500$ when $m=100$. 
The results  in these tables indicate that the elapsed times of the algorithms with constant step-sizes 
varied little from the elapsed times of the algorithms with diminishing step-sizes
and that, for a fixed $m$, all of the algorithms ran at about the same speed.

\begin{table}[htbp]
\centering
\caption{Average time for the algorithms with constant step-sizes applied to consistent case}
{\scriptsize
\begin{tabular}{l||cccccc}
\hline
                     &CSD & CAG & CAM1 & CAM2 & CAD1 & CAD2 \\ \hline
$m=2$      &7.728 & 7.782 & 8.128 & 8.115 & 7.892 & 7.878 \\ \hline
$m=10$     &16.219 & 16.381 & 16.974 & 16.534 & 16.546 & 17.077 \\ \hline
$m=100$    &23.683 & 23.788 & 23.907 & 24.347 & 24.536 & 24.187   \\   \hline     
\end{tabular}\label{table:1}
}
\end{table}

\begin{table}[htbp]
\centering
\caption{Average time for the algorithms with diminishing step-sizes applied to consistent case}
{\scriptsize
\begin{tabular}{l||cccccc}
\hline
                     &DSD & DAG & DAM1 & DAM2 & DAD1 & DAD2 \\ \hline
$m=2$      &7.862 & 7.906 & 8.382 & 8.155 & 7.989 & 8.364\\ \hline
$m=10$     &16.216 & 16.635 & 16.706 & 16.440 & 16.373 & 16.962 \\ \hline
$m=100$    &23.085 & 23.799 & 23.545 & 23.781 & 23.761 & 23.435  \\   \hline     
\end{tabular}\label{table:2}
}
\end{table}

Figures \ref{fig:1} and \ref{fig:2} show the behaviors of $D_n$ and $F_n$ for the algorithms with constant step-sizes, and Figures \ref{fig:3} and \ref{fig:4} 
show the behaviors of $D_n$ and $F_n$ for the algorithms with diminishing step-sizes.
The results in these figures indicate that all algorithms except for CSD, DSD, CAG, and DAG performed well. 
Although CAG and DAG converged to fixed points of $T^i$ faster than the other algorithms,
CAG and DAG did not minimize $f$. 
This is because CAG and DAG used $\beta_n = 0$ (i.e., $m_n = \mathsf{G}(x_n,\xi_n)$),
which means that CAG and DAG attached more weight to converging to a point in $X = \mathrm{Fix}(T^1) \times \mathrm{Fix}(T^2) \times \cdots \times \mathrm{Fix}(T^I)$
than minimizing $f$.
To verify why CSD and DSD did not converge to a fixed point of $T^i$, 
we checked the behaviors of CSD and DSD for ten samplings. 
CSD and DSD were sometimes good and sometimes not within ten samplings.
As a result, the mean value $D_n$ of 
$\sqrt{\sum_{i\in\mathcal{I}} \mathrm{d}(x_n^i(s),T^i (x_n^i(s)))^2}$
for CSD and DSD was not minimized.   

\begin{figure}[htbp]
\centering
\subfigure[$m=2$, $I=5$, $J_i=5$]{\includegraphics[width=39mm]{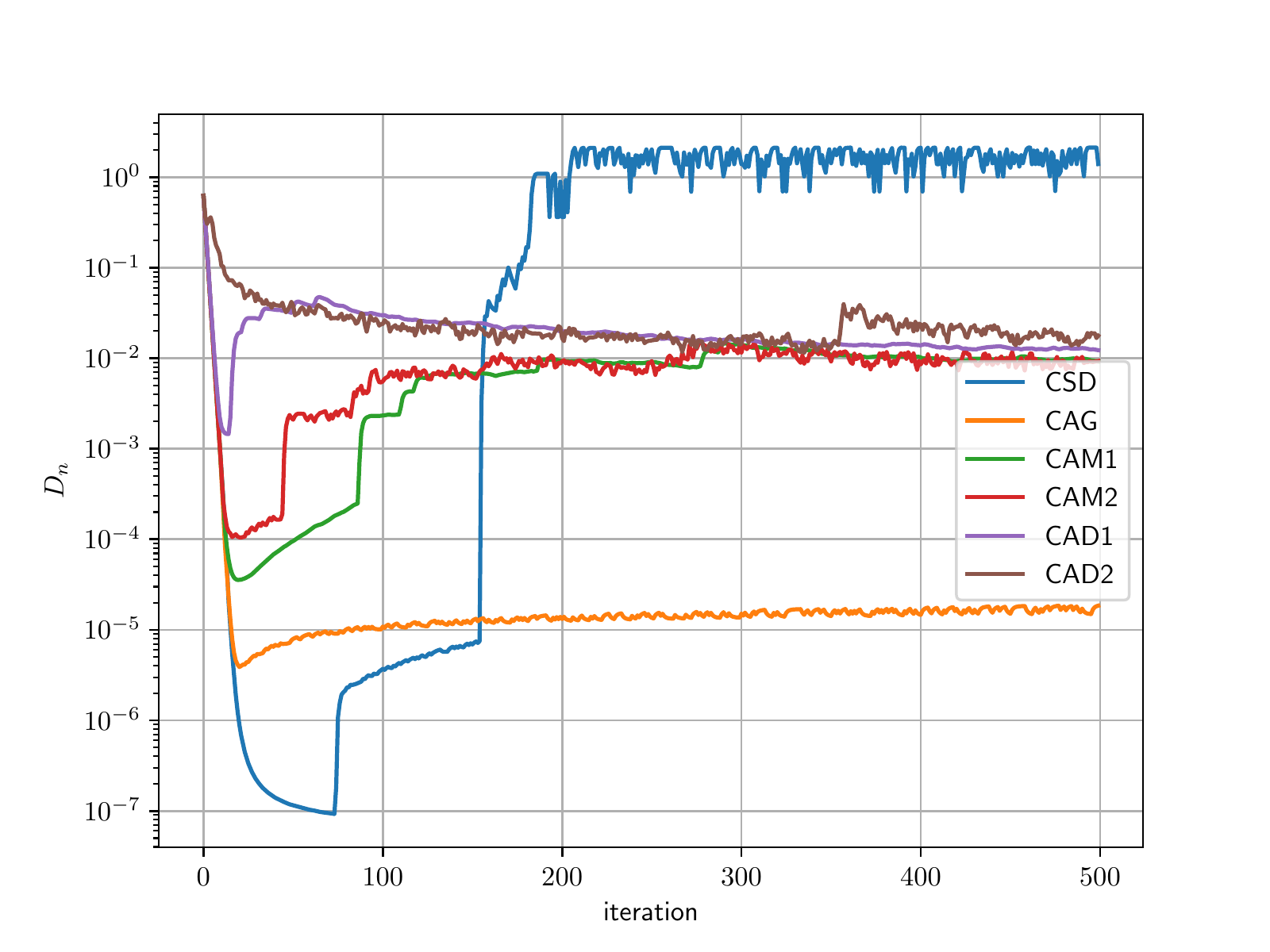}}
\subfigure[$m=10$, $I=5$, $J_i=5$]{\includegraphics[width=39mm]{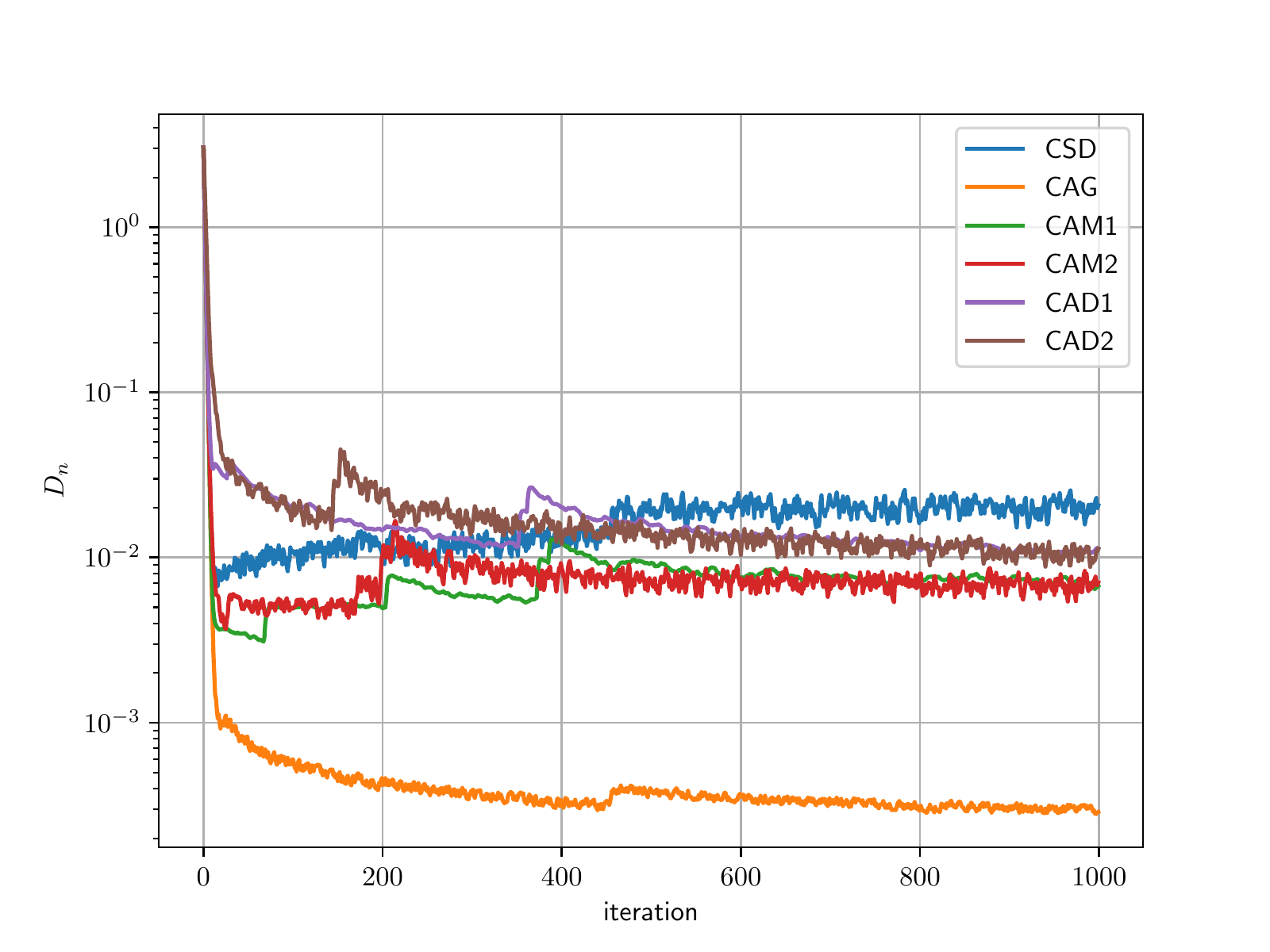}}
\subfigure[$m=100$, $I=5$, $J_i=5$]{\includegraphics[width=39mm]{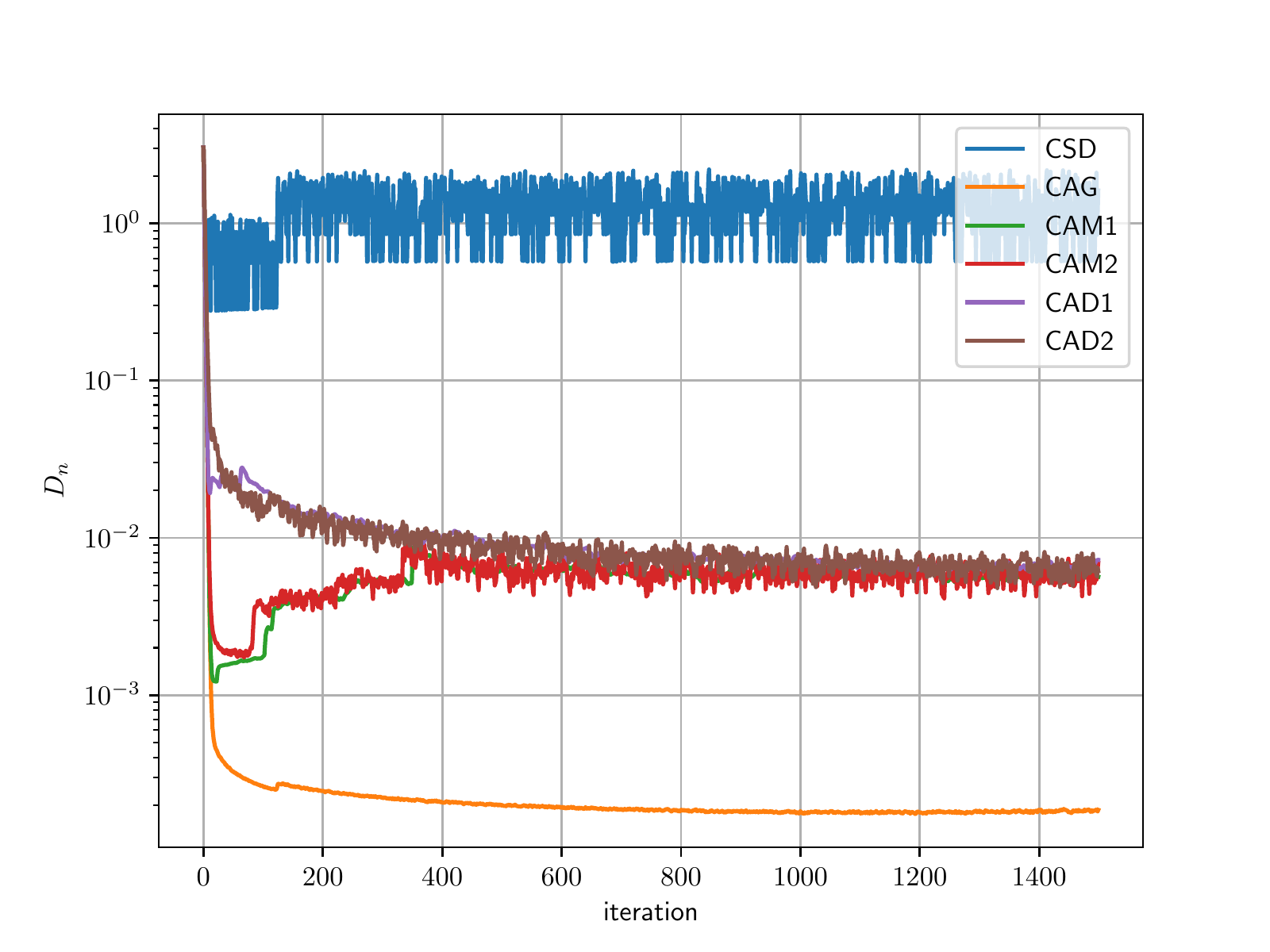}}
\caption{$D_n$ vs. iteration for Algorithm \ref{algo:1} with constant step-sizes (consistent case)}\label{fig:1}
\end{figure}

\begin{figure}[htbp]
\centering
\subfigure[$m=2$, $I=5$, $J_i=5$]{\includegraphics[width=39mm]{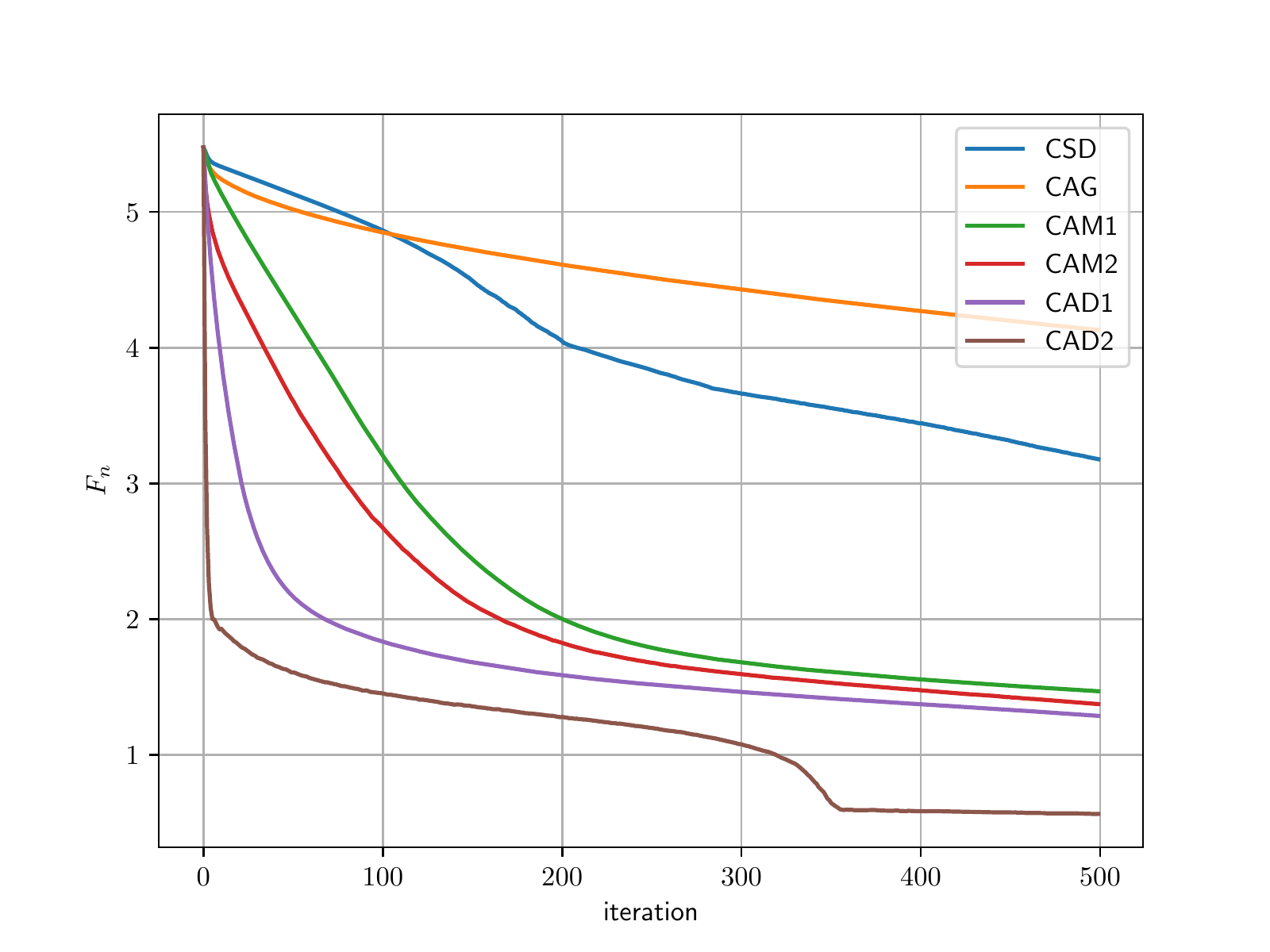}}
\subfigure[$m=10$, $I=5$, $J_i=5$]{\includegraphics[width=39mm]{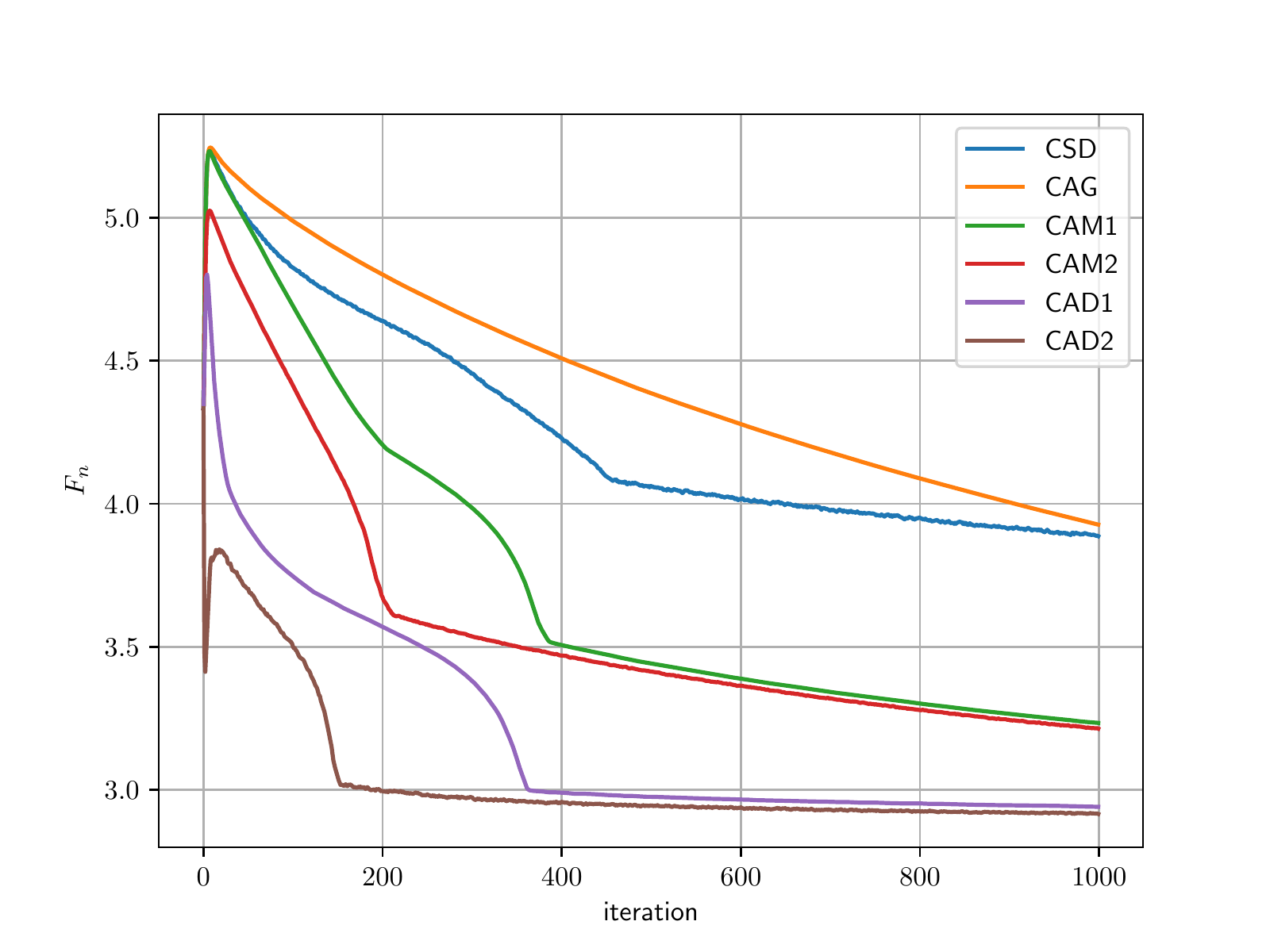}}
\subfigure[$m=100$, $I=5$, $J_i=5$]{\includegraphics[width=39mm]{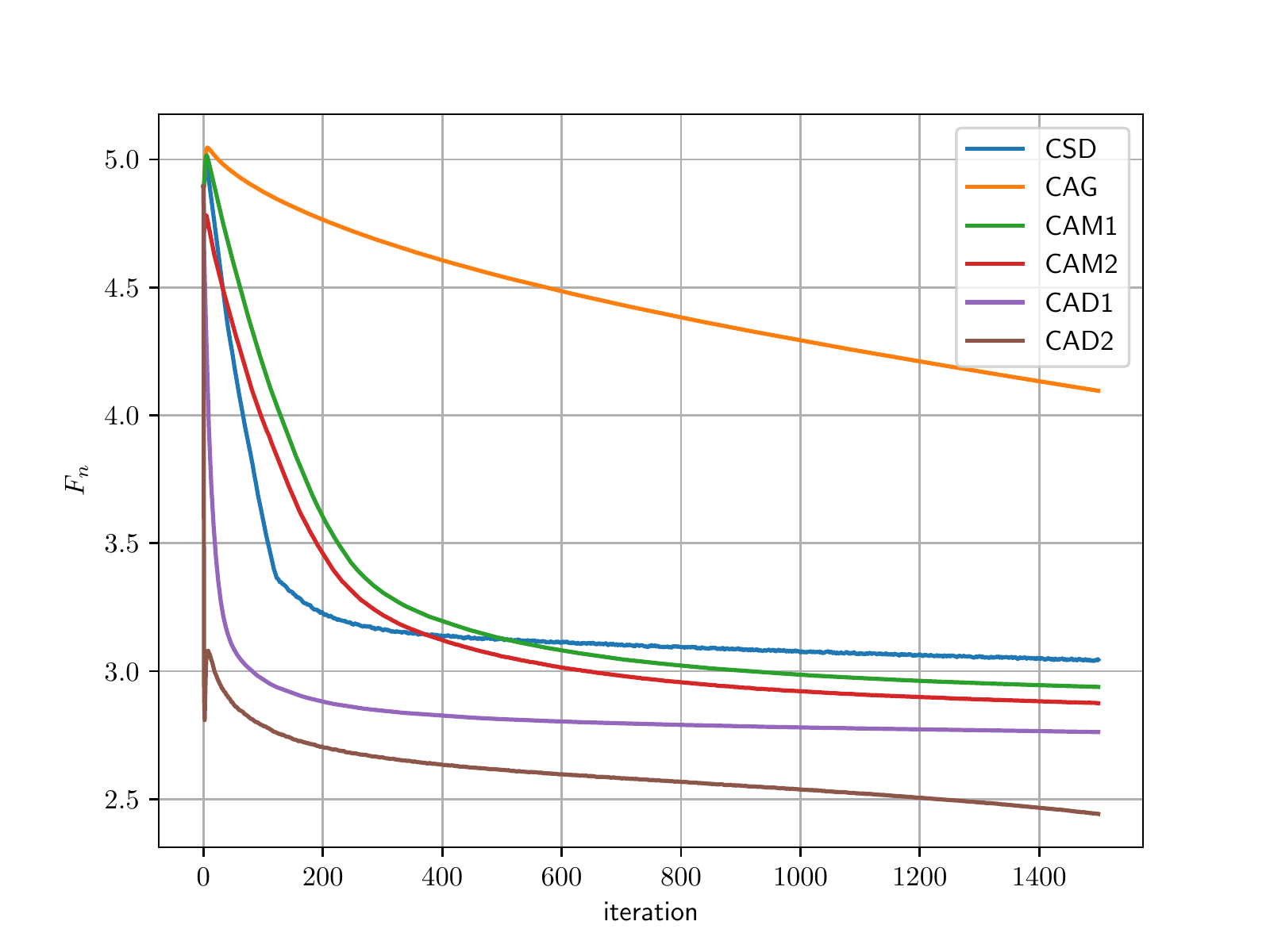}}
\caption{$F_n$ vs. iteration for Algorithm \ref{algo:1} with constant step-sizes (consistent case)}\label{fig:2}
\end{figure}

\begin{figure}[htbp]
\centering
\subfigure[$m=2$, $I=5$, $J_i=5$]{\includegraphics[width=39mm]{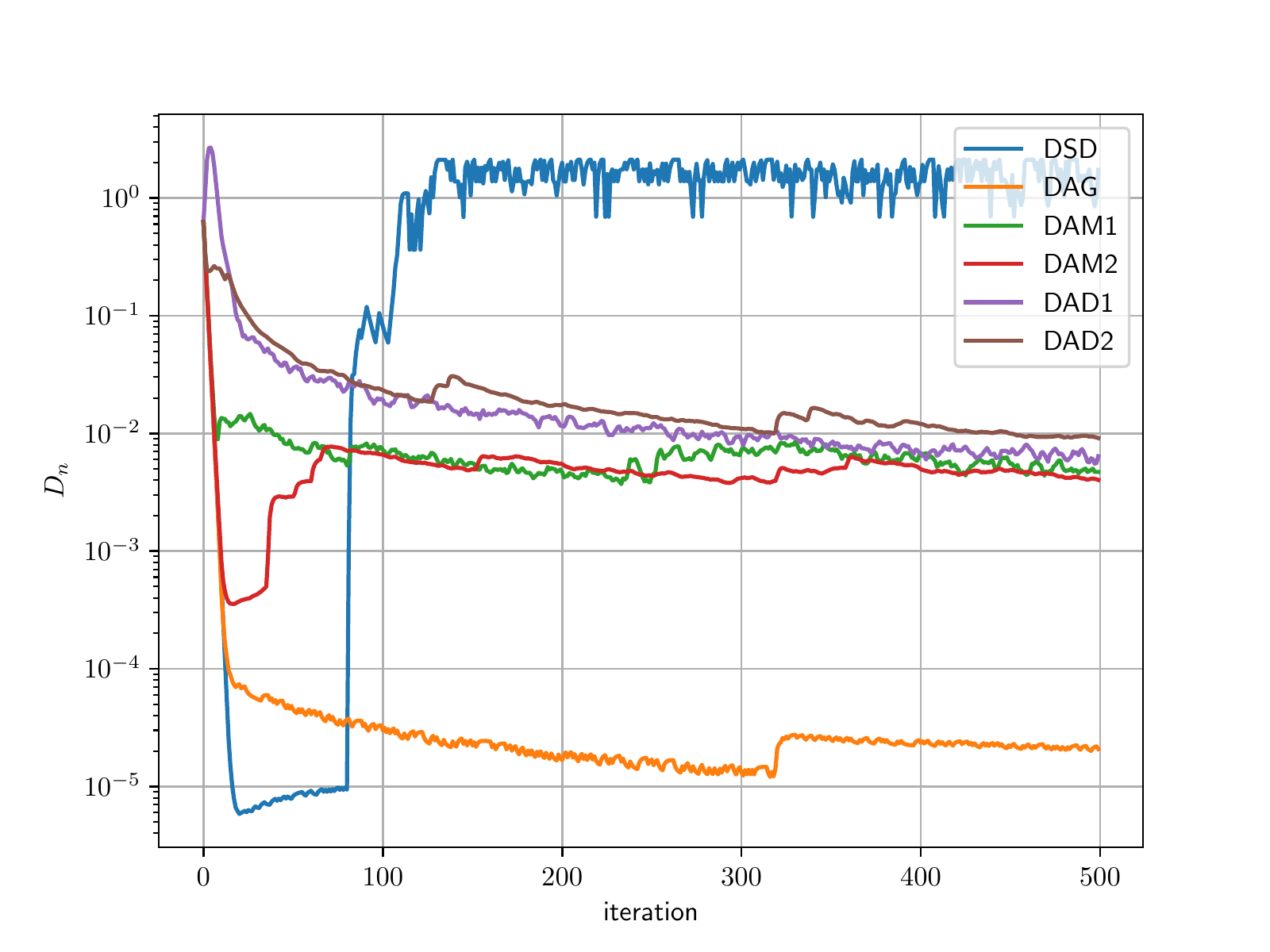}}
\subfigure[$m=10$, $I=5$, $J_i=5$]{\includegraphics[width=39mm]{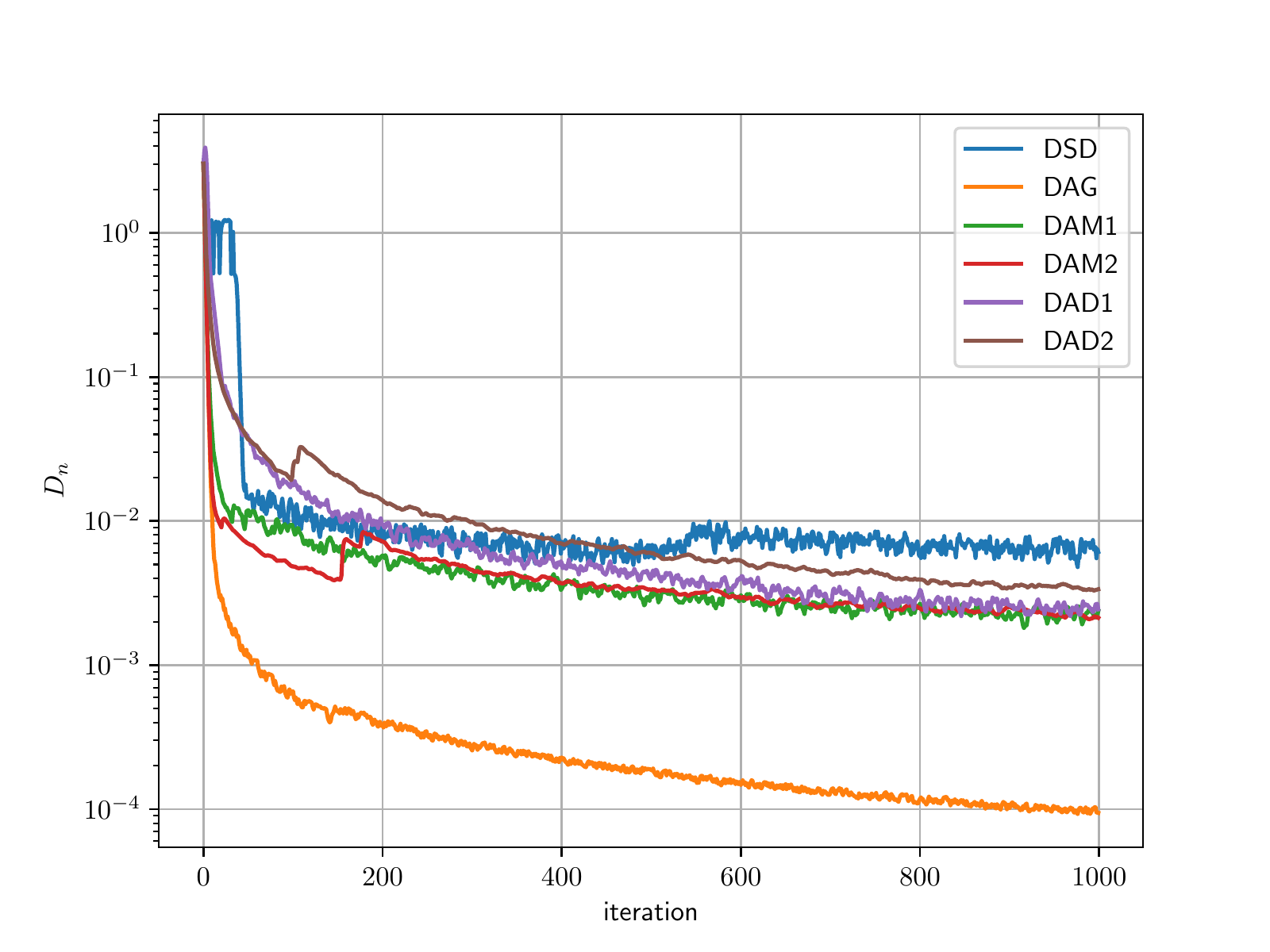}}
\subfigure[$m=100$, $I=5$, $J_i=5$]{\includegraphics[width=39mm]{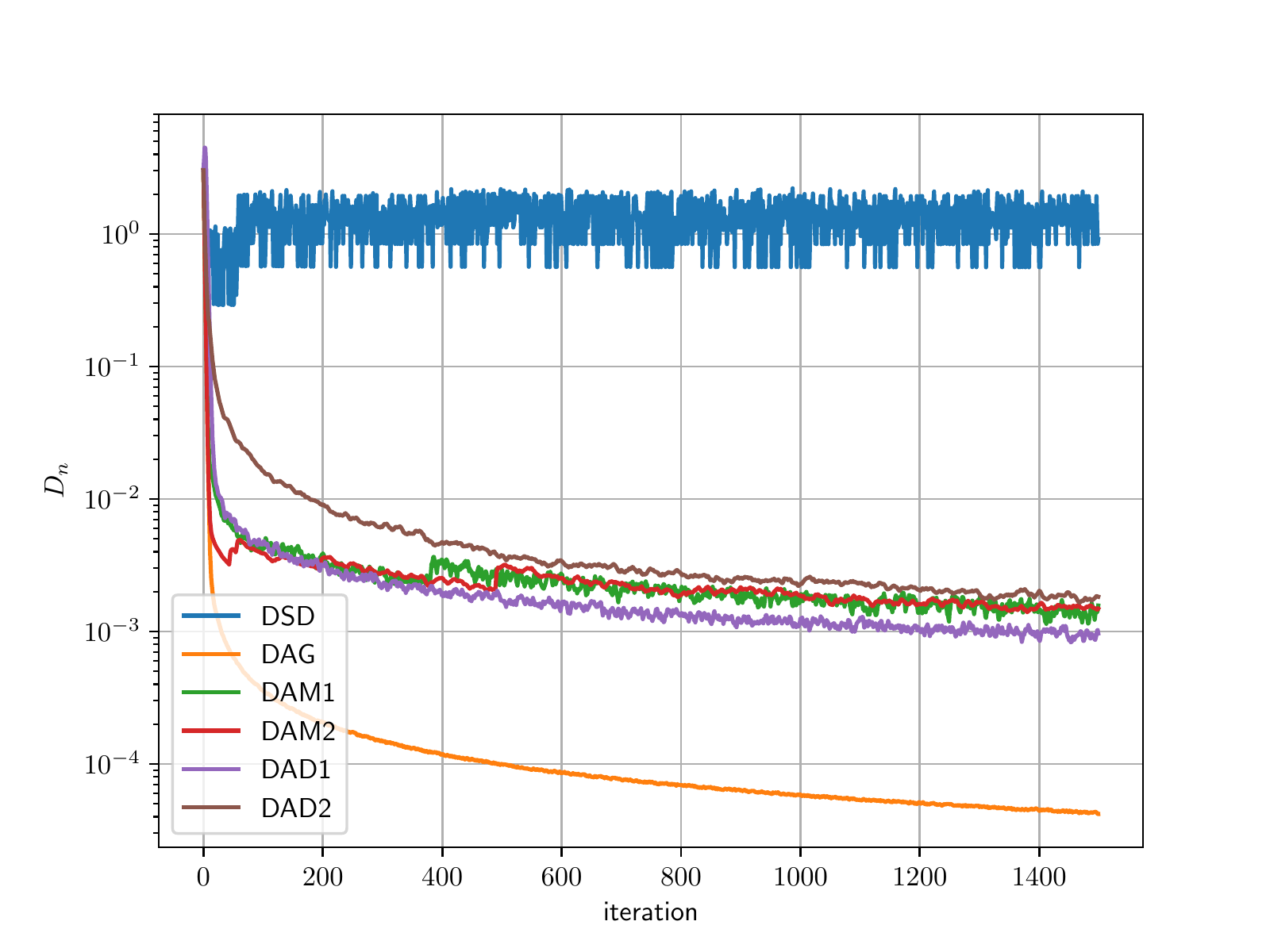}}
\caption{$D_n$ vs. iteration for Algorithm \ref{algo:1} with diminishing step-sizes (consistent case)}\label{fig:3}
\end{figure}

\begin{figure}[htbp]
\centering
\subfigure[$m=2$, $I=5$, $J_i=5$]{\includegraphics[width=39mm]{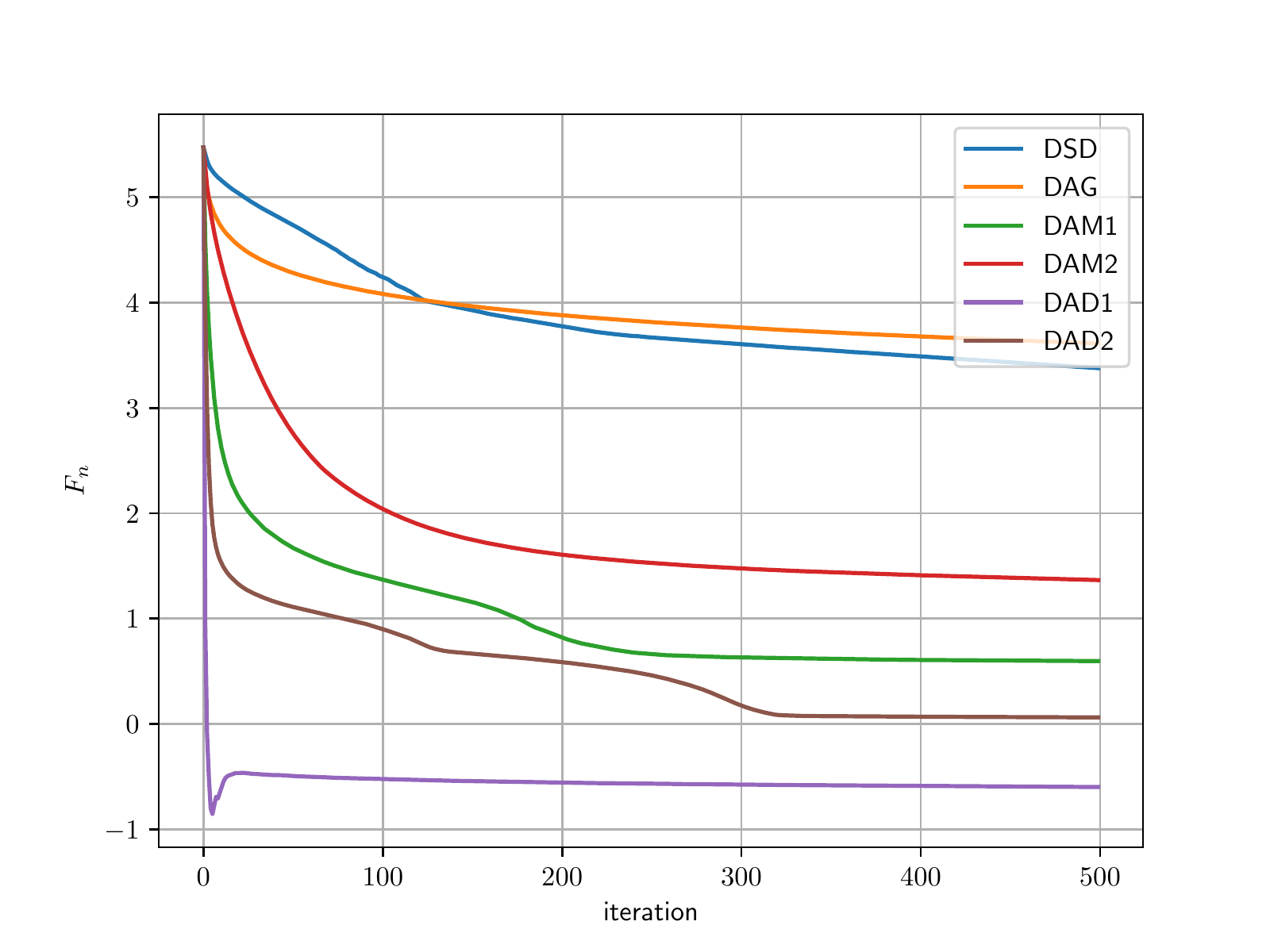}}
\subfigure[$m=10$, $I=5$, $J_i=5$]{\includegraphics[width=39mm]{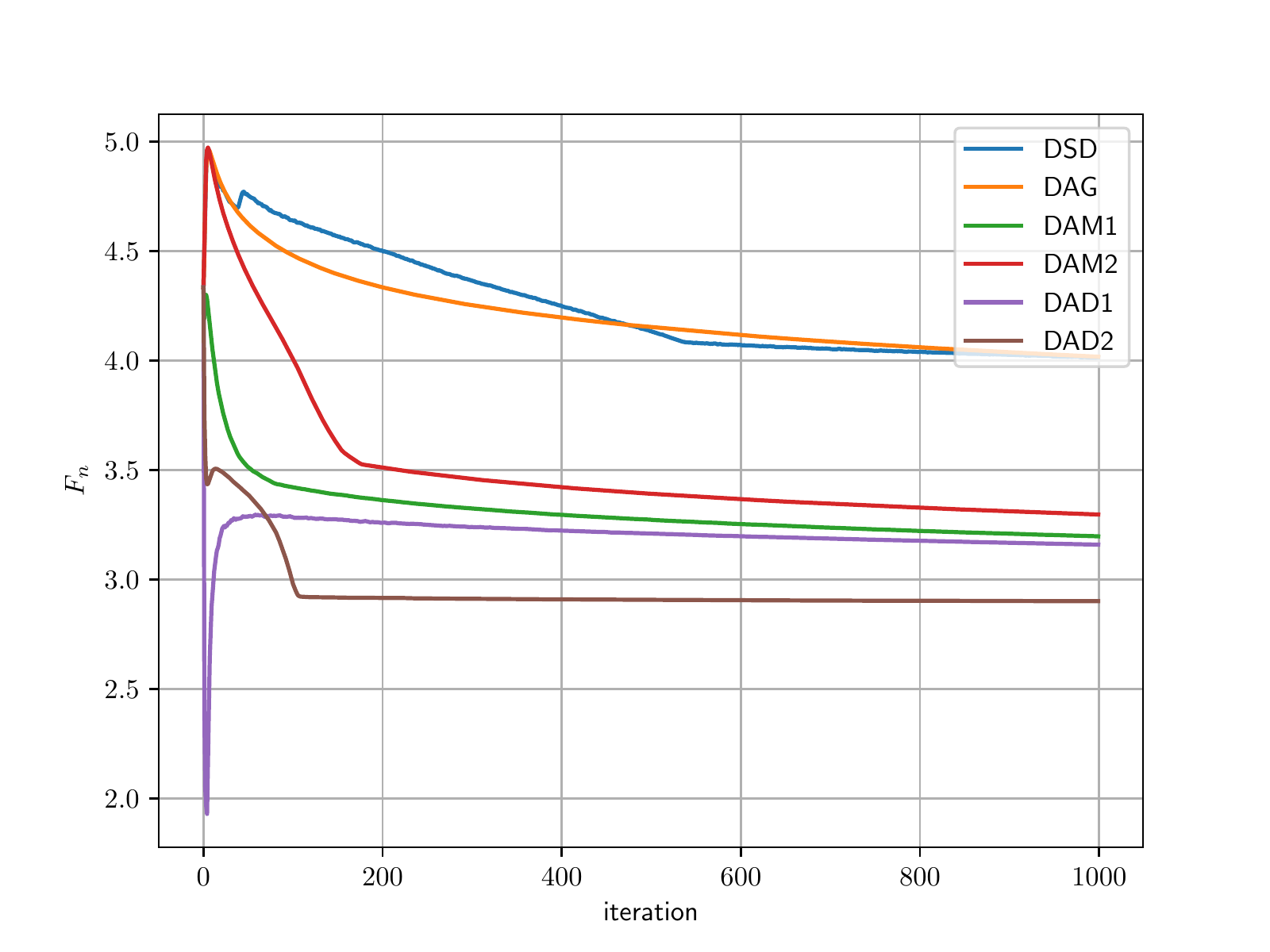}}
\subfigure[$m=100$, $I=5$, $J_i=5$]{\includegraphics[width=39mm]{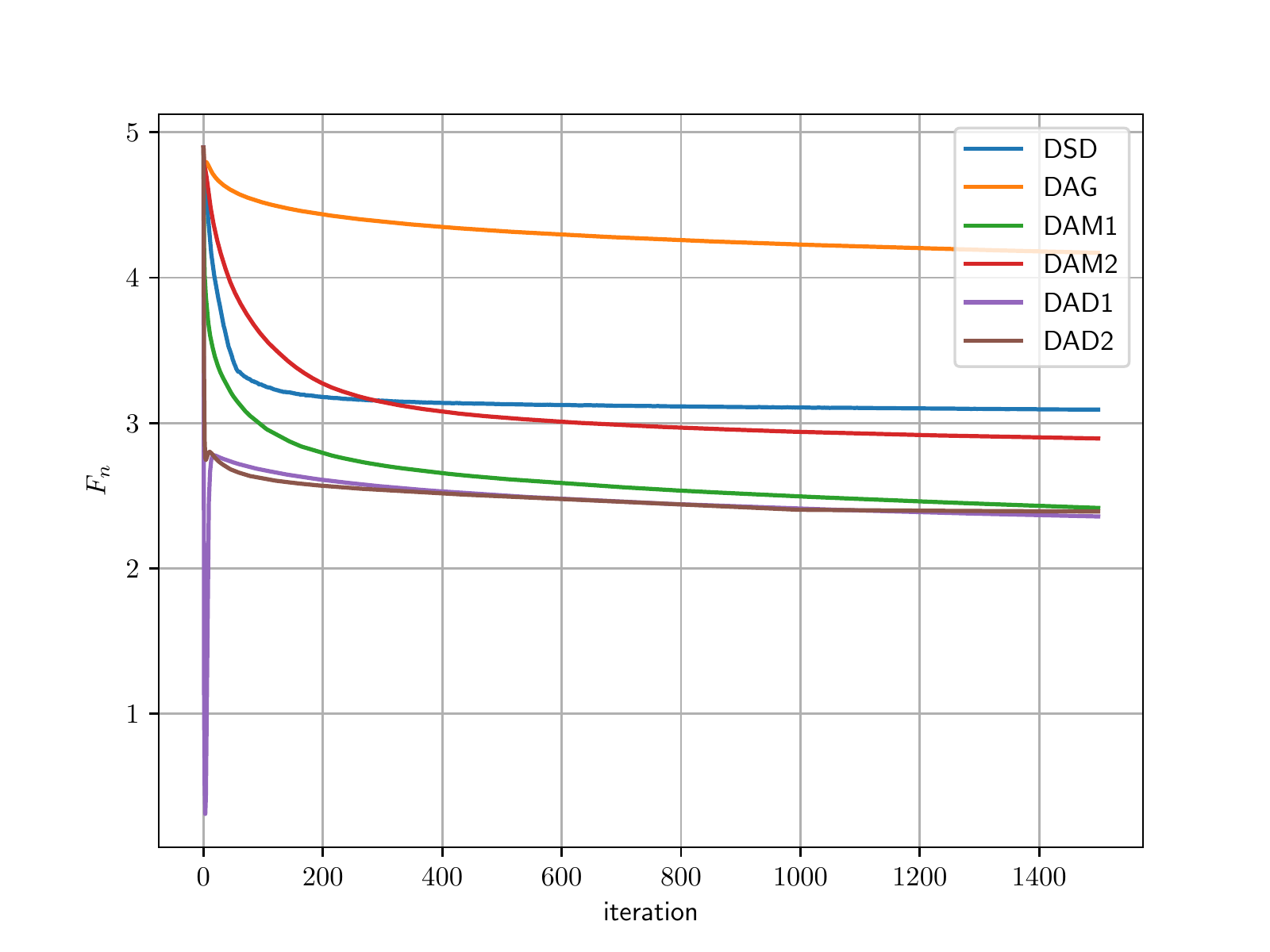}}
\caption{$F_n$ vs. iteration for Algorithm \ref{algo:1} with diminishing step-sizes (consistent case)}\label{fig:4}
\end{figure}

\subsection{Inconsistent case}\label{inconsistent}
We next consider the inconsistent case such that $\bigcap_{j \in \mathcal{J}^i} \mathrm{B}_j^i = \emptyset$, where 
$c_j^i \in M$ and $r_j^i > 0$ in $\mathrm{B}_j^i$ $(m=2,10,100; I= 5; \mathcal{J}^i = \{1,2\})$ defined by \eqref{ball} were randomly chosen
so that $\bigcap_{j \in \mathcal{J}^i} \mathrm{B}_j^i = \emptyset$ was satisfied.
Here, we define a {\em generalized convex feasible set} (see \cite[Section I, Framework 2]{com1999} and \cite[Definition 4.1]{yamada} for the definition under the Hilbert space setting) as follows:
\begin{align}\label{gcfs}
C_\mathrm{d}^i := \left\{ x \in \mathrm{B}_1^i \colon 
\mathrm{d}\left(x, \mathrm{B}_2^i \right)^2 = \inf_{y\in \mathrm{B}_1^i} \mathrm{d}\left(y, \mathrm{B}_2^i \right)^2   \right\}.
\end{align}
The generalized convex feasible set plays an important
role when the constraint set composed of the absolute set and the subsidiary
set is not feasible.
Let $\mathrm{B}_1^i$ be the absolute constrained set and $\mathrm{B}_2^i$ be the subsidiary constrained set.
Then, $C_\mathrm{d}^i$ is feasible (i.e., $C_\mathrm{d}^i \neq \emptyset$) 
even when $\mathrm{B}_1^i \cap \mathrm{B}_2^i = \emptyset$.
Moreover, $C_\mathrm{d}^i$ is a subset of the absolute constrained set $\mathrm{B}_1^i$ with the elements closest to the subsidiary constrained set $\mathrm{B}_2^i$ in terms of the distance function.
Accordingly, it would be reasonable to replace an inconsistent set with the generalized convex feasible set. 
The set $C_\mathrm{d}^i$ defined by \eqref{gcfs} can be expressed as follows:
\begin{align*}
C_\mathrm{d}^i 
= \mathrm{Fix}
\left(P_1^i \left(  \exp \left[- \mathrm{grad} \text{ } 
  \frac{1}{2} \mathrm{d}\left(\cdot, \mathrm{B}_2^i \right)^2 \right]  \right) \right)
= \mathrm{Fix} \left(P_1^i P_2^i \right)
= \mathrm{Fix} \left(T^i \right),
\end{align*}
where the first equation comes from 
\cite[Proposition 3.1, Corollaries 3.1 and 3.2, Theorem 3.3]{li2009_1} 
(see also Proposition \ref{vi}),
the second equation comes from $\mathrm{grad}(1/2)\mathrm{d}(x,y)^2 = - \exp_x^{-1}(y)$ \cite[Proposition 3.3]{ferr2002}, and the third equation comes from \eqref{nonexp}.

Tables \ref{table:3} and \ref{table:4}
%show the average elapsed time (s) 
%for the algorithms used in the experiment under $n = 500$ when $m=2$, $n = 1000$ when $m=10$, and $n=1500$ when $m=100$. 
%These tables 
show that the elapsed times of the algorithms with constant step-sizes 
differed little from the elapsed times of the algorithms with diminishing step-sizes
and that, for a fixed $m$, all of the algorithms ran at about the same speed.

\begin{table}[htbp]
\centering
\caption{Average time for the algorithms with constant step-sizes applied to inconsistent case}
{\scriptsize
\begin{tabular}{l||cccccc}
\hline
                     &CSD & CAG & CAM1 & CAM2 & CAD1 & CAD2 \\ \hline
$m=2$      &3.906 & 3.878 & 4.010 & 4.014 & 3.940 & 3.933 \\ \hline
$m=10$     &8.649 & 8.675 & 8.822 & 9.066 & 8.727 & 8.691 \\ \hline
$m=100$    &13.727 & 13.831 & 14.303 & 14.021 & 14.003 & 14.061   \\   \hline     
\end{tabular}\label{table:3}
}
\end{table}

\begin{table}[htbp]
\centering
\caption{Average time for the algorithms with diminishing step-sizes applied to inconsistent case}
{\scriptsize
\begin{tabular}{l||cccccc}
\hline
                     &DSD & DAG & DAM1 & DAM2 & DAD1 & DAD2 \\ \hline
$m=2$      &3.883 & 3.880 & 4.013 & 4.075 & 4.010 & 3.946\\ \hline
$m=10$     &8.783 & 8.739 & 9.043 & 9.105 & 8.912 & 8.910 \\ \hline
$m=100$    &13.571 & 13.908 & 13.971 & 13.983 & 13.956 & 13.970  \\   \hline     
\end{tabular}\label{table:4}
}
\end{table}

Figures \ref{fig:5} and \ref{fig:6} show the behaviors of $D_n$ and $F_n$ for the algorithms with constant step-sizes, and Figures \ref{fig:7} and \ref{fig:8} 
show the behaviors of $D_n$ and $F_n$ for the algorithms with diminishing step-sizes.
The results shown in these figures indicate that all algorithms except for CSD, DSD, CAG, and DAG performed well, the same as in the consistent case (previous subsection).We checked the behaviors of CSD, DSD, CAG, and DAG and found that the reason they did not perform well was the same as in that case.
\begin{figure}[htbp]
\centering
\subfigure[$m=2$, $I=5$, $J_i=2$]{\includegraphics[width=39mm]{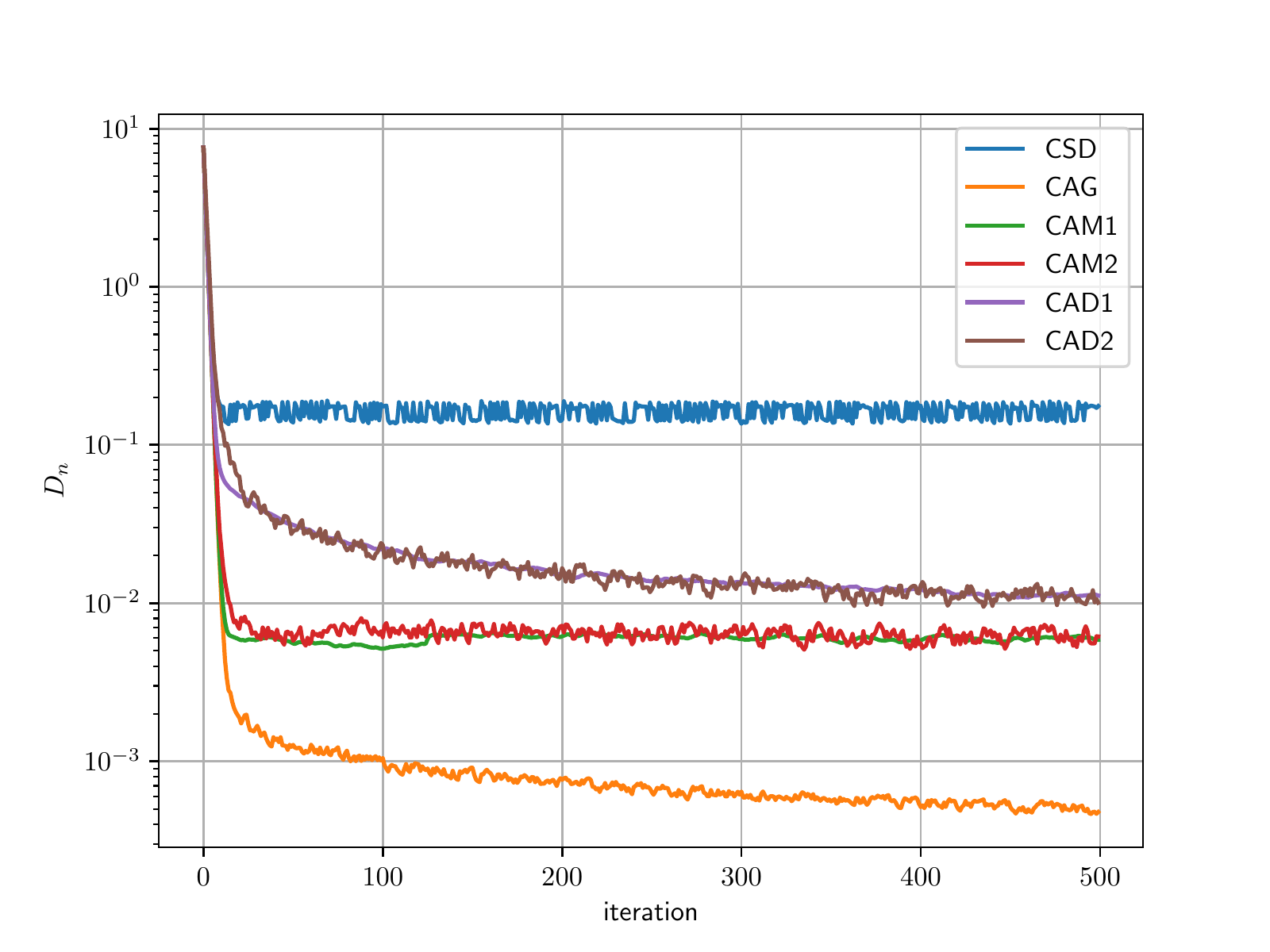}}
\subfigure[$m=10$, $I=5$, $J_i=2$]{\includegraphics[width=39mm]{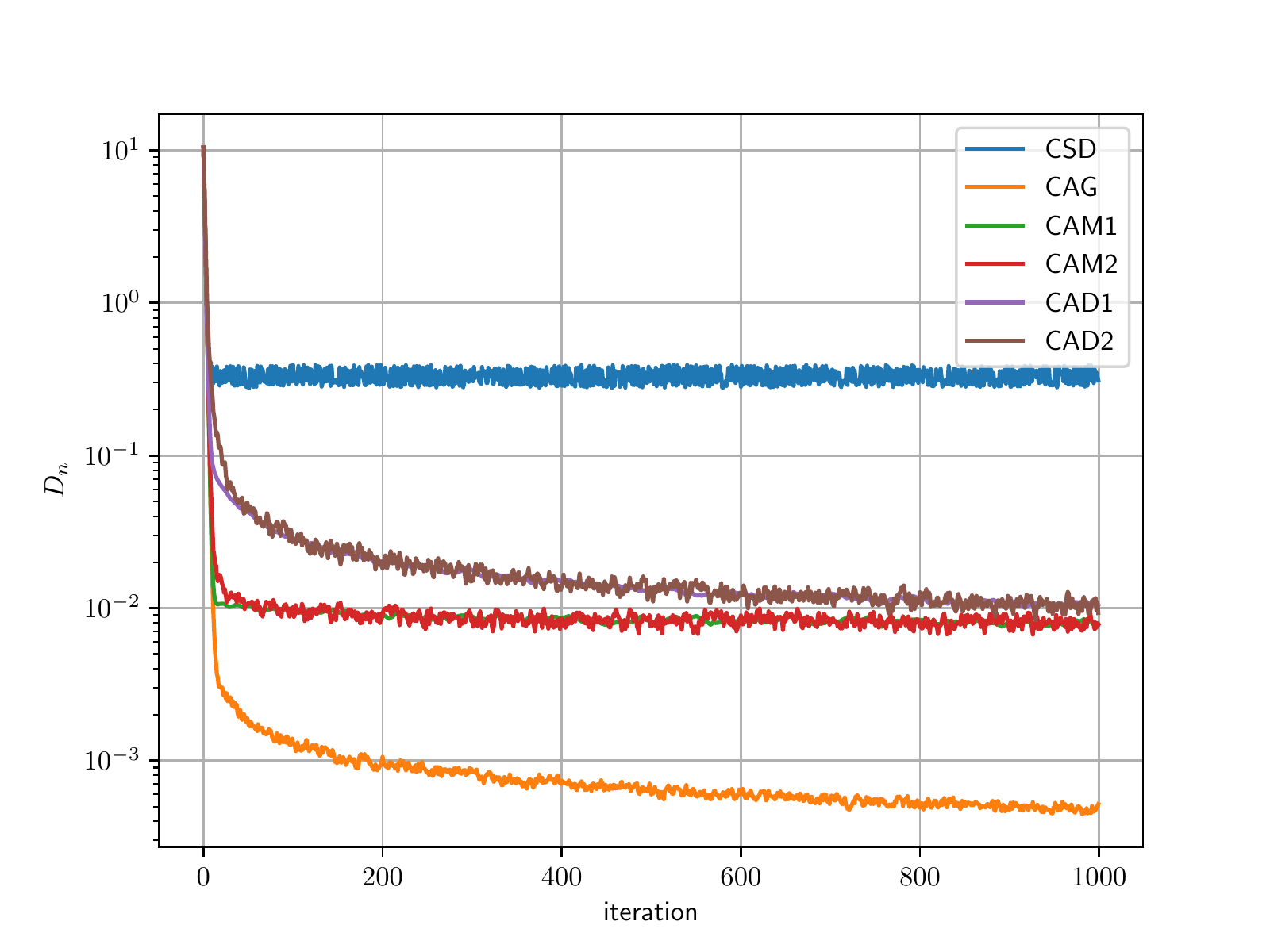}}
\subfigure[$m=100$, $I=5$, $J_i=2$]{\includegraphics[width=39mm]{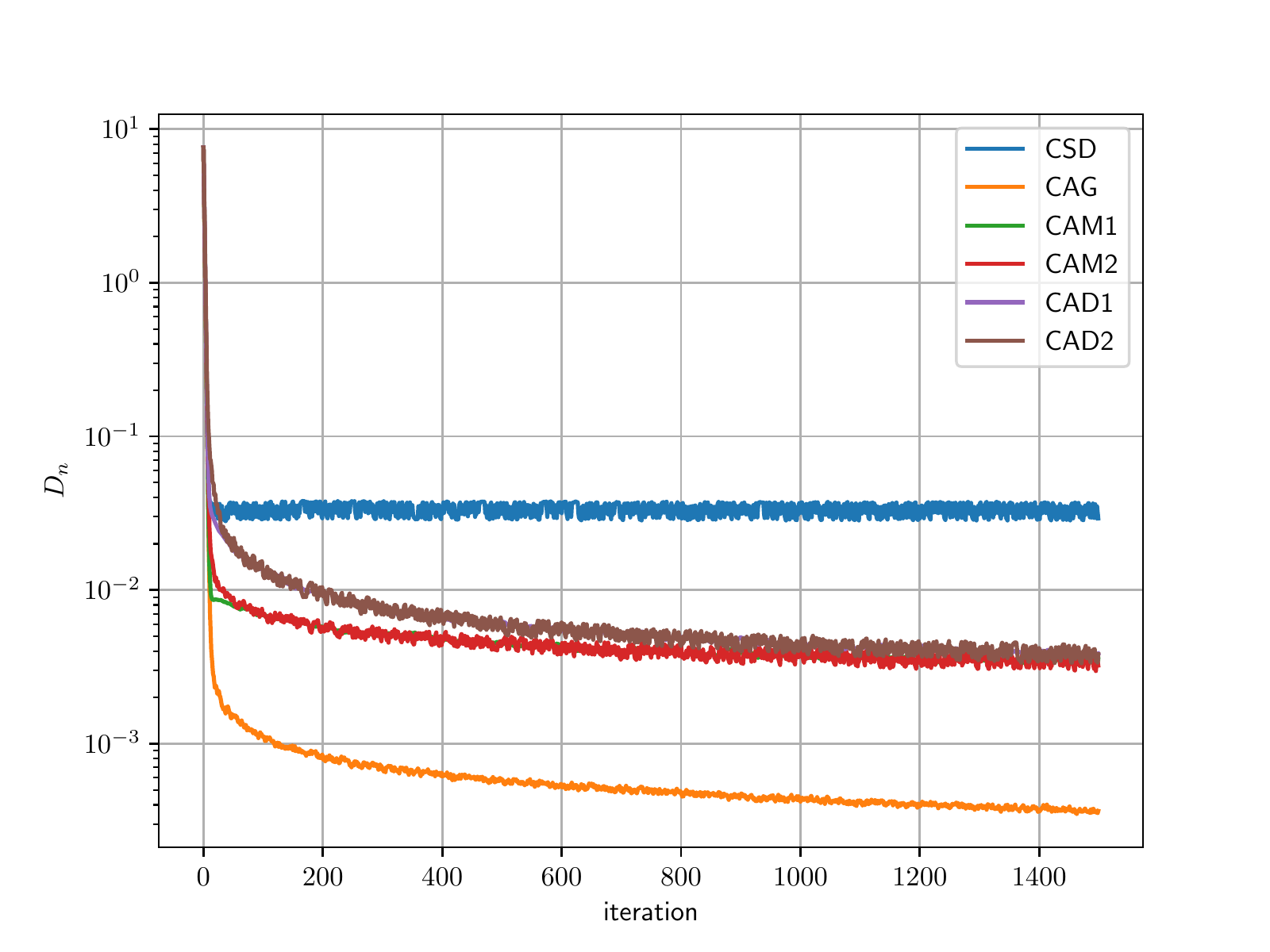}}
\caption{$D_n$ vs. iteration for Algorithm \ref{algo:1} with constant step-sizes (inconsistent case)}\label{fig:5}
\end{figure}

\begin{figure}[htbp]
\centering
\subfigure[$m=2$, $I=5$, $J_i=2$]{\includegraphics[width=39mm]{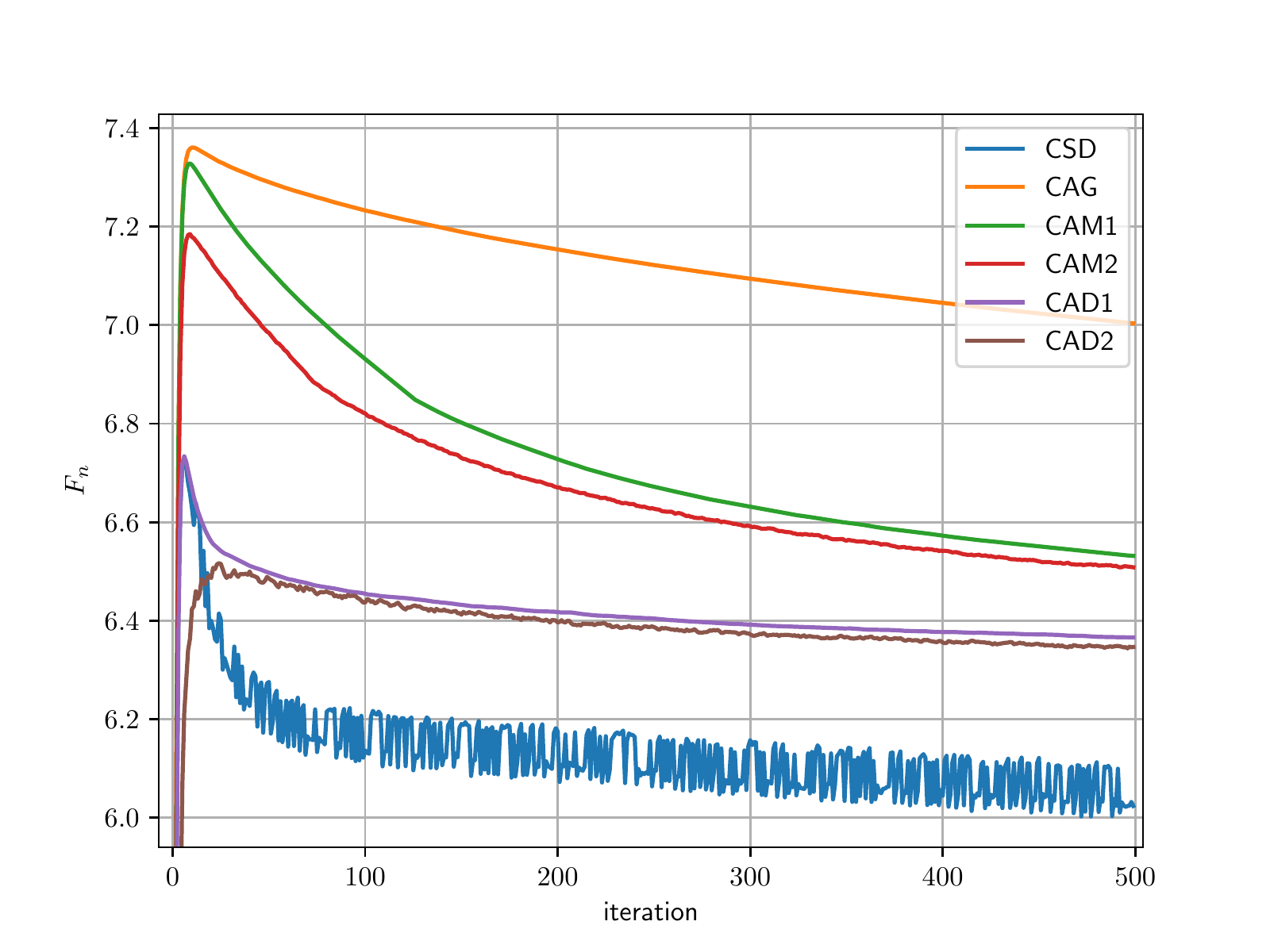}}
\subfigure[$m=10$, $I=5$, $J_i=2$]{\includegraphics[width=39mm]{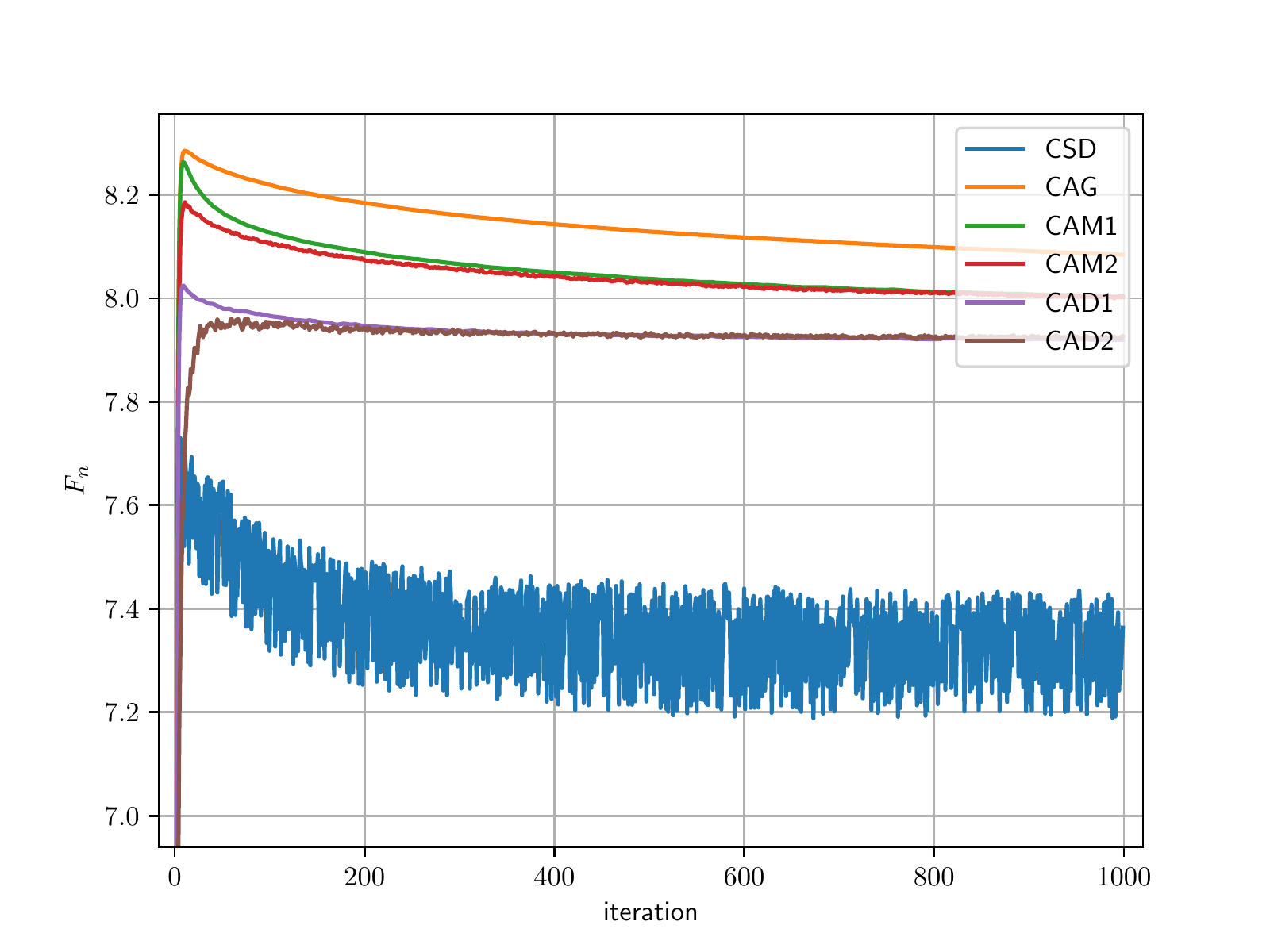}}
\subfigure[$m=100$, $I=5$, $J_i=2$]{\includegraphics[width=39mm]{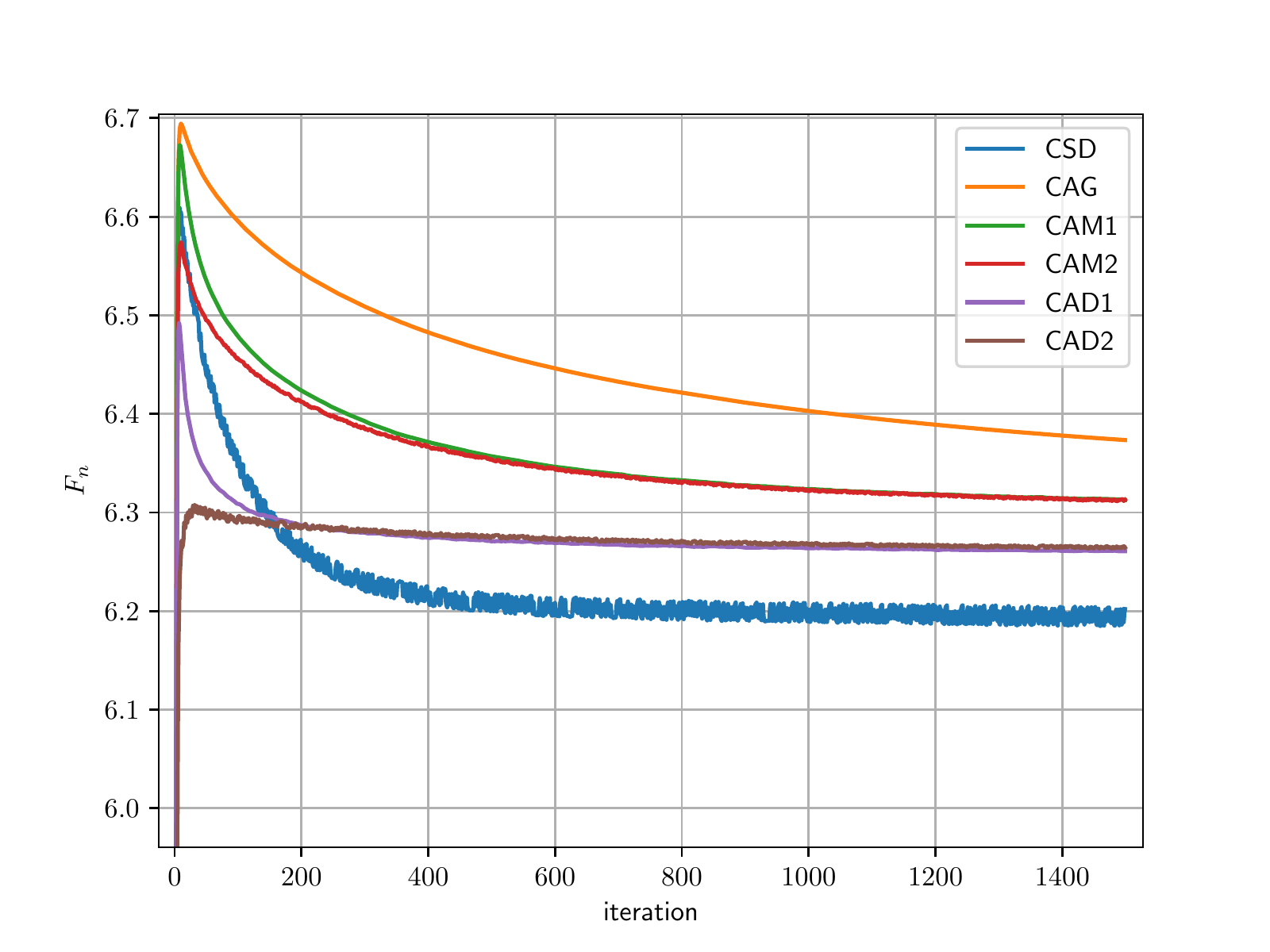}}
\caption{$F_n$ vs. iteration for Algorithm \ref{algo:1} with constant step-sizes (inconsistent case)}\label{fig:6}
\end{figure}

\begin{figure}[htbp]
\centering
\subfigure[$m=2$, $I=5$, $J_i=2$]{\includegraphics[width=39mm]{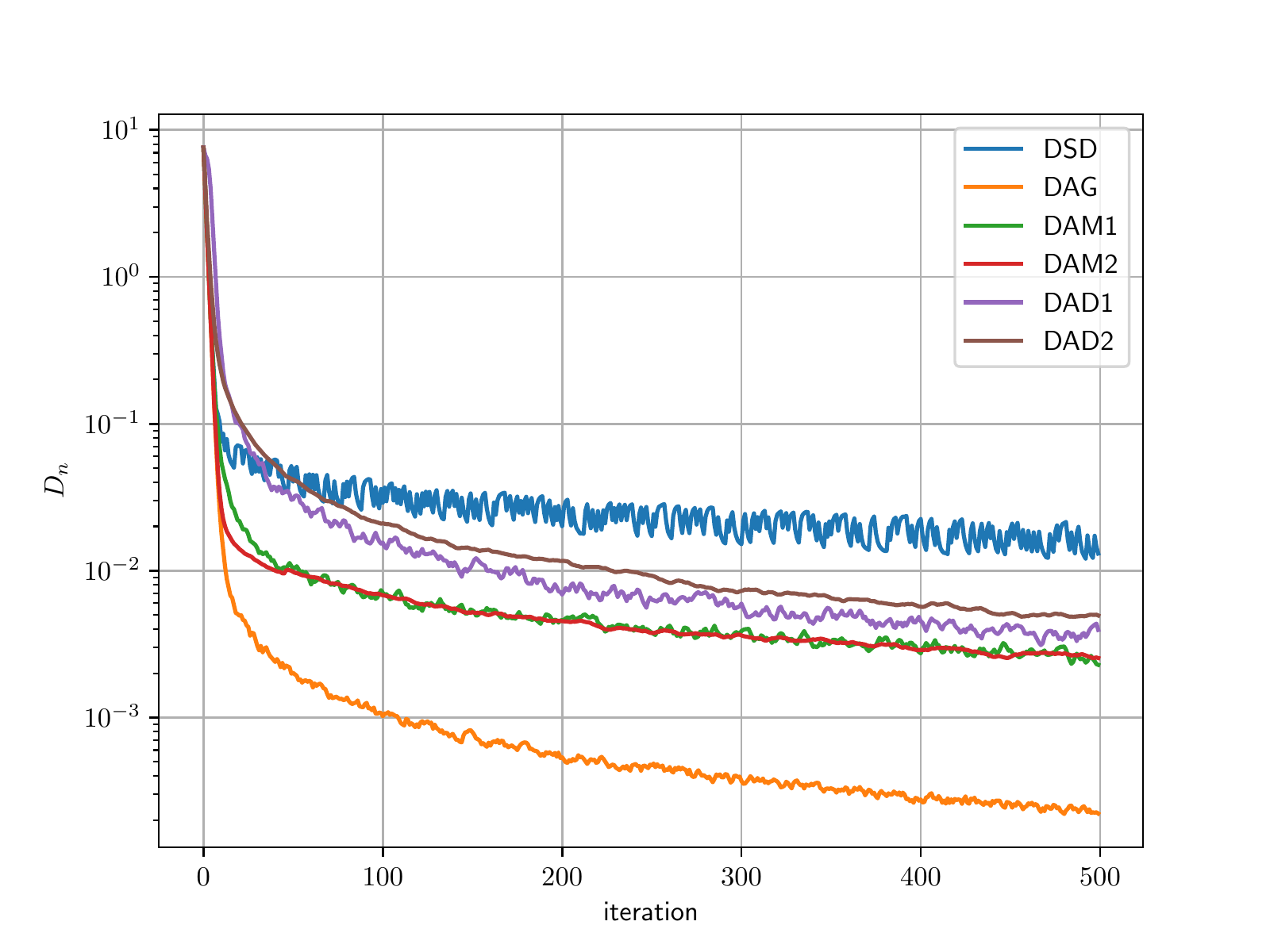}}
\subfigure[$m=10$, $I=5$, $J_i=2$]{\includegraphics[width=39mm]{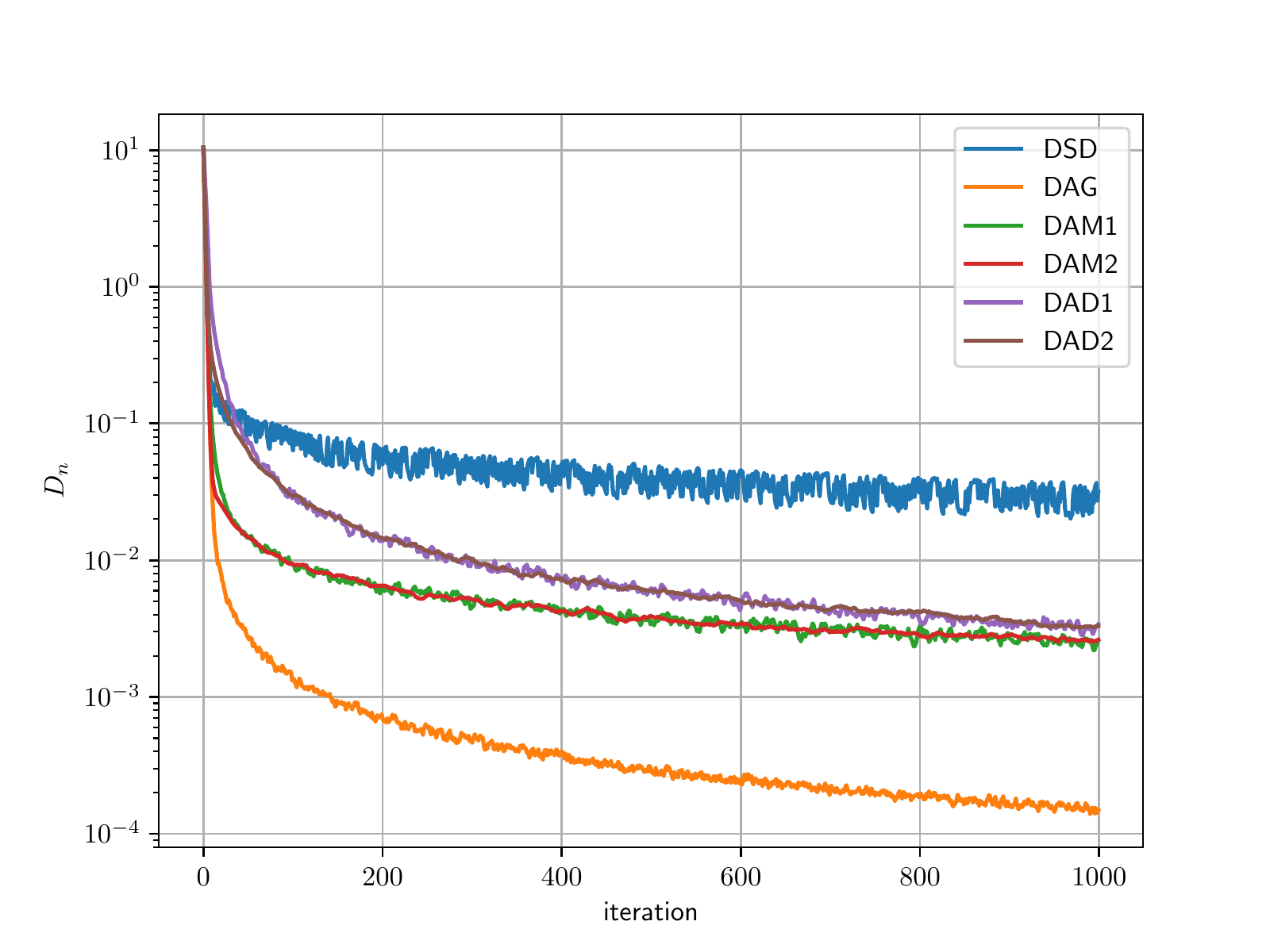}}
\subfigure[$m=100$, $I=5$, $J_i=2$]{\includegraphics[width=39mm]{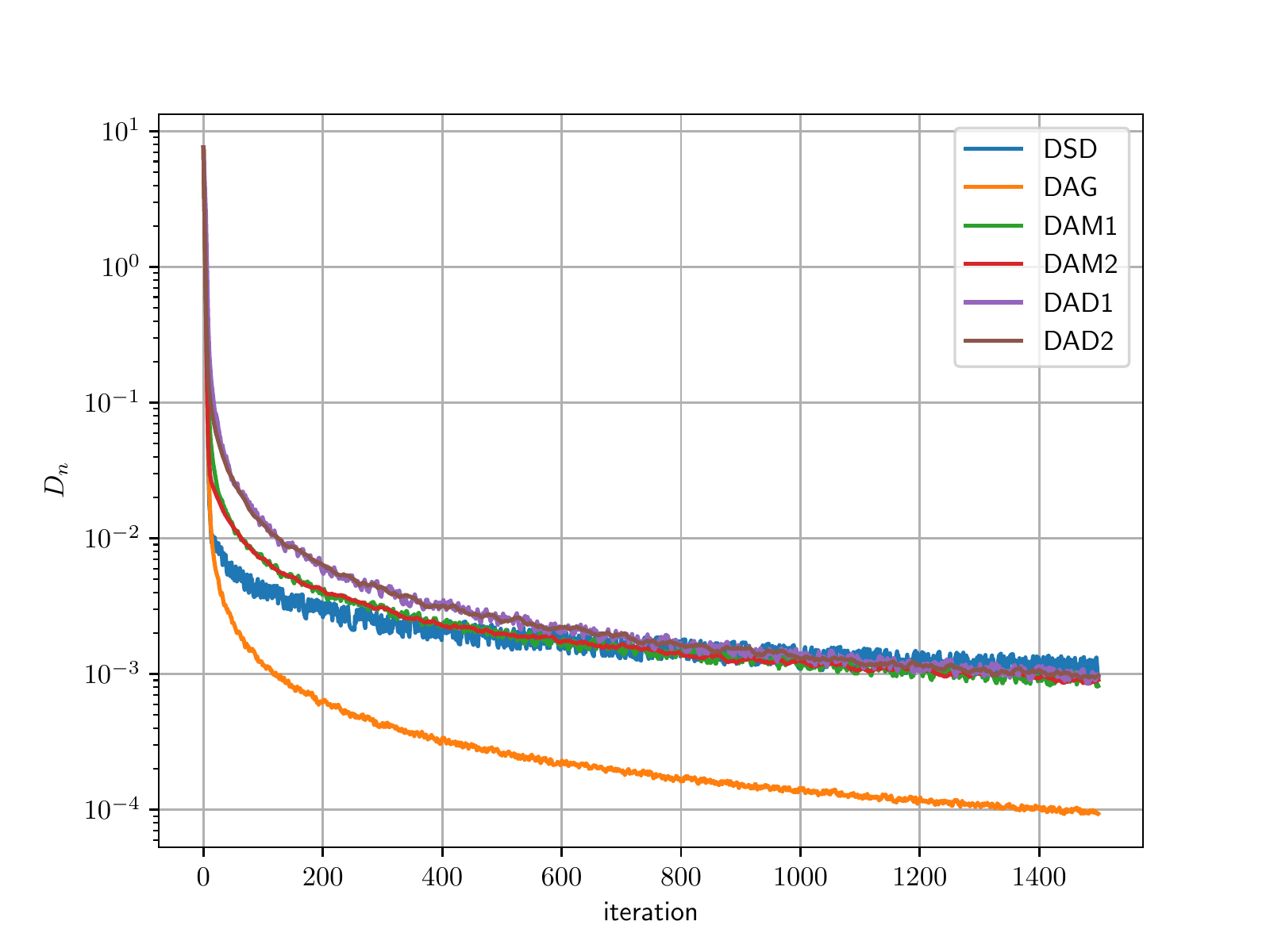}}
\caption{$D_n$ vs. iteration for Algorithm \ref{algo:1} with diminishing step-sizes (inconsistent case)}\label{fig:7}
\end{figure}

\begin{figure}[htbp]
\centering
\subfigure[$m=2$, $I=5$, $J_i=2$]{\includegraphics[width=39mm]{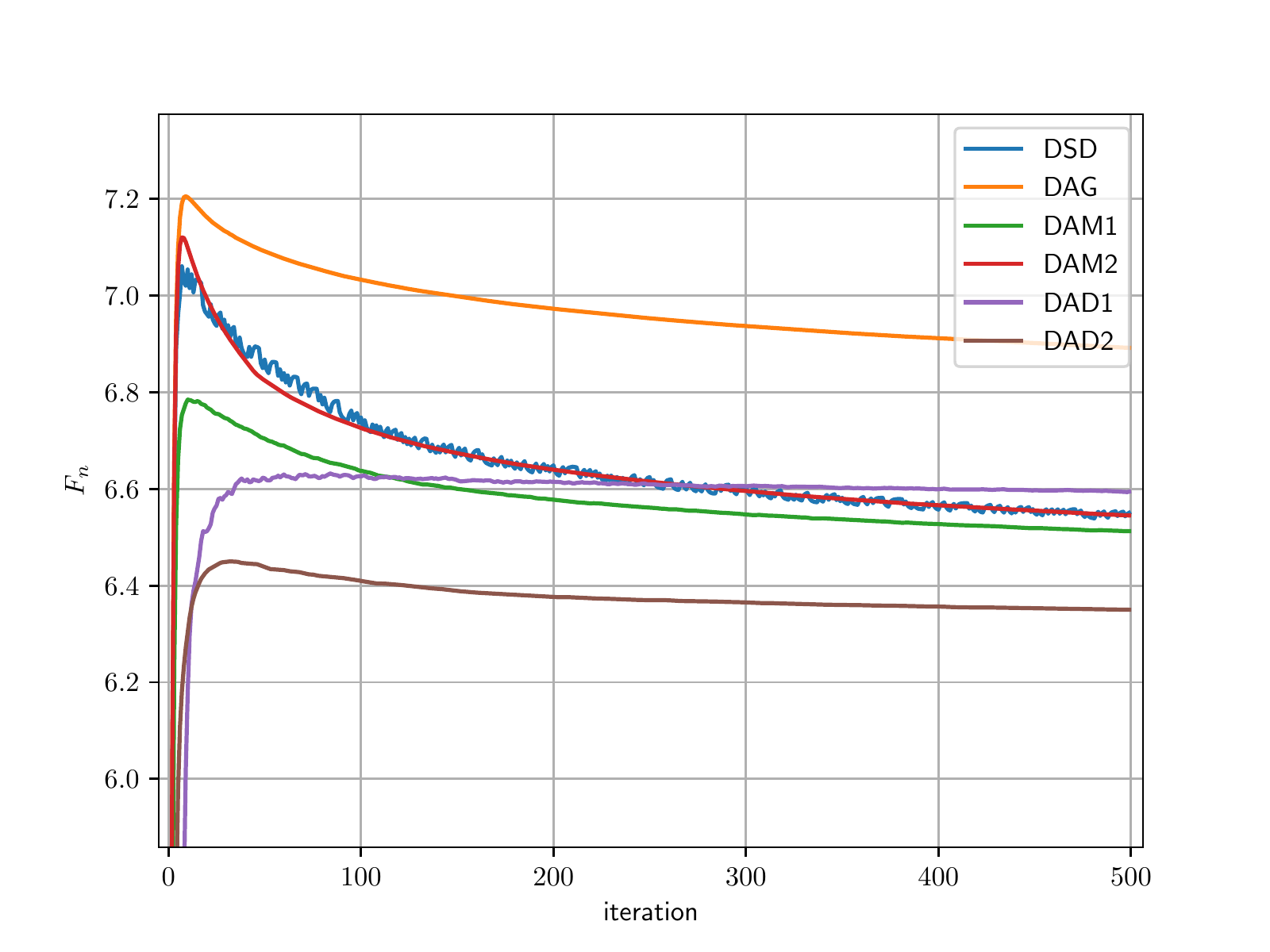}}
\subfigure[$m=10$, $I=5$, $J_i=2$]{\includegraphics[width=39mm]{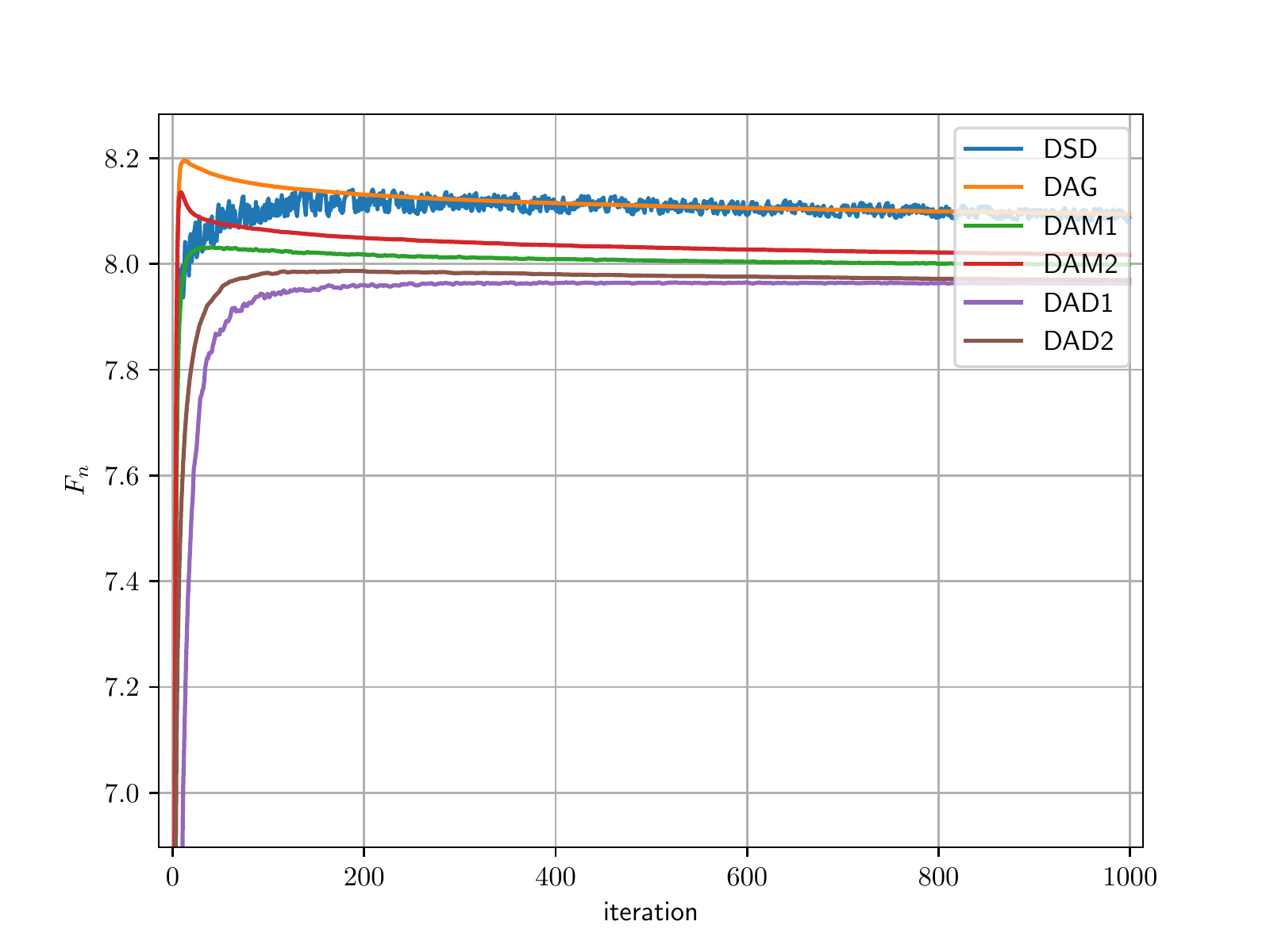}}
\subfigure[$m=100$, $I=5$, $J_i=2$]{\includegraphics[width=39mm]{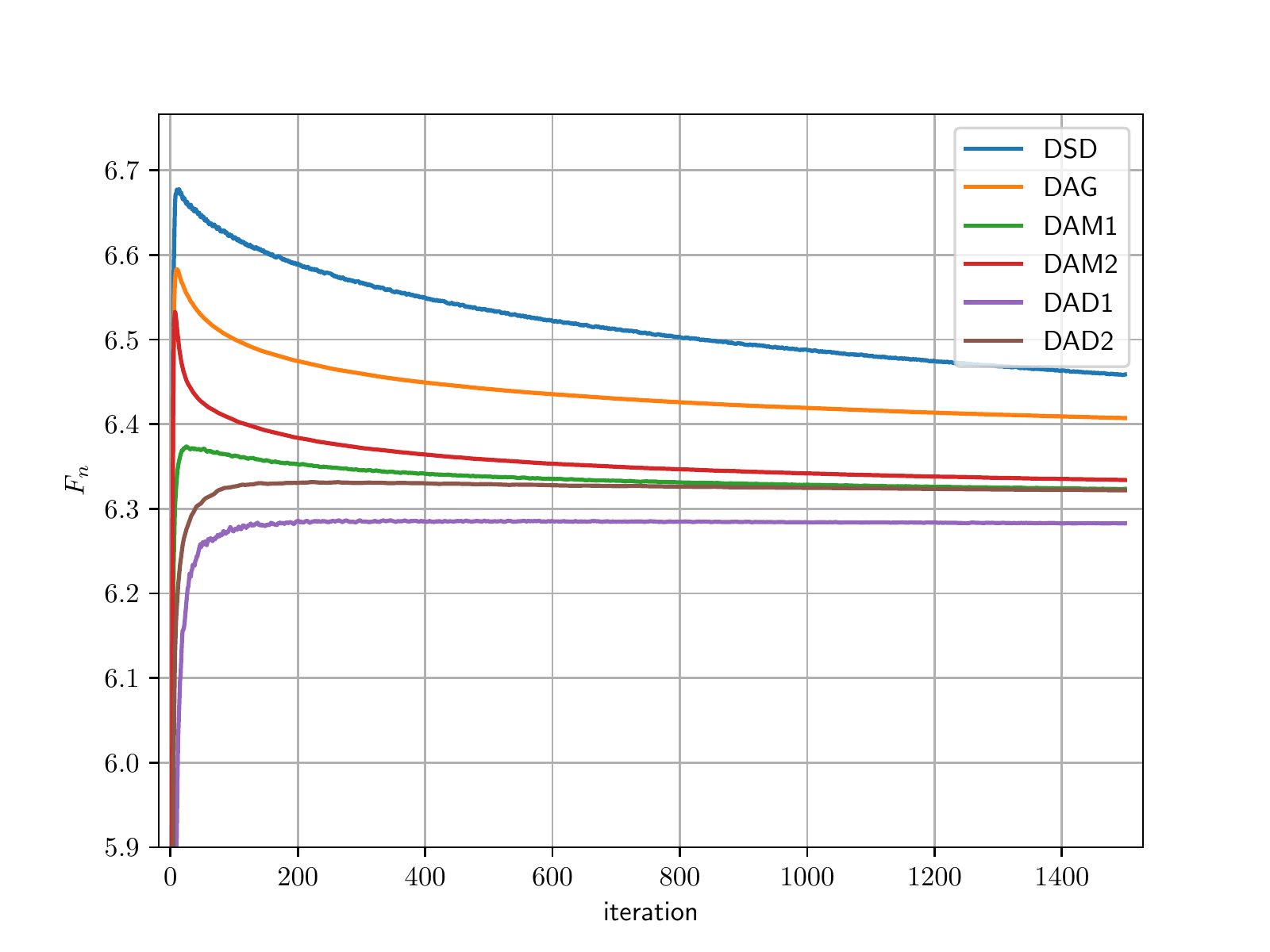}}
\caption{$F_n$ vs. iteration for Algorithm \ref{algo:1} with diminishing step-sizes (inconsistent case)}\label{fig:8}
\end{figure}

\section{Conclusion}\label{sec:7}
This paper proposed the Riemannian stochastic fixed point optimization algorithm
for stochastic optimization with fixed point constraints of quasinonexpansive mappings defined on Riemannian manifolds.
It also gave convergence analyses of the algorithm for both constant 
and diminishing step-sizes and both nonsmooth convex and smooth nonconvex optimization.
For small constant step-sizes, the analyses showed that the algorithm can approximate a solution to the problem. 
For diminishing step-sizes, the analyses suggested the general rate
of convergence of the algorithm.
Finally, the optimality and convergence of the algorithm with each of the formulas 
based on the adaptive learning rate optimization algorithms were demonstrated through numerical comparisons.
In the process, the algorithms with formulas based on Adam and AMSGrad were 
found to be superior for performing stochastic Riemannian optimization with fixed point constraints.

%\section*{Acknowledgments}

%\bibliographystyle{spbasic}     % basic style, author-year citations
\bibliographystyle{spmpsci}      % mathematics and physical sciences
\bibliography{biblio}

\newpage

\section*{Supplementary Material}\label{appen}

\section*{Proofs of Propositions \ref{mappings} and \ref{examples}}

{\em Proof of Proposition \ref{mappings}:}

(i) This follows from the definitions of firmly nonexpansive, nonexpansive,
 and quasinonexpansive mappings. 

(ii)
We prove that \eqref{firm_non} implies \eqref{fquasi}.
The comparison theorem for triangles (see, e.g., \cite[Proposition 2.2]{li2010}),
together with \cite[Proposition 5]{li2011}, ensures that,
for all $x \in C$ and all $y\in \mathrm{Fix}(T)$,
\begin{align*}
&\mathrm{d}(x,T (x))^2 + \mathrm{d}(T (x),y)^2
-2 \left\langle \exp_{T (x)}^{-1} (x), \exp_{T (x)}^{-1}(y) \right\rangle_{T (x)} 
\leq
\mathrm{d}(x,y)^2,\\
&\left\langle \exp_{T (x)}^{-1} (x), \exp_{T (x)}^{-1}(y) \right\rangle_{T (x)} \leq 0,
\end{align*}
which implies \eqref{fquasi}.
From \eqref{quasi}, \eqref{squasi}, and \eqref{fquasi},
we have that 
\eqref{fquasi} implies \eqref{squasi}, and 
\eqref{squasi} implies \eqref{quasi}.

{\em Proof of Proposition \ref{examples}:}

(i)
This follows from \cite[Corollary 1]{li2011}.

(ii)
Proposition \ref{mappings}(i) and Proposition \ref{examples}(i) imply that $P_j$ is nonexpansive. 
Accordingly, $T := P_1 P_2 \cdots P_J$ is nonexpansive. 
Proposition \ref{mappings}(ii) also ensures that 
$P_j$ is strictly quasinonexpansive. 
Hence, the proofs of \cite[Proposition 4.9, Corollary 4.50]{b-c}
lead to Proposition \ref{examples}(ii).

(iii)
This follows from \cite[Theorem 4(i)]{li2011}.

(iv)
The resolvent of $\partial g$ coincides with the Moreau-Yosida regularization of $g$. 
Accordingly, Proposition \ref{examples}(iii) implies Proposition \ref{examples}(iv).

(v)
From the definition of $P_{g,\lambda}$, we have that 
$\mathrm{lev}_{\leq 0}(g) \subset \mathrm{Fix}(P_{g,\lambda})$.
To show that $\mathrm{lev}_{\leq 0}(g) \supset \mathrm{Fix}(P_{g,\lambda})$, 
we assume that $x \in \mathrm{Fix}(P_{g,\lambda})$ and $x \notin \mathrm{lev}_{\leq 0}(g)$.
Then, the definition of $u_x \in \partial g(x)$ 
and the condition $x \notin \mathrm{lev}_{\leq 0}(g)$ guarantee that,
for all $y\in \mathrm{lev}_{\leq 0}(g)$,
\begin{align*}
\left\langle u_x, \exp_x^{-1}(y) \right\rangle_x
\leq g(y) - g(x) \leq - g(x) < 0,
\end{align*}
which implies that $u_x$ is not equal to the zero element $0_x$ 
of $T_x M$.
Accordingly, the definition of $P_{g,\lambda}$ and 
the condition $x \in \mathrm{Fix}(P_{g,\lambda})$ mean that 
\begin{align*}
\exp_x \left( - \lambda \frac{g(x)}{\|u_x\|_x} u_x \right)
= P_{g,\lambda}(x) = x,
\end{align*}
which implies that 
\begin{align*}
0 =
\mathrm{d}
\left(\exp_x \left( - \lambda \frac{g(x)}{\|u_x\|_x} u_x \right),
x
\right)
= 
\lambda \frac{g(x)}{\|u_x\|_x}.
\end{align*}
From $\lambda > 0$ and $u_x \neq 0_x$,
we have that $g(x) = 0$, which is a contradiction from 
$x \notin \mathrm{lev}_{\leq 0}(g)$.
Hence, $\mathrm{lev}_{\leq 0}(g) \supset \mathrm{Fix}(P_{g,\lambda})$.

(vi)
Lemma 5 in \cite{zhang2016} and the definition of $P_{h_j,\lambda}$ ensure that
there exists $\zeta = \zeta (\kappa, D) = \sqrt{|\kappa|}D/\tanh (\sqrt{|\kappa|}D) \in \mathbb{R}_{+}$ such that,
for all $x \in C \backslash \mathrm{lev}_{\leq 0}(h_j)$ and all $y \in \mathrm{lev}_{\leq 0}(h_j) = 
\mathrm{Fix}(P_{h_j,\lambda})$ (by Proposition \ref{examples}(v)),
\begin{align*}
\mathrm{d} (P_{h_j,\lambda}(x),y)^2
&\leq
\zeta \mathrm{d}(P_{h_j,\lambda}(x), x)^2
+ \mathrm{d}(x, y)^2
+ 2 \left\langle \exp_x^{-1}\left(P_{h_j,\lambda}(x)  \right),
\exp_x^{-1}(y) \right\rangle_{x}\\
&=
\zeta \mathrm{d}(P_{h_j,\lambda}(x), x)^2
+ \mathrm{d}(x, y)^2
+ 2 \lambda\frac{h_j(x)}{\|u_{j,x}\|_x^2} 
\left\langle u_{j,x}, \exp_x^{-1}(y) \right\rangle_{x},
\end{align*}
where $(0_x \neq)$ $u_{j,x}\in \partial h_j(x)$,
which, together with the definitions of $P_{h_j,\lambda}$
and $u_{j,x} \in \partial h_j(x)$,
implies that, for $\lambda \in (0,2/\zeta)$,
\begin{align*}
\mathrm{d} (P_{h_j,\lambda}(x),y)^2
&\leq
\mathrm{d}(x, y)^2
+
\zeta \lambda^2 \frac{h_j(x)^2}{\|u_{j,x}\|_x^2} 
- 2 \lambda\frac{h_j(x)^2}{\|u_{j,x}\|_x^2}\\
&=
\mathrm{d}(x, y)^2
+
\lambda ( \zeta \lambda -2 ) \frac{h_j(x)^2}{\|u_{j,x}\|_x^2}.
\end{align*}
From $h_j(x) > 0$, $\mathrm{d} (P_{h_j,\lambda}(x),y) < \mathrm{d} (x,y)$,
i.e.,
$P_{h_j,\lambda}$ is strictly quasinonexpansive. 

(vii)
Proposition \ref{examples}(vi) and the proofs of \cite[Proposition 4.9, Corollary 4.50]{b-c}
lead to Proposition \ref{examples}(vii).

\section*{Proofs of Theorems \ref{cor:1}, \ref{cor:2}, \ref{cor:2_1}, and \ref{cor:2_2}}

The history of the process $\xi_0,\xi_1,\ldots$ up to time $n$ is denoted by 
$\xi_{[n]} = (\xi_0,\xi_1,\ldots,\xi_n)$.
Let $\mathbb{E}[X|\xi_{[n]}]$ denote the conditional expectation of $X$ given 
$\xi_{[n]} = (\xi_0,\xi_1,\ldots,\xi_n)$.
Unless stated otherwise, all relations between random variables are supported to 
hold almost surely.

We prove the following lemma.

\begin{lemma}\label{lem:1}
Suppose that Assumptions (A1)--(A3) and Conditions (C1)--(C2) hold and consider the sequences 
$(x_n)_{n\in\mathbb{N}}$, $(m_n)_{n\in\mathbb{N}}$, and $(\mathsf{d}_n)_{n\in\mathbb{N}}$ defined by Algorithm \ref{algo:1}.
Define $y_n^i$ for all $i\in\mathcal{I}$ and all $n\in \mathbb{N}$ 
by  
\begin{align*}
y_n^i 
:= \exp_{x_n^i}^i \left(\alpha_n \mathsf{d}_n^i \right)
= \exp_{x_n^i}^i \left(- \alpha_n \frac{\hat{m}_n^i}{\mathsf{h}_n^i} \right).
\end{align*}
Then, for all $i\in \mathcal{I}$, there exists a positive number 
$\zeta^i$ such that, for all $x^i\in X^i$ 
and all $n\in\mathbb{N}$, almost surely
\begin{align}\label{INEQ:1}
\begin{split}
\mathrm{d}^i (x_{n+1}^i, x^i )^2
&\leq
\mathrm{d}^i (x_n^i, x^i )^2 
+ \frac{2 \alpha_n}{(1-\hat{\beta}^{n+1})\mathsf{h}_n^i} \left\langle m_n^i, 
\left(\exp_{x_n^i}^i \right)^{-1} (x^i) \right\rangle_{x_n^i}\\
&\quad +
\frac{\zeta^i \alpha_n^2}{(1-\hat{\beta}^{n+1})^2} \frac{\left\| {m}_n^i \right\|_{x_n^i}^2}{(\mathsf{h}_n^i)^2}
- \alpha^i (1-\alpha^i) 
  \mathrm{d}^i \left(T^i (y_n^i), y_n^i \right)^2.
\end{split}   
\end{align}
Moreover, under (C3), for all $i\in \mathcal{I}$, 
there exists a positive number 
$\tilde{{B}^i}^2 := \max \{ \|\tau_{-1}\|_{x_0^i}^2,
{B^i}^2 \}$ such that, for all $n\in\mathbb{N}$,
$\mathbb{E} [ \|m_n^i \|_{x_n^i}^2 ] 
\leq 
\tilde{{B}^i}^2$.  
\end{lemma}

{\em Proof:}
Lemma 5 in \cite{zhang2016}, together with Assumption (A3) (see also \eqref{Di}), guarantees that, for all $i\in \mathcal{I}$,
there exists $\zeta^i = \zeta (\kappa^i, D^i) \in \mathbb{R}_+$ such that,
for all $x^i \in X^i$ and all $n\in \mathbb{N}$,
\begin{align*}
\mathrm{d}^i (y_n^i, x^i )^2
\leq
\zeta^i \mathrm{d}^i \left(y_n^i, x_n^i\right)^2
+ \mathrm{d}^i (x_n^i, x^i)^2
+ 2 \alpha_n \left\langle \frac{\hat{m}_n^i}{\mathsf{h}_n^i}, 
\left(\exp_{x_n^i}^i \right)^{-1} (x^i) \right\rangle_{x_n^i},
\end{align*}
where $\kappa^i$ denotes the lower bound of curvature of $M^i$ and $\zeta(\kappa^i, c) := \sqrt{|\kappa^i|}c/\tanh (\sqrt{|\kappa^i|}c)$
for $c \in \mathbb{R}_+$.
From the definitions of $y_n^i$ and $\hat{m}_n^i$, we have that 
\begin{align}\label{INEQ:2}
\mathrm{d}^i \left(y_n^i, x_n^i \right)^2 
= 
\left\| \left(\exp_{x_n^i}^i \right)^{-1} (y_n^i) \right\|_{x_n^i}^2
= 
\alpha_n^2 \frac{\left\| \hat{m}_n^i \right\|_{x_n^i}^2}{(\mathsf{h}_n^i)^2}
= 
\frac{\alpha_n^2}{(1-\hat{\beta}^{n+1})^2} \frac{\left\| {m}_n^i \right\|_{x_n^i}^2}{(\mathsf{h}_n^i)^2}.
\end{align}
Accordingly, for all $i\in \mathcal{I}$ and all $x^i \in X^i$,
\begin{align*}
&\mathrm{d}^i (y_n^i, x^i)^2\\
&\leq
\mathrm{d}^i (x_n^i, x^i)^2
+
\frac{\zeta^i \alpha_n^2}{(1-\hat{\beta}^{n+1})^2} \frac{\left\| {m}_n^i \right\|_{x_n^i}^2}{(\mathsf{h}_n^i)^2} 
+ \frac{2 \alpha_n}{(1-\hat{\beta}^{n+1})\mathsf{h}_n^i} \left\langle m_n^i, 
\left(\exp_{x_n^i}^i \right)^{-1} (x^i) \right\rangle_{x_n^i}.
\end{align*}
Meanwhile, from $Q_{\alpha^i}^i := P^i S_{\alpha^i}^i$ (see \eqref{2})
and $x_{n+1}^i = Q_{\alpha^i}^i (y_n^i)$,
Proposition \ref{examples}(i) ensures that
\begin{align*}
\mathrm{d}^i (x_{n+1}^i, x^i)^2
\leq 
\mathrm{d}^i ( S_{\alpha^i}^i (y_n^i), x^i)^2,
\end{align*}
which, together with \eqref{keyineq} and \eqref{1}, implies that
\begin{align*}
\mathrm{d}^i (x_{n+1}^i, x^i)^2
\leq 
\mathrm{d}^i ( y_n^i, x^i)^2 - \alpha^i (1-\alpha^i) 
  \mathrm{d}^i \left(T^i (y_n^i), y_n^i \right)^2.
\end{align*} 
Therefore, \eqref{INEQ:1} holds.

The definitions of $m_n$ and $\tau_n$, together with the convexity of 
$\|\cdot\|_{x_n^i}^2$, guarantee that, for all $i\in\mathcal{I}$ and all $n\in\mathbb{N}$,
\begin{align*}
\mathbb{E}\left[ \left\|m_n^i\right\|_{x_n^i}^2   \right]
&\leq \beta_n \mathbb{E}\left[ 
\left\| \varphi_{x_{n-1}^i \to x_n^i}^i (m_{n-1}^i)\right\|_{x_n^i}^2   \right] + (1-\beta_n)\mathbb{E}\left[ \left\|\mathsf{G}^i(x_n, \xi_n) \right\|_{x_n^i}^2   \right]\\
&\leq 
\beta_n \mathbb{E}\left[ \left\|m_{n-1}^i \right\|_{x_{n-1}^i}^2   \right] + (1-\beta_n) {B^i}^2.
\end{align*}
Induction thus ensures that, for all $i\in\mathcal{I}$ and all $n\in\mathbb{N}$,
\begin{align}\label{induction}
\mathbb{E} \left[ \left\|m_n^i \right\|_{x_n^i}^2 \right] 
\leq 
\tilde{{B}^i}^2 := \max \left\{ \left\|\tau_{-1}\right\|_{x_0^i}^2,
{B^i}^2 \right\} < + \infty.
\end{align}
This completes the proof.

Lemma \ref{lem:1} also leads to the following lemma, 
which is used to show the main theorems.

\begin{lemma}\label{lem:2}
Suppose that Assumptions (A1)--(A5) and Conditions (C1)--(C3) hold
and $X_n(x)$ is defined for all $x\in X$ and all $n\in\mathbb{N}$ by
\begin{align*}
X_n(x) := \mathbb{E} \left[\sum_{i\in\mathcal{I}} \mathsf{h}_n^i 
\mathrm{d}^i (x_n^i,x^i) \right].
\end{align*} 
Then, for all $x\in X$ and all $n \in \mathbb{N}$,
\begin{align*}
X_{n+1}(x) 
&\leq
X_n(x)
+ \frac{2\alpha_n (1-\beta_n)}{1 - \hat{\beta}^{n+1}}
\mathbb{E}\left[\left\langle \exp_{x_n}^{-1}(x), \mathsf{g}(x_n)  \right\rangle_{x_n}  \right]\\
&\quad + \mathbb{E}\left[ \sum_{i\in\mathcal{I}} D^i 
\left(\mathsf{h}_{n+1}^i - \mathsf{h}_{n}^i\right)  \right] 
+ \frac{2\alpha_n \beta_n}{1 - \hat{\beta}} 
  \sum_{i\in\mathcal{I}} \tilde{B}^i D^i
  + \frac{\alpha_n^2}{(1-\hat{\beta})^2} \sum_{i\in\mathcal{I}}
  \frac{\zeta^i \tilde{{B}^i}^2}{\mathsf{h}_0^i},
\end{align*}
where $\zeta^i$ and $\tilde{B}^i$ are defined as in Lemma \ref{lem:1}.
\end{lemma}

{\em Proof:}
Condition (C1) and $x_n = x_n(\xi_{[n-1]})$ mean that,
for all $i\in \mathcal{I}$, all $x^i \in X^i$, and all $n\in\mathbb{N}$,
\begin{align*}
&\mathbb{E} \left[ \left\langle \mathsf{G}^i (x_n,\xi_n), 
\left(\exp_{x_n^i}^i \right)^{-1} (x^i) \right\rangle_{x_n^i} \right]\\
&= 
\mathbb{E} \left[ \mathbb{E} \left[ \left\langle \mathsf{G}^i (x_n,\xi_n), 
\left(\exp_{x_n^i}^i \right)^{-1} (x^i) \right\rangle_{x_n^i} \Big| \xi_{[n-1]} \right] \right]\\
&=
\mathbb{E} \left[  \left\langle \mathbb{E} \left[ \mathsf{G}^i (x_n,\xi_n)
\Big| \xi_{[n-1]} \right], 
\left(\exp_{x_n^i}^i \right)^{-1} (x^i) \right\rangle_{x_n^i}  \right]\\
&=
\mathbb{E} \left[  \left\langle  \mathsf{g}^i (x_n),
\left(\exp_{x_n^i}^i \right)^{-1} (x^i) \right\rangle_{x_n^i}  \right].
\end{align*}
Since Condition (C2) implies that, for all $x \in X$ and all $n\in\mathbb{N}$, 
\begin{align*}
%f(x) - f(x_n)
%\geq 
\left\langle  \mathsf{g}(x_n), \exp_{x_n}^{-1} (x) \right\rangle_{x_n}
= 
\sum_{i\in\mathcal{I}}
  \left\langle  \mathsf{g}^i (x_n),
\left(\exp_{x_n^i}^i \right)^{-1} (x^i) \right\rangle_{x_n^i},  
\end{align*}
we have that, for all $x \in X$ and all $n\in\mathbb{N}$, 
\begin{align}\label{f}
\mathbb{E}\left[\left\langle  \mathsf{g}(x_n), \exp_{x_n}^{-1} (x) \right\rangle_{x_n}\right]
=
\mathbb{E} \left[ \sum_{i\in\mathcal{I}} \left\langle \mathsf{G}^i (x_n,\xi_n), 
\left(\exp_{x_n^i}^i \right)^{-1} (x^i) \right\rangle_{x_n^i} \right].
\end{align}
The Cauchy-Schwarz inequality ensures that, 
for all $x \in X$ and all $n\in \mathbb{N}$,
\begin{align*}
\mathbb{E}\left[ \sum_{i\in\mathcal{I}} 
\left\langle \tau_{n-1}^i, (\exp_{x_n^i}^i )^{-1} (x^i) \right\rangle_{x_n^i} \right]
\leq
\mathbb{E}\left[
\sum_{i\in\mathcal{I}} 
\left\| \tau_{n-1}^i \right\|_{x_n^i} 
\left\|\left(\exp_{x_n^i}^i \right)^{-1} (x^i) \right\|_{x_n^i} \right],
\end{align*}
which, together with \eqref{Di} and Lemma \ref{lem:1}, implies that
\begin{align}\label{ii}
\mathbb{E}\left[ \sum_{i\in\mathcal{I}} 
\left\langle \tau_{n-1}^i, (\exp_{x_n^i}^i )^{-1} (x^i) \right\rangle_{x_n^i} \right]
\leq
\sum_{i\in\mathcal{I}} \tilde{B}^i D^i.
\end{align}
Moreover, from Lemma \ref{lem:1}, $\hat{\beta} \in [0,1)$, and 
$1/\mathsf{h}_n^i \leq 1/\mathsf{h}_0^i$ (by Assumption (A4)),
\begin{align}\label{iii}
\begin{split}
\mathbb{E}\left[
\sum_{i\in\mathcal{I}}
\frac{\zeta^i \alpha_n^2}{(1-\hat{\beta}^{n+1})^2} \frac{\left\| {m}_n^i \right\|_{x_n^i}^2}{\mathsf{h}_n^i}
\right]
\leq
\frac{\alpha_n^2}{(1-\hat{\beta})^2} 
\sum_{i\in\mathcal{I}}  \frac{\zeta^i \tilde{{B}^i}^2}
{\mathsf{h}_0^i}.
\end{split}
\end{align}
Accordingly, Lemma \ref{lem:1}, together with \eqref{f}, \eqref{ii}, and \eqref{iii}, leads to the assertion in Lemma \ref{lem:2}.

The following is a convergence analysis of Algorithm \ref{algo:1}.

\begin{theorem}\label{theorem:1}
Suppose that Assumptions (A1)--(A5) and Conditions (C1)--(C3) hold.
Then, Algorithm \ref{algo:1} satisfies that, for all $i\in \mathcal{I}$
and all $n \geq 1$, 
\begin{align}\label{main1}
\begin{split}
&\mathbb{E} \left[ \frac{1}{n} \sum_{k=1}^n
\mathrm{d}^i \left(T^i (y_k^i), y_k^i \right)^2 \right]
\leq 
\frac{1}{\hat{\alpha}^i}
\left\{
\frac{D^i}{n}
+
\frac{2 \tilde{B}^i D^i}{\hat{\mathsf{h}}_0^i} \frac{1}{n} \sum_{k=1}^n \alpha_k
+ 
\frac{\zeta^i \tilde{{B}^i}^2}{(\hat{\mathsf{h}}_0^i)^2}
\frac{1}{n} \sum_{k=1}^n \alpha_k^2
\right\},\\
&\mathbb{E} \left[ \frac{1}{n} \sum_{k=1}^n \mathrm{d}^i \left(y_k^i, x_k^i \right)^2 \right] 
\leq 
\frac{\tilde{{B}^i}^2}{(1-\hat{\beta})^2(\hat{\mathsf{h}}_0^i)^2}
\frac{1}{n} \sum_{k=1}^n \alpha_k^2,
\end{split}
\end{align}
where $\hat{\alpha}^i := \alpha^i (1-\alpha^i)$
and $\hat{\mathsf{h}}_0^i := (1-\hat{\beta})\mathsf{h}_0^i$.
Moreover, if $(\alpha_n (1-\beta_n))_{n\in\mathbb{N}}$ and 
$(\beta_n)_{n\in\mathbb{N}}$ are monotone decreasing, then,
for all $n \geq 1$,
\begin{align}\label{main2}
\begin{split}
&\mathbb{E} \left[ \frac{1}{n} \sum_{k=1}^n 
\left\langle \exp_{x_k}^{-1} (x), \mathsf{g}(x_k) \right\rangle_{x_k} \right]\\
&\quad\leq 
\frac{\sum_{i\in \mathcal{I}} \hat{B}^i {D^i}^2}{2 (1-\beta_1)} \frac{1}{n \alpha_n}
+
\frac{\sum_{i\in\mathcal{I}} \zeta^i \tilde{{B}^i}^2 (\mathsf{h}_0^i)^{-1}}{2(1-\hat{\beta})(1-\beta_1)
} 
\frac{1}{n} \sum_{k=1}^n \alpha_k
 +
\frac{\sum_{i\in\mathcal{I}} \tilde{B}^i D^i}{1 - \beta_1} 
\frac{1}{n} \sum_{k=1}^n \beta_k.
\end{split}
\end{align}
If (A1)'
$T^i \colon M^i \to M^i$ ($i\in \mathcal{I}$) is nonexpansive with $\mathrm{Fix}(T^i) \neq
\emptyset$, then, for all $n \geq 1$, 
\begin{align}\label{triangle_1}
\begin{split}
&\mathbb{E} \left[ \frac{1}{n} \sum_{k=1}^n
\mathrm{d}^i \left(T^i (x_k^i), x_k^i \right)^2 \right]\\
&\quad \leq 
\frac{2}{\hat{\alpha}^i}
\frac{D^i}{n}
+
\frac{4 \tilde{B}^i D^i}{\hat{\alpha}^i \hat{\mathsf{h}}_0^i} \frac{1}{n} \sum_{k=1}^n \alpha_k
+ 
\frac{2 \tilde{{B}^i}^2}{(\hat{\mathsf{h}}_0^i)^2}
\left\{
\frac{\zeta^i}{\hat{\alpha}^i}
+ \frac{4}{(1-\hat{\beta})^2}
\right\}
\frac{1}{n} \sum_{k=1}^n \alpha_k^2.
\end{split}
\end{align}
\end{theorem}

{\em Proof:}
The Cauchy-Schwarz inequality, together with Lemma \ref{lem:1} and Assumption (A3) (see \eqref{Di}), ensures that,
for all $i\in \mathcal{I}$, all $x^i \in X^i$, and 
all $n\in\mathbb{N}$, 
\begin{align*}
\mathbb{E} \left[ \left\langle m_n^i, 
(\exp_{x_n^i}^i )^{-1} (x^i) \right\rangle_{x_n^i} \right]
\leq
\mathbb{E} \left[ \left\|m_n^i \right\|_{x_n^i} \left\| (\exp_{x_n^i}^i )^{-1} (x^i)\right\|_{x_n^i} \right]
\leq 
\tilde{B}^i D^i.
\end{align*} 
Lemma \ref{lem:1}, together with $\hat{\beta}\in [0,1)$ and Assumption (A4), 
guarantees that, 
for all $i\in \mathcal{I}$, all $x^i \in X^i$, and 
all $k\geq 1$,
\begin{align}\label{Tx}
\begin{split}
&\alpha^i (1-\alpha^i) \mathbb{E} \left[
\mathrm{d}^i \left( T^i (y_k^i), y_k^i \right)^2 \right]\\
&\quad\leq
\mathbb{E} \left[ \mathrm{d}^i (x_k^i, x^i)^2 \right]
-
\mathbb{E} \left[ \mathrm{d}^i (x_{k+1}^i, x^i)^2 \right]
+
\frac{\zeta^i  \tilde{{B}^i}^2}{(1-\hat{\beta})^2 (\mathsf{h}_0^i)^2} \alpha_k^2
+ \frac{2 \tilde{B}^i D^i}{(1-\hat{\beta})\mathsf{h}_0^i} \alpha_k.
\end{split}
\end{align}
Accordingly, we have that, for all $i\in \mathcal{I}$
and all $n\geq 1$,
\begin{align*}
\hat{\alpha}^i \mathbb{E} \left[ \sum_{k=1}^n
\mathrm{d}^i \left(T^i (y_k^i), y_k^i \right)^2 \right]
&\leq
D^i
+
\frac{\zeta^i  \tilde{{B}^i}^2}{(1-\hat{\beta})^2 (\mathsf{h}_0^i)^2} \sum_{k=1}^n \alpha_k^2
 + \frac{2 \tilde{B}^i D^i}{(1-\hat{\beta})\mathsf{h}_0^i} \sum_{k=1}^n \alpha_k,
\end{align*}
where 
$\hat{\alpha}^i := \alpha^i (1-\alpha^i)$
and
\eqref{Di} implies that $\mathbb{E} \left[ \mathrm{d}^i (x_1^i, x^i)^2 \right] \leq D^i$.
From \eqref{INEQ:2}, for all $i\in \mathcal{I}$
and all $n\geq 1$,
\begin{align*}
\mathbb{E} \left[ \sum_{k=1}^n \mathrm{d}^i \left(y_k^i, x_k^i \right)^2 \right] 
\leq 
\frac{\tilde{{B}^i}^2}{(1-\hat{\beta})^2 (\mathsf{h}_0^i)^2}
\sum_{k=1}^n \alpha_k^2,
\end{align*}
which implies that \eqref{main1} holds.
Lemma \ref{lem:1}, together with the definition of $m_n^i$, implies that, for all $i\in\mathcal{I}$,
all $x^i \in X^i$, and all $n \in \mathbb{N}$, 
\begin{align*}
&\left\langle - \mathsf{G}^i (x_n,\xi_n), 
\left(\exp_{x_n^i}^i \right)^{-1} (x^i) \right\rangle_{x_n^i}\\
&\leq
\underbrace{\frac{(1 - \hat{\beta}^{n+1}) \mathsf{h}_n^i}{2 \alpha_n (1-\beta_n)} \left\{ \mathrm{d}^i (x_n^i, x^i)^2 
-
\mathrm{d}^i (x_{n+1}^i, x^i)^2 \right\}}_{H_n^i(x^i)}\\
&\quad + 
\underbrace{\frac{\beta_n}{1-\beta_n} 
\left\langle \tau_{n-1}^i, 
\left(\exp_{x_n^i}^i \right)^{-1} (x^i) \right\rangle_{x_n^i}}_{B_n^i(x^i)}
 +
\underbrace{\frac{\zeta^i \alpha_n}{2(1-\hat{\beta}^{n+1})(1-\beta_n)} \frac{\left\| {m}_n^i \right\|_{x_n^i}^2}{\mathsf{h}_n^i}}_{A_n^i(x^i)},
\end{align*}
which, together with \eqref{f}, implies that, for all $x_\star \in X_\star$ and all $n \geq 1$,
\begin{align}\label{iv}
\begin{split}
&\mathbb{E} \left[ \frac{1}{n} \sum_{k=1}^n 
\left\langle \exp_{x_k}^{-1} (x), \mathsf{g}(x_k) \right\rangle_{x_k} \right]\\
&\leq
\frac{1}{n} \mathbb{E}\left[ \sum_{k=1}^n \sum_{i\in\mathcal{I}} H_k^i (x_\star^i) \right]
+
\frac{1}{n} \mathbb{E}\left[ \sum_{k=1}^n \sum_{i\in\mathcal{I}} B_k^i 
(x_\star^i) \right]
+
\frac{1}{n} \mathbb{E}\left[ \sum_{k=1}^n \sum_{i\in\mathcal{I}} A_k^i 
(x_\star^i)\right].
\end{split}
\end{align}
The definition of $H_n^i (x^i)$ ($i\in\mathcal{I},n\in\mathbb{N}$)
and \eqref{Di} guarantee that, for all $i\in\mathcal{I}$, all $x_\star \in X_\star$, and all $n\geq 1$,
\begin{align*}
&\sum_{k=1}^n  H_k^i (x_\star^i)\\
&\leq
\frac{(1 - \hat{\beta}^{2}) \mathsf{h}_1^i}{2 \alpha_1 (1-\beta_1)} {D^i}^2 
+ 
\sum_{k=2}^n \left\{
\frac{(1 - \hat{\beta}^{k+1}) \mathsf{h}_k^i}{2 \alpha_k (1-\beta_k)}   
-
\frac{(1 - \hat{\beta}^{k}) \mathsf{h}_{k-1}^i}{2 \alpha_{k-1} (1-\beta_{k-1})}   
\right\}
\mathrm{d}^i (x_k^i, x_\star^i)^2.
\end{align*}
Since $\hat{\beta} \in [0,1)$ and Assumption (A4) hold 
and $(\alpha_n(1-\beta_n))_{n\in\mathbb{N}}$ is monotone decreasing, 
we have that, for all $k \geq 2$,
\begin{align*}
\frac{(1 - \hat{\beta}^{k+1}) \mathsf{h}_k^i}{2 \alpha_k (1-\beta_k)}   
-
\frac{(1 - \hat{\beta}^{k}) \mathsf{h}_{k-1}^i}{2 \alpha_{k-1} (1-\beta_{k-1})} \geq 0.
\end{align*}
Accordingly, for all $i\in\mathcal{I}$ and all $x_\star \in X_\star$,
\begin{align}\label{i}
\begin{split}
&\mathbb{E}\left[\sum_{k=1}^n  H_k^i (x_\star^i) \right]\\
&\leq
\mathbb{E}\left[\frac{(1 - \hat{\beta}^{2}) \mathsf{h}_1^i}{2 \alpha_1 (1-\beta_1)} {D^i}^2 
+ 
\sum_{k=2}^n \left\{
\frac{(1 - \hat{\beta}^{k+1}) \mathsf{h}_k^i}{2 \alpha_k (1-\beta_k)}   
-
\frac{(1 - \hat{\beta}^{k}) \mathsf{h}_{k-1}^i}{2 \alpha_{k-1} (1-\beta_{k-1})}   
\right\}
{D^i}^2 \right]\\
&=
\mathbb{E}\left[\frac{(1 - \hat{\beta}^{n+1}) \mathsf{h}_n^i}{2 \alpha_n (1-\beta_n)}
{D^i}^2 \right]\\
&\leq
\frac{\hat{B}^i {D^i}^2}{2 (1-\beta_1) \alpha_n},
\end{split}
\end{align}
where the second inequality comes from $\hat{\beta} \in [0,1)$, 
Assumption (A5),
and $\beta_n \leq \beta_1$ ($n\geq 1$).
The Cauchy-Schwarz inequality ensures that, 
for all $x_\star \in X_\star$ and all $n\geq 1$,
\begin{align*}
\mathbb{E}\left[ \sum_{k=1}^n \sum_{i\in\mathcal{I}} B_k^i (x_\star^i) \right]
\leq
\mathbb{E}\left[
\sum_{i\in\mathcal{I}} \sum_{k=1}^n
\frac{\beta_k}{1-\beta_k} 
\left\| \tau_{k-1}^i \right\|_{x_k^i} 
\left\|(\exp_{x_k^i}^i )^{-1} (x_\star^i) \right\|_{x_k^i} \right],
\end{align*}
which, together with \eqref{Di}, Lemma \ref{lem:1}, and 
$\beta_n \leq \beta_1$ ($n\geq 1$), implies that
\begin{align}\label{ii_1}
\mathbb{E}\left[ \sum_{k=1}^n \sum_{i\in\mathcal{I}} B_k^i (x_\star^i) \right]
\leq
\frac{\sum_{i\in\mathcal{I}} \tilde{B}^i D^i}{1 - \beta_1} \sum_{k=1}^n \beta_k.
\end{align}
Moreover, from Lemma \ref{lem:1}, $\hat{\beta} \in [0,1)$, 
Assumption (A4), and $\beta_n \leq \beta_1$ ($n\geq 1$),
\begin{align}\label{iii_1}
\begin{split}
\mathbb{E}\left[ \sum_{k=1}^n \sum_{i\in\mathcal{I}} A_k^i 
(x_\star^i)\right]
&=
\mathbb{E}\left[
\sum_{k=1}^n \sum_{i\in\mathcal{I}}
\frac{\zeta^i \alpha_k}{2(1-\hat{\beta}^{k+1})(1-\beta_k)} \frac{\left\| {m}_k^i \right\|_{x_k^i}^2}{\mathsf{h}_k^i}
\right]\\
&\leq 
\frac{1}{2(1-\hat{\beta})(1-\beta_1)}
\sum_{i\in\mathcal{I}} \frac{\zeta^i \tilde{{B}^i}^2}{\mathsf{h}_0^i} 
\sum_{k=1}^n \alpha_k.
\end{split}
\end{align}
Hence, \eqref{iv}, \eqref{i}, \eqref{ii_1}, and \eqref{iii_1}
lead to \eqref{main2}. 

Suppose that Assumption (A1)' holds.
Since the triangle inequality implies that, for all $i\in\mathcal{I}$ and all $n\in\mathbb{N}$,
\begin{align*}
\mathrm{d}^i \left(T^i(x_n^i),x_n^i \right) 
\leq 
\mathrm{d}^i \left(T^i(x_n^i),T^i(y_n^i) \right)
+
\mathrm{d}^i \left(T^i(y_n^i),y_n^i \right)
+
\mathrm{d}^i \left(y_n^i,x_n^i \right),
\end{align*}
Assumption (A1)' ensures that 
\begin{align*}
\mathrm{d}^i \left(T^i(x_n^i),x_n^i \right) 
\leq 
\mathrm{d}^i \left(T^i(y_n^i),y_n^i \right)
+
2 \mathrm{d}^i \left(y_n^i,x_n^i \right),
\end{align*}
which implies that, for all $i\in\mathcal{I}$ and all $n\in\mathbb{N}$,
\begin{align}\label{triangle}
\mathrm{d}^i \left(T^i(x_n^i),x_n^i \right)^2 
\leq 
2 \mathrm{d}^i \left(T^i(y_n^i),y_n^i \right)^2
+
8 \mathrm{d}^i \left(y_n^i,x_n^i \right)^2.
\end{align}
\eqref{main1} and \eqref{triangle} thus lead to \eqref{triangle_1},
which completes the proof.

{\em Proof of Theorem \ref{cor:1}:}
Let $i\in\mathcal{I}$ be fixed arbitrarily.
From \eqref{INEQ:2} and Lemma \ref{lem:1}, together with 
Assumption (A4) and 
$\alpha_n := \alpha$ ($n\in\mathbb{N}$), 
\eqref{limsup} holds.
If \eqref{liminf} does not hold, then there exists $\delta > 0$ 
such that
\begin{align*}
\alpha^i (1-\alpha^i) \liminf_{n \to + \infty}
\mathbb{E} \left[
\mathrm{d}^i \left( T^i (y_n^i), y_n^i \right)^2 \right]
>
\frac{2 \tilde{B}^i D^i}{(1-\hat{\beta})\mathsf{h}_0^i} \alpha
+
\frac{\zeta^i  \tilde{{B}^i}^2}{(1-\hat{\beta})^2 (\mathsf{h}_0^i)^2} \alpha^2
+ \delta.
\end{align*}
The definition of the limit inferior of 
$(\mathbb{E}[\mathrm{d}^i ( T^i (y_n^i), y_n^i)^2 ])_{n\in\mathbb{N}}$
ensures that there exists $n_0\in \mathbb{N}$ such that,
for all $n \geq n_0$,
\begin{align*}
\alpha^i (1-\alpha^i)\liminf_{n \to + \infty}
\mathbb{E} \left[
\mathrm{d}^i \left( T^i (y_n^i), y_n^i \right)^2 \right] 
 - \frac{1}{2} \delta 
\leq 
\alpha^i (1-\alpha^i) \mathbb{E} \left[
\mathrm{d}^i \left( T^i (y_n^i), y_n^i \right)^2 \right],
\end{align*}
which implies that, for all $n \geq n_0$,
\begin{align*}
\alpha^i (1-\alpha^i) \mathbb{E} \left[
\mathrm{d}^i \left( T^i (y_n^i), y_n^i \right)^2 \right]
>
\frac{2 \tilde{B}^i D^i}{(1-\hat{\beta})\mathsf{h}_0^i} \alpha
+
\frac{\zeta^i  \tilde{{B}^i}^2}{(1-\hat{\beta})^2 (\mathsf{h}_0^i)^2} \alpha^2
+ \frac{1}{2} \delta.
\end{align*}
From \eqref{Tx} with $\alpha_n := \alpha$ and $\beta_n := \beta$ ($n\in\mathbb{N}$),
\begin{align*}
\mathbb{E} \left[ \mathrm{d}^i (x_{n+1}^i, x^i)^2 \right]
&\leq
\mathbb{E} \left[ \mathrm{d}^i (x_n^i, x^i)^2 \right]
-
\alpha^i (1-\alpha^i) \mathbb{E} \left[
\mathrm{d}^i \left( T^i (y_n^i), y_n^i \right)^2 \right]\\
&\quad 
+ \frac{2 \tilde{B}^i D^i}{(1-\hat{\beta})\mathsf{h}_0^i} \alpha
+
\frac{\zeta^i  \tilde{{B}^i}^2}{(1-\hat{\beta})^2 (\mathsf{h}_0^i)^2} \alpha^2,
\end{align*}
which implies that, for all $n \geq n_0$,
\begin{align*}
\mathbb{E} \left[ \mathrm{d}^i (x_{n+1}^i, x^i)^2 \right]
&<
\mathbb{E} \left[ \mathrm{d}^i (x_n^i, x^i)^2 \right]
-
\left\{ \frac{2 \tilde{B}^i D^i}{(1-\hat{\beta})\mathsf{h}_0^i} \alpha
+
\frac{\zeta^i  \tilde{{B}^i}^2}{(1-\hat{\beta})^2 (\mathsf{h}_0^i)^2} \alpha^2
+ \frac{1}{2} \delta \right\} \\
&\quad 
+ \frac{2 \tilde{B}^i D^i}{(1-\hat{\beta})\mathsf{h}_0^i} \alpha
+
\frac{\zeta^i  \tilde{{B}^i}^2}{(1-\hat{\beta})^2 (\mathsf{h}_0^i)^2} \alpha^2\\
&= 
\mathbb{E} \left[ \mathrm{d}^i (x_n^i, x^i)^2 \right] 
- \frac{1}{2} \delta\\ 
&< 
\mathbb{E} \left[ \mathrm{d}^i (x_{n_0}^i, x^i)^2 \right] 
- \frac{1}{2} \delta (n+1 -n_0).
\end{align*}
Since the right-hand side of the above inequality approaches
minus infinity when $n$ diverges, we have a contradiction.
Hence, \eqref{liminf} holds.

Assumptions (A4) and (A5) and 
the conditions, $\lim_{n\to + \infty} \hat{\beta}^{n+1} = 0$
and $X_n^\star := X_n(x_\star) \leq \sum_{i\in\mathcal{I}} \hat{B}^i D^i < + \infty$
($x_\star \in X_\star$) (by Assumptions (A3) and (A5)),
guarantee that, for all $\epsilon > 0$,
there exists $n_1 \in \mathbb{N}$ such that, for all $n \in \mathbb{N}$,
$n \geq n_1$ implies that 
\begin{align}\label{epsilon}
\mathbb{E}\left[ \sum_{i\in\mathcal{I}} D^i 
\left(\mathsf{h}_{n+1}^i - \mathsf{h}_{n}^i\right)  \right]
+ 
\hat{\beta}^{n+1} \left(X_{n+1}^\star - X_n^\star \right)
\leq \alpha (1-\beta)\epsilon.
\end{align}
Let us show that, for all $\epsilon > 0$,
\begin{align}\label{liminf_1}
\liminf_{n\to + \infty} \mathbb{E}\left[ f(x_n) - f_\star \right]
\leq 
\frac{\sum_{i\in\mathcal{I}} \zeta^i \tilde{{B}^i}^2 (\mathsf{h}_0^i)^{-1}}{2(1-\beta)(1-\hat{\beta})} \alpha
+
\frac{\sum_{i\in\mathcal{I}} \tilde{B}^i D^i}{(1-\beta)(1-\hat{\beta})} \beta
+
\frac{3}{2} \epsilon.
\end{align}
If \eqref{liminf_1} does not hold, then there exists $\epsilon_0 > 0$
such that
\begin{align*}
\liminf_{n\to + \infty} \mathbb{E}\left[ f(x_n) - f_\star \right]
> 
\frac{\sum_{i\in\mathcal{I}} \zeta^i \tilde{{B}^i}^2 (\mathsf{h}_0^i)^{-1}}{2(1-\beta)(1-\hat{\beta})} \alpha
+
\frac{\sum_{i\in\mathcal{I}} \tilde{B}^i D^i}{(1-\beta)(1-\hat{\beta})} \beta
+
\frac{3}{2} \epsilon_0.
\end{align*}
From the definition of the limit inferior of 
$(\mathbb{E}[ f(x_n) - f_\star ])_{n\in\mathbb{N}}$,
there exists $n_2 \in \mathbb{N}$ such that, for all $n \geq n_2$,
\begin{align*}
\liminf_{n\to + \infty} \mathbb{E}\left[ f(x_n) - f_\star \right]
- \frac{1}{2} \epsilon_0
\leq
\mathbb{E}\left[ f(x_n) - f_\star \right].
\end{align*}
Hence, we have that, for all $n \geq n_2$,
\begin{align*}
\mathbb{E}\left[ f(x_n) - f_\star \right]
> 
\frac{\sum_{i\in\mathcal{I}} \zeta^i \tilde{{B}^i}^2 (\mathsf{h}_0^i)^{-1}}{2(1-\beta)(1-\hat{\beta})} \alpha
+
\frac{\sum_{i\in\mathcal{I}} \tilde{B}^i D^i}{(1-\beta)(1-\hat{\beta})} \beta
+ \epsilon_0.
\end{align*}
The convexity of $f$ implies that, for all $n\in\mathbb{N}$,
\begin{align}\label{convex}
\mathbb{E}\left[\left\langle \exp_{x_n}^{-1} (x_\star), \mathsf{g}(x_n) \right\rangle_{x_n} \right]
\leq f_\star - f(x_n).
\end{align}
Since Lemma \ref{lem:2}, together with $\alpha_n := \alpha$, $\beta_n := \beta$ ($n\in\mathbb{N}$), \eqref{epsilon}, and \eqref{convex},
ensures that, for all $n \geq n_1$,
\begin{align*}
X_{n+1}^\star
&\leq
X_{n}^\star
+ \alpha (1-\beta)\epsilon
- 2\alpha (1-\beta)
\mathbb{E}\left[f(x_n) - f_\star\right]
 + \frac{2\alpha \beta}{1 - \hat{\beta}} 
  \sum_{i\in\mathcal{I}} \tilde{B}^i D^i\\
&\quad + \frac{\alpha^2}{1-\hat{\beta}} \sum_{i\in\mathcal{I}}
  \frac{\zeta^i \tilde{{B}^i}^2}{\mathsf{h}_0^i},
\end{align*}
we find that, for all $n \geq n_3 := \max \{ n_1,n_2 \}$,
\begin{align*}
X_{n+1}^\star 
&<
X_{n}^\star
+ \alpha (1-\beta)\epsilon_0
+ \frac{2\alpha \beta}{1 - \hat{\beta}} 
  \sum_{i\in\mathcal{I}} \tilde{B}^i D^i
  + \frac{\alpha^2}{1-\hat{\beta}} \sum_{i\in\mathcal{I}}
  \frac{\zeta^i \tilde{{B}^i}^2}{\mathsf{h}_0^i}\\
&\quad - 2\alpha (1-\beta)
\left\{
\frac{\sum_{i\in\mathcal{I}} \zeta^i \tilde{{B}^i}^2 (\mathsf{h}_0^i)^{-1}}{2(1-\beta)(1-\hat{\beta})} \alpha
+
\frac{\sum_{i\in\mathcal{I}} \tilde{B}^i D^i}{(1-\beta)(1-\hat{\beta})} \beta
+ \epsilon_0
\right\}\\
&= 
X_{n}^\star
- \alpha (1-\beta)\epsilon_0\\
&< 
X_{n_3}^\star
- \alpha (1-\beta)\epsilon_0 \left(n+1 - n_3 \right),
\end{align*}
which is a contradiction. 
Hence, \eqref{liminf_1} holds for all $\epsilon > 0$.
This implies that \eqref{liminf_2} holds.
Obviously, \eqref{ave_1} holds from \eqref{main1} with 
$\alpha_n := \alpha$ and $\beta_n := \beta$ ($n\in\mathbb{N}$).

The conditions $\alpha_n := \alpha$ and $\beta_n := \beta$ ($n\in\mathbb{N}$)
satisfy that 
$(\alpha_n(1-\beta_n))_{n\in\mathbb{N}}$ and $(\beta_n)_{n\in\mathbb{N}}$ are monotone decreasing.
Since $f$ is convex, induction shows that
\begin{align}\label{conv}
f(\bar{x}_n) \leq \frac{1}{n} \sum_{k=1}^n f(x_k).
\end{align}
Therefore, \eqref{main2} in Theorem \ref{theorem:1} leads to \eqref{ave_2}.
If Assumption (A1)' holds, then Theorem \ref{theorem:1} leads to \eqref{av_0}.

{\em Proof of Theorem \ref{cor:2}:}
From \eqref{INEQ:2} and Lemma \ref{lem:1}, together with Assumption (A4) and 
$\lim_{n\to + \infty} \alpha_n = 0$ (by $\sum_{n=0}^{+\infty} \alpha_n^2 < + \infty$), we have that $\lim_{n \to +\infty} \mathbb{E}[ \mathrm{d}^i (y_n^i, x_n^i )^2 ] = 0$.
Define $Y_n^i$ for all $x\in X$, all $i\in\mathcal{I}$, 
and all $n\in\mathbb{N}$ by 
\begin{align*}
Y_n^i (x) := \alpha_n \mathbb{E}\left[ \mathrm{d}^i (x_n^i,x^i)^2\right].
\end{align*}
Inequality \eqref{Tx} then ensures that, 
for all $x\in X$, all $i\in\mathcal{I}$, 
and all $k\in\mathbb{N}$,
\begin{align*}
&\alpha^i (1-\alpha^i) \alpha_k \mathbb{E} \left[
\mathrm{d}^i \left( T^i (y_k^i), y_k^i \right)^2 \right]\\
&\quad\leq
Y_k^i (x)
-
\alpha_k \mathbb{E} \left[ \mathrm{d}^i (x_{k+1}^i, x^i)^2 \right]
+
\frac{\zeta^i  \tilde{{B}^i}^2}{(1-\hat{\beta})^2 (\mathsf{h}_0^i)^2} \alpha_k^3
+ \frac{2 \tilde{B}^i D^i}{(1-\hat{\beta})\mathsf{h}_0^i} \alpha_k^2,
\end{align*}
which, together with $\alpha_{n+1} \leq \alpha_n$ ($n\in\mathbb{N}$),
implies that
\begin{align*}
&\alpha^i (1-\alpha^i) \alpha_k \mathbb{E} \left[
\mathrm{d}^i \left( T^i (y_k^i), y_k^i \right)^2 \right]\\
&\quad\leq
Y_k^i (x)
-
Y_{k+1}^i (x)
+
\frac{\zeta^i  \tilde{{B}^i}^2}{(1-\hat{\beta})^2 (\mathsf{h}_0^i)^2} \alpha_k^3
+ \frac{2 \tilde{B}^i D^i}{(1-\hat{\beta})\mathsf{h}_0^i} \alpha_k^2.
\end{align*}
Summing the above inequality from $k=0$ to $k=n$ means that,
for all $x\in X$, all $i\in\mathcal{I}$, 
and all $n\in\mathbb{N}$,
\begin{align*}
\hat{\alpha}^i \sum_{k=0}^n \alpha_k \mathbb{E} \left[
\mathrm{d}^i \left( T^i (y_k^i), y_k^i \right)^2 \right]
\leq 
Y_0^i (x)
+
\frac{\zeta^i  \tilde{{B}^i}^2}{(1-\hat{\beta})^2 (\mathsf{h}_0^i)^2} 
\sum_{k=0}^n \alpha_k^3
+ \frac{2 \tilde{B}^i D^i}{(1-\hat{\beta})\mathsf{h}_0^i} \sum_{k=0}^n \alpha_k^2,
\end{align*}
where $\hat{\alpha}^i := \alpha^i (1-\alpha^i)$.
Since $(\alpha_n)_{n\in\mathbb{N}} \subset (0,1)$ satisfies 
$\sum_{n=0}^{+\infty} \alpha_n^2 < + \infty$,
we have that, for all $x\in X$ and all $i\in\mathcal{I}$, 
\begin{align}\label{SUM}
\sum_{n=0}^{+\infty} \alpha_n \mathbb{E} \left[
\mathrm{d}^i \left( T^i (y_n^i), y_n^i \right)^2 \right] 
< + \infty.
\end{align}
We prove that, for all $i\in \mathcal{I}$, 
\begin{align}\label{main:3}
\liminf_{n\to +\infty} \mathbb{E} \left[
\mathrm{d}^i \left( T^i (y_n^i), y_n^i \right)^2 \right] \leq 0.
\end{align}
Assume that \eqref{main:3} does not hold. 
Then, there exist $i\in\mathcal{I}$, $\gamma > 0$, and 
$m_0 \in \mathbb{N}$ such that, for all $n\geq m_0$,
\begin{align*}
\mathbb{E} \left[
\mathrm{d}^i \left( T^i (y_n^i), y_n^i \right)^2 \right] \geq \gamma,
\end{align*}
which, together with $\sum_{n=0}^{+\infty} \alpha_n = + \infty$ and \eqref{SUM}, implies that
\begin{align*}
+\infty
=
\gamma \sum_{n=m_0}^{+\infty} \alpha_n
\leq
\sum_{n=m_0}^{+\infty} \alpha_n \mathbb{E} \left[
\mathrm{d}^i \left( T^i (y_n^i), y_n^i \right)^2 \right] 
< + \infty.
\end{align*}
This is a contradiction. Hence, we have \eqref{main:3},
which implies that, for all $i\in\mathcal{I}$, 
$\liminf_{n\to +\infty} \mathbb{E} [
\mathrm{d}^i ( T^i (y_n^i), y_n^i )^2 ] = 0$.
Lemma \ref{lem:2} with \eqref{convex} guarantees that, for all $k\in\mathbb{N}$,
\begin{align*}
\frac{2\alpha_k}{1 - \hat{\beta}^{k+1}}
\mathbb{E}\left[f(x_k) - f_\star \right]
&\leq
X_k^\star
-
X_{k+1}^\star
+ \mathbb{E}\left[ \sum_{i\in\mathcal{I}} D^i 
\left(\mathsf{h}_{k+1}^i - \mathsf{h}_{k}^i\right)  \right]\\
&\quad + \frac{2\alpha_k \beta_k}{1 - \hat{\beta}} 
  \left( \sum_{i\in\mathcal{I}} \tilde{B}^i D^i + F \right)
  + \frac{\alpha_k^2}{(1-\hat{\beta})^2} \sum_{i\in\mathcal{I}}
  \frac{\zeta^i \tilde{{B}^i}^2}{\mathsf{h}_0^i},
\end{align*}
where $X_n^\star := X_n(x^\star)$ and $F := \sup \{|\mathbb{E}[f(x_n) - f_\star ]| \colon n\in \mathbb{N}\}$ 
is finite from Assumptions (A2) and (A3).
Summing the above inequality from $k=0$ to $k=n$ implies that
\begin{align*}
2 \sum_{k=0}^n \frac{\alpha_k}{1 - \hat{\beta}^{k+1}}
\mathbb{E}\left[f(x_k) - f_\star \right]
&\leq
X_0^\star
+ \sum_{i\in\mathcal{I}} D^i B^i
+ \frac{1}{(1-\hat{\beta})^2} \sum_{i\in\mathcal{I}}
  \frac{\zeta^i \tilde{{B}^i}^2}{\mathsf{h}_0^i}
  \sum_{k=0}^n \alpha_k^2\\
&\quad + \frac{2}{1 - \hat{\beta}} 
  \left( \sum_{i\in\mathcal{I}} \tilde{B}^i D^i + F \right)
  \sum_{k=0}^n \alpha_k \beta_k,
\end{align*}
where $\sup \{\mathbb{E}[\mathsf{h}_n^i] \colon n\in \mathbb{N} \} \leq \hat{B}^i$ holds 
from Assumption (A5).
From $\sum_{n=0}^{+\infty} \alpha_n \beta_n < + \infty$ and 
$\sum_{n=0}^{+\infty} \alpha_n^2 < + \infty$, 
we have that
\begin{align*}
\sum_{k=0}^{+\infty} \frac{\alpha_k}{1 - \hat{\beta}^{k+1}}
\mathbb{E}\left[f(x_k) - f_\star \right] 
< + \infty.
\end{align*}
We also have that 
$\sum_{n=0}^{+\infty} \alpha_k/(1 - \hat{\beta}^{n+1})
\geq 
\sum_{n=0}^{+\infty} \alpha_k = + \infty$.
Accordingly, a discussion similar to the one for proving \eqref{main:3} 
leads to the finding that
\begin{align*}
\liminf_{n\to + \infty} \mathbb{E}\left[f(x_k) - f_\star \right] 
\leq 0.
\end{align*}

Theorem \ref{theorem:1}, together with \eqref{conv} and \eqref{step},
leads to \eqref{main:4} and \eqref{main:4_1} 
with rate of convergence \eqref{main1}, \eqref{main2}, and \eqref{triangle_1}.
This completes the proof.

{\em Proofs of Theorems \ref{cor:2_1} and \ref{cor:2_2}:}
A discussion similar to the one for showing Theorem \ref{cor:1}
(resp. Theorem \ref{cor:2}) with $\mathsf{g} := \mathrm{grad}f$ 
leads to Theorem \ref{cor:2_1} (resp. Theorem \ref{cor:2_2}).

\section*{Proof of Corollary \ref{COR:1_1}}

{\em Proof:}
The step-sizes $\alpha_n := 1/n^\eta$ ($\eta \in (1/2,1], n \geq 1$) and  
$(\beta_n)_{n\in\mathbb{N}}$ such that $\sum_{n=1}^{+\infty} \alpha_n \beta_n < + \infty$ satisfy \eqref{step:0}.
Hence, Theorem \ref{cor:2} leads to \eqref{main:5}.

Let $\alpha_n := 1/n^\eta$ ($\eta \in [1/2,1), n \geq 1$) and  
$(\beta_n)_{n\in\mathbb{N}}$ such that $\sum_{n=1}^{+\infty} \beta_n < + \infty$. 
We have that 
\begin{align}\label{CS2}
\lim_{n \to + \infty} \frac{1}{n\alpha_n} 
= \lim_{n \to + \infty} \frac{1}{n^{1-\eta}} = 0.
\end{align}
Moreover,  
\begin{align}\label{CS1}
\begin{split}
\frac{1}{n} \sum_{k=1}^n \alpha_k^2
&\leq
\frac{1}{n} \sum_{k=1}^n \alpha_k
\leq 
\frac{1}{n} \left\{ 1 + \int_1^n \frac{\mathrm{d}t}{t^\eta}   \right\}
= 
\frac{1}{n} \left\{ \frac{n^{1-\eta}}{1-\eta} - \frac{\eta}{1-\eta}   \right\}\\
&\leq
\frac{1}{1-\eta}\frac{1}{n^\eta}.
%\leq 
%\frac{1}{1-\eta}\frac{1}{n^{1-\eta}}.
\end{split}
\end{align}
Hence,
$\lim_{n\to +\infty} (1/n) \sum_{k=1}^n \alpha_k 
= \lim_{n\to +\infty} (1/n) \sum_{k=1}^n \alpha_k^2 = 0$.
From $\sum_{n=1}^{+\infty} \beta_n < + \infty$, 
$\lim_{n\to +\infty} (1/n) \sum_{k=1}^n \beta_k = 0$.
Hence, 
$\alpha_n := 1/n^\eta$ and 
$(\beta_n)_{n\in\mathbb{N}}$ such that 
$\sum_{n=1}^{+\infty} \beta_n < + \infty$ satisfy 
\eqref{step}.
Accordingly, from Theorem \ref{cor:2} with \eqref{conv}, \eqref{CS2}, and \eqref{CS1},
we have the assertions in Corollary \ref{COR:1_1}.
%In the case where $\alpha_n = 1/n^{\eta}$ ($\eta \in (1/2,1), n \geq 1$), we have that
%$\lim_{n\to + \infty} 1/(n\alpha_n) = \lim_{n\to +\infty} 1/n^{1-\eta} = 0$.
%The Cauchy-Schwarz inequality ensures that 
%\begin{align}\label{CS}
%\frac{1}{n} \sum_{k=1}^n \alpha_k^2
%\leq
%\frac{1}{n} \sum_{k=1}^n \alpha_k
%\leq 
%\frac{1}{n} \sqrt{\sum_{k=1}^n 1^2} 
%\sqrt{\sum_{k=1}^n \left(\frac{1}{k^{\eta}}\right)^2 }
%\leq \frac{\bar{B}}{\sqrt{n}},
%\end{align}
%where $\bar{B}^2 := \sum_{n=1}^{+\infty} \alpha_n^2 < + \infty$.
%Hence, $\lim_{n\to +\infty} (1/n) \sum_{k=1}^n \alpha_k 
%= \lim_{n\to +\infty} (1/n) \sum_{k=1}^n \alpha_k^2 = 0$.
%The condition $\sum_{n=1}^{+\infty} \beta_n < + \infty$
%implies that 
%$\lim_{n\to +\infty} (1/n) \sum_{k=1}^n \beta_k = 0$.
%Accordingly, 
%$\lambda_n := 1/n^{\eta}$ ($\eta \in (1/2, 1)$) and 
%$(\beta_n)_{n\in\mathbb{N}}$ such that 
%$\sum_{n=1}^{+\infty} \beta_n < + \infty$ satisfy 
%\eqref{step}.
%Theorem \ref{cor:2} with \eqref{CS} and \eqref{conv} thus implies the assertions in Corollary \ref{COR:1_1}.
%An argument similar to the one for obtaining Corollary \ref{COR:1_1},
%together with \eqref{triangle_1} and \eqref{triangle}, leads to Corollary \ref{COR:1_2}.

\end{document}